\documentclass[a4paper,12pt]{amsart}
\usepackage[utf8]{inputenc}
\usepackage[english]{babel}
\usepackage{amsmath,amsthm,amssymb,amsfonts,mathtools}
\usepackage{tikz}
\usetikzlibrary{calc,decorations.pathreplacing}
\usepackage{caption}
\usepackage{float}
\usepackage[margin=0.8in]{geometry}
\usepackage[colorlinks=true]{hyperref}
\usepackage[nameinlink,noabbrev]{cleveref}
\hypersetup{linkcolor=blue,citecolor=red,urlcolor=blue}
\usepackage{enumitem}
\usepackage{aliascnt}

\numberwithin{equation}{section}

\newcommand{\C}{\mathcal{C}}
\newcommand{\F}{\mathcal{F}}
\newcommand{\R}{\mathbb{R}}
\newcommand{\Z}{\mathbb{Z}}
\newcommand{\N}{\mathbb{N}}

\theoremstyle{definition}
\newtheorem{definition}{Definition}[section]

\theoremstyle{plain}
\newaliascnt{theorem}{definition}
\newtheorem{theorem}[theorem]{Theorem}
\aliascntresetthe{theorem}
\newaliascnt{proposition}{definition}
\newtheorem{proposition}[proposition]{Proposition}
\aliascntresetthe{proposition}
\newaliascnt{lemma}{definition}
\newtheorem{lemma}[lemma]{Lemma}
\aliascntresetthe{lemma}
\newaliascnt{corollary}{definition}
\newtheorem{corollary}[corollary]{Corollary}
\aliascntresetthe{corollary}
\newaliascnt{example}{definition}
\newtheorem{example}[example]{Example}
\aliascntresetthe{example}

\theoremstyle{remark}
\newaliascnt{remark}{definition}
\newtheorem{remark}[remark]{Remark}
\aliascntresetthe{remark}

\crefname{definition}{definition}{definitions}
\Crefname{definition}{Definition}{Definitions}
\crefname{theorem}{theorem}{theorems}
\Crefname{theorem}{Theorem}{Theorems}
\crefname{proposition}{proposition}{propositions}
\Crefname{proposition}{Proposition}{Propositions}
\crefname{lemma}{lemma}{lemmas}
\Crefname{lemma}{Lemma}{Lemmas}
\crefname{corollary}{corollary}{corollaries}
\Crefname{corollary}{Corollary}{Corollaries}
\crefname{remark}{remark}{remarks}
\Crefname{remark}{Remark}{Remarks}

\renewcommand{\le}{\leqslant}
\renewcommand{\leq}{\leqslant}
\renewcommand{\ge}{\geqslant}

\title[Geometry of ground states in the long-range Ising model]
{Geometry of ground states \\ in the long-range Ising model}

\date{\today}

\begin{document}

\author{Serena Dipierro}
\author{Riccardo Sabbatini}
\author{Enrico Valdinoci}

\address{Serena Dipierro: Department of Mathematics and Statistics, The University of Western Australia, 35 Stirling Highway, Crawley, Perth, WA 6009, Australia}
\email{serena.dipierro@uwa.edu.au}

\address{Riccardo Sabbatini: Department of Mathematics and Statistics, The University of Western Australia, 35 Stirling Highway, Crawley, Perth, WA 6009, Australia}
\email{riccardo.sabbatini@research.uwa.edu.au}

\address{Enrico Valdinoci: Department of Mathematics and Statistics, The University of Western Australia, 35 Stirling Highway, Crawley, Perth, WA 6009, Australia}
\email{enrico.valdinoci@uwa.edu.au}

\begin{abstract}
We study the ground states of the ferromagnetic Ising model on $\mathbb Z^d$ with
long-range interactions and prescribed exterior data. We first show that the half-space
configuration is a ground state, by means of a discrete calibration (which is the lattice
counterpart of the calibrations of Cabr\'e and Pagliari for nonlocal minimal surfaces).
In dimension $d\ge2$ and for positive kernels, the half-space configuration is the unique minimizer in every
finite domain, and it is ``stable'' with respect to perturbations of the exterior datum whose energy
is below half of a positive excitation gap. This rigidity has no counterpart in the continuum.

For general exterior data and positive kernels, we show that along a monotone family the minimizer is unique except
for finitely many data. This can be considered
a ``generic uniqueness result'' that
provides the discrete counterpart of a recent work of Dipierro, Fern\'andez-Real and Valdinoci in the continuum.

We then introduce a pointwise notion, called perception, which is the lattice counterpart of
the instability and stickiness phenomena for nonlocal minimal surfaces discovered by
Dipierro, Savin and Valdinoci. We give a local criterion for the half-space configuration
to be perceptive at a given site. 

Finally,
we consider a box with exterior datum $+1$ above and below a slab and~$-1$ on its sides.
For fractional kernels with small exponent, the minimizer is~$+1$ in the whole box, and
its positive phase is connected. For short-range kernels, instead, the negative phase
may enter the box, and the positive phase may be disconnected. Moreover, if the datum is~$-1$ around the box and~$+1$ only at distance two from it, the configuration identically~$-1$ is perceptive at every boundary site of the box for fractional kernels with small exponent.

\end{abstract}
\maketitle 

\setcounter{tocdepth}{1}
\tableofcontents

\section{Introduction}
\subsection{Physical description of the ferromagnetic Ising model}
In this paper, we study the geometric features of
ground states of a discrete statistical mechanics
spin model describing phase transitions and ferromagnetism.

At the end of the nineteenth century, the physics community was particularly interested
in understanding how certain materials change their physical properties as the
temperature and the external magnetic field vary. In this context, in 1895 Pierre Curie
observed that some materials are magnetized at low temperature, and lose their
magnetization above a critical temperature, now called the \emph{Curie
temperature}.\footnote{Below this temperature the material is in an ordered phase, called
\emph{ferromagnetic}, and above it in a disordered one, called \emph{paramagnetic}.}

At the beginning of the twentieth century, this type of phase transition (from a
magnetized phase to a non-magnetized one) was described by the phenomenological theory
of Weiss~\cite{Weiss1907}. A
rigorous proof of the existence of a phase transition, starting from a microscopic model
containing no a priori information about the transition itself, was therefore still
missing. One of the first models of statistical mechanics for which a phase
transition could be proved is the Ising model.

The Ising model was introduced by Lenz in 1920~\cite{Lenz1920}, 
and studied in 1925 by Lenz's student Ising~\cite{Ising1925}
in dimension one 
and in
1936 by Peierls~\cite{Peierls1936} in higher dimensions.
The model deals with a lattice in which
each atom carries a magnetic moment that tends to align with those of the neighbouring
atoms.

It was soon understood that phase transitions in the Ising model
were affected both by dimensionality and interaction range.
Indeed, while Ising showed that no phase transition
occurs in the one-dimensional case, the work of Peierls
detected phase transitions in higher dimensions.
Also, while the original system focused on nearest-neighbor interactions, the model easily allows for generalizations
to longer-range interactions: in this spirit,
Kac and Thompson~\cite{KacThompson1969} conjectured that, if the interaction between the~$i$-th and~$j$-th sites of the lattice
decays like~$|i-j|^{-\alpha}$ with $1<\alpha\le 2$, the one-dimensional Ising model undergoes a phase
transition. Dyson~\cite{Dyson1969} proved this conjecture for~$\alpha\in(1,2)$,
and Fr\"ohlich and Spencer~\cite{MR660541} settled the case~$\alpha=2$.

Although the Ising model was introduced to describe how a material behaves as the
external temperature varies, the case where the temperature is zero is of interest in
itself.\footnote{Here the temperature $T$ is measured in kelvin, so that $T=0$
corresponds to absolute zero.} At zero temperature, thermal fluctuations vanish, and the
study of the Ising model reduces to that of the minimizers of its Hamiltonian, called
\emph{ground states}.

More precisely, each site $i \in \Z^d$ carries a magnetic moment, called a
\emph{spin}, taking the values~$-1$ and $+1$. A state of the system, called a
\emph{configuration}, is a function $u\colon\Z^d\to\{-1,+1\}$, and we write
$u_i:=u(i)$ for the spin at the site $i$. Two sites $i,j\in\Z^d$ are
\emph{nearest neighbours}, and we write $i\sim j$, if $|i-j|=1$. In the
classical model, spins interact only with their nearest neighbours. Given a
finite set $\Gamma\subset\Z^d$ and an external magnetic field $B\in\R$, the
ferromagnetic Hamiltonian in $\Gamma$ is
\begin{equation*}
H_\Gamma(u)=-\sum_{\substack{\{i,j\}\subseteq\Gamma\\ i\sim j}}u_iu_j
\;-\;B\sum_{i\in\Gamma}u_i,
\end{equation*}
the first sum being taken over the unordered pairs of nearest neighbours
contained in $\Gamma$. The first term of~$H_\Gamma$ favours the alignment of
neighbouring spins, whence the name \emph{ferromagnetic},\footnote{The
terminology refers to the sign of the interaction, and does not by itself entail
that the model behaves as a ferromagnet.} while the second favours the alignment
of each spin with the field $B$. We refer to~\cite{FriedliVelenik} for a
comprehensive introduction to the Ising model.
\subsection{A long-range Ising model and goals of this paper}
A natural extension of the classical Ising model
consists in
letting sites at any distance interact.
Already in dimension one, the range of the interaction changes the behaviour of the
model.
From the variational point of view, the Ising energy can be regarded as a discrete
approximation of the perimeter when the interactions are short range~\cite{MR2123651},
and of the fractional perimeter when they are long range~\cite{MR3652519}. To see this
parallel, one associates with every configuration an interface separating the region
where the spins equal $+1$ from the one where they equal $-1$.

The goal of this paper is to understand the geometry of the ground states of the
long-range Ising model under prescribed exterior data.

We first prove in~\Cref{th:rigidity} that the half-space configuration, defined in~\eqref{eq:flatconfig}, is a ground state.
The proof relies on a calibration, which is the discrete analogue of the calibrations for
nonlocal minimal surfaces introduced by Cabré
in~\cite{Cabr__2020} and by Pagliari in~\cite{Pagliari_2019}. 
We also show that the half-space configuration is the unique
minimizer in every finite domain, and that it remains so under sufficiently small
perturbations of the exterior datum. Uniqueness also follows from the calibration, as
in~\cite{Pagliari_2019}.

For a general exterior datum, instead, uniqueness may fail, see~\Cref{ex:nonunique}.
However, it fails only for exceptional data. Along a monotone family of exterior data,
the minimizer is unique except for finitely many data, see~\Cref{prop:generic}. This
result relies on the maximum principle in~\Cref{cor:max}, and it is the lattice
counterpart of a recent result establishing the ``generic uniqueness'' of nonlocal minimal surfaces~\cite{DFV2026} (see also~\cite{MR5036406, MR4649881, BERNO} for related results for free boundary problems).

A further question is how the interface of a ground state responds to a perturbation of
the exterior datum. In the continuum, this question is the subject of the stickiness
theory for nonlocal minimal surfaces, initiated in~\cite{DIPIERRO20171791}. There, flat
nonlocal minimal surfaces are called \textit{unstable}, because an arbitrarily small perturbation
of the exterior datum can change the minimizer. For instance, in the plane, a small
compactly supported perturbation of a half-plane datum makes the minimizer detach from
the half-plane, see~\cite[Theorem~1.4]{DIPIERRO20171791}. Since the perturbation is
placed away from the domain, the minimizer then does not attain the exterior datum
continuously at the boundary, but \textit{sticks} at it. 

In the discrete setting,
we introduce a pointwise notion,
called \textit{perception}, which describes whether a
perturbation of the exterior datum reverses the spin at a given site in every minimizer,
see~\Cref{def:sticky}. Perception at the boundary is the case in which the site lies on the
boundary of the domain and the perturbation lies at distance greater than one from it,
see~\Cref{def:pointsticky}. In~\Cref{th:stickiness} we give a local criterion for the
half-space configuration to be perceptive at a site, and in~\Cref{cor:stickiness} we deduce
perception at the boundary.
This is the discrete analogue of stickiness.

Finally, motivated by~\cite{dipierro2020disconnectednessnonlocalminimalsurfaces}, we
study a problem in which the spins may change in a vertical cylinder. Outside the
cylinder, we prescribe the exterior datum~$-1$ in a horizontal slab and~$+1$ above and
below it, see~\Cref{fig:regions}. We then consider the minimizers of the Ising model in the
cylinder, and we ask whether their positive phase joins the two sides of the slab through
the cylinder. As we show in~\Cref{th:fumihiko}, the answer to this connectedness question
depends on the range of the interaction and on the geometry.

For simplicity, we prove~\Cref{th:fumihiko} when the spins may change only in the box
obtained by intersecting the cylinder with the slab. Outside this box the spins are given,
and they are equal to~$+1$ in the parts of the cylinder above and below the slab. The case
of the cylinder reduces to that of the box, and we postpone this technical step to
Appendix~\ref{sec:reduction}.

We also consider the box with exterior datum identically $-1$, whose minimizer is identically~$-1$. The datum of~\Cref{th:fumihiko} is obtained by reversing the spins in the two
half-spaces above and below the box. Hence, for fractional kernels with small exponent, the configuration identically~$-1$ is
perceptive at every site of the box. However, it is not boundary-perceptive at any site,
since the perturbation touches the box.

If, instead, a layer of width one around the box is left equal to~$-1$, then for
fractional kernels with small exponent the configuration identically $-1$ is
boundary-perceptive at every boundary site of the box, see~\Cref{prop:collar}. This is the lattice counterpart of the complete adherence observed
in~\cite{dipierro2020disconnectednessnonlocalminimalsurfaces}.

\subsection{Notation}

We say that a configuration of the Ising model  is a function $u: \Z^d \to \{-1,1\}$. The set of all configurations is denoted by $\C$. Given two configurations $u,v \in \C$, we define
$$ (u\vee v)_i:=\max\{u_i,v_i\}\qquad {\mbox{and}}\qquad(u\wedge v)_i:=\min\{u_i,v_i\}.$$

Throughout the paper, $|\cdot|$ denotes the Euclidean norm and $\|\cdot\|_{\infty}$ the sup norm. We use the convention~$\operatorname{sign}(0):=0$.

We let $\mathbb Z^d$ carry the nearest-neighbour graph structure: $i\sim j$
if and only if $|i-j|=1$. A set~$A\subseteq\mathbb Z^d$ is \emph{connected} if any
two of its sites are joined by a path in $A$ with steps of this form.

The \emph{positive phase} of a configuration~$u$ is
$S(u):=\{i\in\mathbb Z^d : u_i=+1\}$. The \emph{interface} of~$u$ is the boundary in $\mathbb R^d$ of the union\footnote{The interface is the
union of the faces shared by the cubes centred at $i$ and $j$, for $i\sim j$ with
$u_i\neq u_j$. In~$\mathbb Z^d$, one can define it as the set of pairs $\{i,j\}$ with
$i\sim j$ and $u_i\neq u_j$, or as the set of sites of~$S(u)$ having a nearest neighbour
outside $S(u)$. The first notion of interface in $\Z^d$ is invariant under the flip $u\mapsto-u$, while the
second one is not.} of the unit
cubes centred at the sites of~$S(u)$, namely of
$$\bigcup_{i\in S(u)}\left(i+[-\tfrac12,\tfrac12]^d\right).$$

\medskip

Given a finite set $\Gamma \subset \Z^d$, we define the Hamiltonian as the formal sum
\begin{equation}\label{eq:hamiltonian}
H_{\Gamma}(u;K):= \frac{1}{2} \sum_{i \in \Gamma,\, j \in \Gamma} K(i,j) (1-u_iu_j) + \sum_{i \in \Gamma , j \in \Z^d \setminus \Gamma} K(i,j) (1-u_iu_j), 
\end{equation}
where  $K: \Z^d \times \Z^d \to \R$ is the interaction kernel. For the classical Ising model with nearest-neighbour interactions, the kernel is
\begin{equation}\label{eq:NNkernel}
K_{\rm nn}(i,j):=
\begin{cases}
1, & {\mbox{ if }}|i-j|=1,\\
0, & {\mbox{ otherwise.}}
\end{cases}
\end{equation}
The assumptions on the kernels considered in this paper are listed in~\Cref{sec:assumptions}.

In what follows, we will omit denoting the kernel in the Hamiltonian when it is clear from the context, thus we write $H_\Gamma(u)$ instead of~$H_\Gamma(u;K)$.

Given a finite set $\Gamma \subset \Z^d$, a summable kernel $K$, a configuration~$u$ on $\Gamma$, and a boundary datum $\eta \in \mathcal C$ on $\Z^d \setminus \Gamma$, we write 
\begin{equation*}
   H_\Gamma(u \mid \eta):= \frac{1}{2} \sum_{i , j \in \Gamma} K(i,j) (1-u_iu_j) + \sum_{i \in \Gamma, \, j \in \Z^d \setminus \Gamma} K(i,j) (1-u_i\eta_j)
\end{equation*}
for the Hamiltonian evaluated on the configuration equal to~$u$ in $\Gamma$ and to~$\eta$ outside.

\begin{definition}\label{def:minimizer}
We say that a configuration~$u$ is a \emph{minimizer} for $H$ in a finite set $\Gamma \subseteq \Z^d$ if
$$
H_\Gamma(u) \le H_\Gamma(v),
$$
for every configuration~$v$ that agrees with $u$ outside of~$\Gamma$.
\end{definition}

\begin{definition}\label{def:ground}
We say that a configuration~$u$ is a \emph{ground state} for $H$ if it is a minimizer for~$H$ in every finite set $\Gamma \subset \Z^d$.
\end{definition}

\subsection{The assumptions of this paper}\label{sec:assumptions}
Throughout the paper, we assume that the interaction kernel $K$ satisfies the following properties:
\begin{enumerate}[label=(\roman*),ref=(\roman*)]
    \item\label{it:transinv} $K(i,j)=K(j-i)$, namely, the interaction is translation invariant;
    \item\label{it:even} $K$ is even, that is, $K(i)=K(-i)$;
    \item\label{it:noself} $K(0)=0$, namely, each site does not interact with itself;
    \item\label{it:ferro} $K$ is nonnegative, that is, $K(i) \ge 0$ for every $i \in \mathbb Z^d$,
    namely, the interaction is ferromagnetic;
    \item\label{it:summable} $K$ is summable, that is,
    \[
    \sum_{i\in\mathbb Z^d} K(i)<+\infty.
    \]
\end{enumerate}

In some results, we further assume that $K$ is a fractional kernel, namely we use the notation
\begin{equation}\label{eq:frachyp}
    K_s(i,j) :=\frac{1}{|i-j|^{d+2s}},
    \qquad \text{for} \, s\in(0,1).
\end{equation}

We point out that
the notion of minimizer was introduced 
in~\Cref{def:minimizer}
for finite sets $\Gamma$. When $\Gamma$ is
infinite, the Hamiltonian~$H_\Gamma$ is in general not finite. Nevertheless, the
difference of the energies of two configurations is well defined whenever they differ
at finitely many sites. More precisely, let $u$ and $v$ be such that the set
$E:=\{i\in\Z^d\ \text{such that}\ u_i\neq v_i\}$ is finite and contained in $\Gamma$. We set
\[
H_\Gamma(u)-H_\Gamma(v)
:=\sum_{ l\in E,\, j\in\mathbb Z^d\setminus E}
K(l,j)\,\big[(1-u_lu_j)-(1-v_lv_j)\big].
\]
The sum is absolutely convergent by assumption~\ref{it:summable}. This allows us to
extend the notion of minimizer to infinite sets.

\begin{definition}\label{rmk:gammainfinite}
Let $\Gamma\subseteq\Z^d$. We say that $u$ is a \emph{minimizer in} $\Gamma$ if
$H_\Gamma(u)-H_\Gamma(v)\le0$ for every $v$ coinciding with $u$ outside a finite subset
of~$\Gamma$.
\end{definition}
When $\Gamma$ is finite, \Cref{rmk:gammainfinite} agrees with~\Cref{def:minimizer}.
In particular, a ground state is a minimizer in every $\Gamma\subseteq\Z^d$. Indeed, the
competitors allowed in $\Gamma$ differ from it at finitely many sites, and a ground state
has lower energy than all such configurations. Conversely, the minimizers in $\Z^d$ are the ground states.
\subsection{Perception under perturbations of the exterior datum}
We introduce notions of \textit{perception} and of \textit{perception at the boundary} for the minimizers
of the Ising model with respect to perturbations of the exterior datum.

Let $\Gamma\subset\Z^d$ be finite, let $\eta$ be an exterior datum, and let
$\Lambda\subseteq\Z^d\setminus\Gamma$. Set\footnote{We denote by $\chi_A$ the
characteristic function of a set $A\subseteq\Z^d$, that is, $\chi_A(i)=1$ if $i\in A$ and
$\chi_A(i)=0$ otherwise.}
\begin{equation}\label{ndjkwhriuy32234567SDFGHVdye4ut}
\eta_\Lambda:=\eta\,(1-2\chi_\Lambda),
\end{equation}
that is, the datum obtained from $\eta$ by reversing the spins of~$\Lambda$.

\begin{definition}[Perception]\label{def:sticky}
Let $u$ be a minimizer of~$H_\Gamma(\cdot\mid\eta)$. We say that $u$ is \emph{perceptive at}
$i\in\Gamma$ under the perturbation $\Lambda$ if $v_i=-u_i$ for every minimizer $v$ of
$H_\Gamma(\cdot\mid\eta_\Lambda)$.
\end{definition}

In other words, $u$ is perceptive at $i$ if the perturbation $\Lambda$, although located
outside $\Gamma$, changes the spin at $i$ of every minimizer.

To visualize the notion of
perception, see~\Cref{fig:collar}.
More precisely, in~\Cref{fig:collar}(a), the domain is~$\Gamma=\{-2,\ldots,2\}\times\{-1,0,1\}$ and the initial configuration is~$u\equiv-1$.
In~\Cref{fig:collar}(b), the exterior spins are reversed on~$\Lambda=\{j\in\Z^2:|j_2|\ge3\}$.
For the fractional kernel~\eqref{eq:frachyp} with~$s$ sufficiently small,
\Cref{prop:collar} shows that every minimizer for the perturbed datum equals~$+1$
in~$\Gamma$. In particular, at the 
central site~$k=(0,0)$, every such
minimizer has spin~$+1=-u_k$
and accordingly~$u$ is perceptive at~$k$.
\medskip

Now, given a configuration~$u$, we set
\begin{equation}\label{3456cmeuiwdef59ui54bordo}
\partial_u\Gamma:=\big\{i\in\Gamma\ \text{such that}\ u_i=u_{i'}\ \text{for some}\
i'\in\Z^d\setminus\Gamma\ \text{with}\ |i-i'|=1\big\}.
\end{equation}
A site $i\in\Gamma$ belongs to~$\partial_u\Gamma$ if it has a nearest neighbour
$i'\in\Z^d\setminus\Gamma$ such that $u$ takes the same value at $i$ and at $i'$.

\begin{figure}[H]
\centering
\definecolor{pmmark}{RGB}{230,0,0}
\newcommand{\pmplus}[2]{\draw[fill=white,line width=0.4pt] (#1,#2) circle (0.2);
  \draw[line width=0.5pt] (#1-0.11,#2)--(#1+0.11,#2) (#1,#2-0.11)--(#1,#2+0.11);}
\newcommand{\pmminus}[2]{\draw[fill=white,line width=0.4pt] (#1,#2) circle (0.2);
  \draw[line width=0.5pt] (#1-0.11,#2)--(#1+0.11,#2);}
\newcommand{\pmlabel}[5]{\draw[pmmark,line width=1.1pt] (#1,#2) circle (0.25);
  \node[font=\scriptsize,text=pmmark,inner sep=0.6pt] at (#1+#4,#2+#5) {$#3$};}
\newcommand{\pmpanel}[4]{%
\begin{tikzpicture}[x=0.55cm,y=0.55cm]
\useasboundingbox (-3.5,-2.5) rectangle (6.5,5.2);
\fill[gray!22] (-3.5,-2.5) rectangle (0.5,4.5);
#3
\draw[gray!55,thin,step=1] (-3.5,-2.5) grid (6.5,4.5);
\draw[line width=1.3pt,line cap=round] #2;
\foreach \x in {-3,...,6}{\foreach \y in {-2,...,4}{\pmminus{\x}{\y}}}
\foreach \x/\y in {#1}{\pmplus{\x}{\y}}
\node[font=\footnotesize] at (-1.5,4.95) {$\Gamma$};
#4
\end{tikzpicture}}
\begin{minipage}{0.49\textwidth}\centering
\pmpanel{99/99}{}{}%
{\pmlabel{0}{1}{i}{-0.5}{0.5}\pmlabel{1}{1}{i'}{0.5}{0.5}}
\par\smallskip{\footnotesize (a)}
\end{minipage}\hfill
\begin{minipage}{0.49\textwidth}\centering
\pmpanel{1/0,1/1,1/2,2/0,2/1,2/2,3/0,3/1,3/2,4/0,4/1,4/2,5/0,5/1,5/2,6/0,6/1,6/2}%
{(6.5,-0.5)--(0.5,-0.5)--(0.5,2.5)--(6.5,2.5)}%
{\fill[orange!25] (0.5,-0.5) rectangle (6.5,2.5);}%
{\node[font=\footnotesize] at (3.5,4.95) {$\Lambda$};
\pmlabel{0}{1}{i}{-0.5}{0.5}\pmlabel{1}{1}{i'}{0.5}{0.5}}
\par\smallskip{\footnotesize (b)}
\end{minipage}

\caption{\footnotesize The spins $+1$ and $-1$ are drawn as $\oplus$ and $\ominus$. The sites of~$\Gamma$ lie in the grey region, and those of~$\Lambda$ in the orange one, bounded by a thick line. The sites named in the text are circled in red.}
\label{fig:instability}
\end{figure}
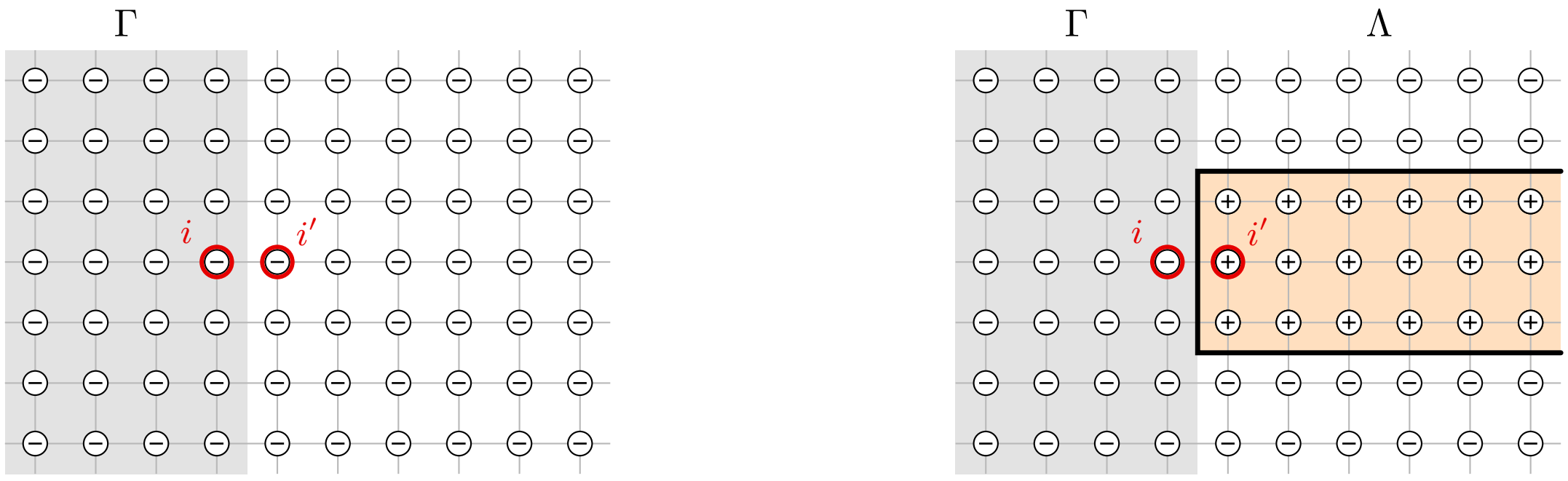

In~\Cref{fig:instability}(a), we have $u_i=u_{i'}$, so that $i\in\partial_u\Gamma$. In~\Cref{fig:instability}(b), the configuration~$u$ is equal to~$+1$ on $\Lambda$. Hence $u_{i'}\neq u_i$, and $i$ does not belong to~$\partial_u\Gamma$.

\begin{definition}[Perception at the boundary]\label{def:pointsticky}
Assume that $\operatorname{dist}(\Lambda,\Gamma)\ge2$. Let $u$ be a minimizer
of~$H_\Gamma(\cdot\mid\eta)$ and fix $i\in\partial_u\Gamma$. We say that $u$
\emph{is boundary-perceptive at} $i$ under the perturbation $\Lambda$ if $u$ is perceptive at $i$.
\end{definition}

Roughly speaking, $u$ is boundary-perceptive at $i$ if $u$ is perceptive at $i$ and $i$ is aligned with a nearest neighbour $i'\in\Z^d\setminus\Gamma$ that the perturbation does not reach.

To visualize an example
of perception at the boundary, take $i=(2,0)\in\Gamma$
and $i'=(3,0)\notin\Gamma$, as in~\Cref{fig:collar}.
Since $|i-i'|=1$ and $u_i=u_{i'}=-1$, we have $i\in\partial_u\Gamma$.
For $x\in\Gamma$ and $j\in\Lambda$, we have $|x_2|\le1$ and $|j_2|\ge3$,
so $|j-x|\ge|j_2-x_2|\ge2$. Equality holds for $x=(0,1)$ and $j=(0,3)$.
Thus $\operatorname{dist}(\Lambda,\Gamma)=2$. By~\Cref{prop:collar}, every
perturbed minimizer has spin $+1$ at $i$. Hence, in this case, $u$ is boundary-perceptive at $i$.\medskip

We observe that perception requires an interaction between $\Lambda$ and $\Gamma$.
Indeed, if $K(l,j)=0$ for every $l\in\Gamma$ and $j\in\Lambda$, then
$H_\Gamma(\cdot\mid\eta_\Lambda)=H_\Gamma(\cdot\mid\eta)$. Hence $u$ is also a minimizer of
the perturbed problem, and it is not perceptive at any site. Moreover, for the
nearest-neighbour kernel, the requirement $\operatorname{dist}(\Lambda,\Gamma)\ge2$
in~\Cref{def:pointsticky} implies that no site of~$\Lambda$ interacts with a site of
$\Gamma$. Hence no site is boundary-perceptive for the nearest-neighbour kernel.

\subsection{Statement of the main results}
Let
\[
\Pi^+:=\{i=(i_1,\dots,i_d)\in\mathbb Z^d : i_d> 0\}
\qquad{\mbox{and}}\qquad
\Pi^-:=\mathbb Z^d\setminus\Pi^+,
\]
and define the half-space configuration~$\omega$ as
\begin{equation}\label{eq:flatconfig}
\omega_i:=
\begin{cases}
-1, & {\mbox{ if }}i\in\Pi^+,\\
+1, &{\mbox{ if }} i\in\Pi^-.
\end{cases}
\end{equation}

Given a finite set $\Gamma \subset\mathbb Z^d$, we consider the problem
\begin{equation}\label{eq:flatproblem}
\min H_\Gamma(\cdot \mid \omega).
\end{equation}
Given a set $\Lambda \subset\Pi^+\setminus\Gamma$ and $\omega_\Lambda:=\omega+2\chi_\Lambda$, 
we consider the problem
\begin{equation}\label{eq:perturbedproblem}
\min H_\Gamma(\cdot \mid \omega_\Lambda).
\end{equation}
We refer to~\eqref{eq:flatproblem} as the flat problem and to~\eqref{eq:perturbedproblem} as the perturbed problem.

Finally, we measure the energy of the perturbation by
\begin{equation}\label{eq:alpha}
\alpha_\Gamma(\Lambda):=2\sum_{i\in\Gamma,\,j\in\Lambda}K(i,j).
\end{equation}
This is the largest change of energy that the perturbation produces on a configuration,
see the forthcoming formula~\eqref{eq:datumchange}.

\begin{definition}[Excitation gap]\label{def:gap}
Let $\Gamma\subset\mathbb Z^d$ be finite and nonempty. The \emph{excitation gap} of~$\omega$ in $\Gamma$ is
\begin{equation}\label{eq:gap}
\delta_\Gamma:=\min\big\{H_\Gamma(v\mid\omega)-H_\Gamma(\omega\mid\omega)\  :
v\neq\omega\ \text{in}\ \Gamma\big\}.
\end{equation}
\end{definition}

We point out that
the minimum in~\eqref{eq:gap} is taken over finitely many configurations, since~$\Gamma$ is finite. Thus~$\delta_\Gamma$ is well defined, and it is the smallest energy needed to excite the flat
configuration in~$\Gamma$.

Our first result establishes the rigidity of the half-space configuration. We prove that it is the unique solution of the half-space boundary value problem~\eqref{eq:flatproblem} and that it remains unchanged inside the domain under sufficiently small perturbations of the boundary conditions.

\begin{theorem}\label{th:rigidity}
The half-space configuration~$\omega$ is a ground state, and it achieves the minimum in~\eqref{eq:flatproblem} for every finite $\Gamma\subset\mathbb Z^d$.

Assume in addition that $d\ge2$ and $K(i,j)>0$ for every $i\neq j$. Then $\omega$ is the \emph{unique}
minimizer of~\eqref{eq:flatproblem} for every finite $\Gamma$.
Moreover $\delta_\Gamma>0$ and, if
\[
2\,\alpha_\Gamma(\Lambda)<\delta_\Gamma,
\]
then $\sigma=\omega$ in $\Gamma$ for every minimizer $\sigma$ of~\eqref{eq:perturbedproblem}.
\end{theorem}

The excitation gap separates rigidity from perception. Indeed, by~\Cref{th:rigidity}, the flat interface is rigid under perturbations with $2\,\alpha_\Gamma(\Lambda)<\delta_\Gamma$, and
perception can occur only when~$2\,\alpha_\Gamma(\Lambda)\ge\delta_\Gamma$. For a single site, this threshold separates the two regimes exactly, see~\Cref{rmk:sharpsite}.

Uniqueness for the half-space datum is a special feature. For a general exterior datum
it may fail, see~\Cref{ex:nonunique}. On the other hand, the next result shows that, along a monotone
family of exterior data, this happens only for finitely many data. This is the lattice
counterpart of~\cite[Theorem~1.1]{DFV2026}.
In what follows, for a finite set~$X$, we denote by~$\sharp X$ its cardinality.

\begin{theorem}[Generic uniqueness]\label{prop:generic}
Assume that $K(i,j)>0$ for every $i\neq j$, and let $\Gamma\subset\mathbb Z^d$ be finite.

Let $I\subseteq\mathbb R$,
and let $\{\eta^{(t)}\}_{t\in I}$ be exterior data such that
$\eta^{(t_1)}\le\eta^{(t_2)}$ and $\eta^{(t_1)}\neq\eta^{(t_2)}$ for every~$t_1<t_2$ in~$I$. Then $H_\Gamma(\cdot\mid\eta^{(t)})$ has a unique minimizer for every~$t\in I$,
except for at most~$\sharp\Gamma$ values of~$t$.
\end{theorem}

We stress that the index set~$I$ is not required to be discrete. An example of a
family of exterior data indexed by the whole of~$\mathbb R$ is given
in~\Cref{ex:irrational}.

Our third result is a criterion for perception. If the perturbation is strong
enough near a site of~$\Gamma$, then every minimizer of the perturbed problem differs
from $\omega$ at that site.

\begin{theorem}\label{th:stickiness}
Let $\sigma$ be any minimizer of~\eqref{eq:perturbedproblem}. Assume that there exists $i\in\Gamma$ such that
\begin{equation}\label{eq:instabilitycondition}
\omega_i=-1
\qquad\text{and}\qquad
2 \sum_{j\in\Lambda}K(i,j)
>
-\sum_{j \neq i}K(i,j)\,\omega_j.
\end{equation}
Then $\sigma_i=+1$. In particular, $\sigma\neq\omega$ in $\Gamma$, and $\omega$ is perceptive at $i$ under the
perturbation $\Lambda$.
\end{theorem}

\begin{corollary}\label{cor:stickiness}
Assume that $\operatorname{dist}(\Lambda,\Gamma)\ge2$, and that
$i\in\partial_\omega\Gamma$ satisfies~\eqref{eq:instabilitycondition}. Then $\omega$ is boundary-perceptive
at~$i$ under the perturbation $\Lambda$.
\end{corollary}

\begin{figure}[H]
\centering
\definecolor{pmmark}{RGB}{230,0,0}
\definecolor{pmface}{RGB}{0,100,255}
\newcommand{\pmplus}[2]{\draw[fill=white,line width=0.4pt] (#1,#2) circle (0.2);
  \draw[line width=0.5pt] (#1-0.11,#2)--(#1+0.11,#2) (#1,#2-0.11)--(#1,#2+0.11);}
\newcommand{\pmminus}[2]{\draw[fill=white,line width=0.4pt] (#1,#2) circle (0.2);
  \draw[line width=0.5pt] (#1-0.11,#2)--(#1+0.11,#2);}
\newcommand{\pmring}[2]{\draw[pmmark,line width=1.1pt] (#1,#2) circle (0.25);}
\begin{tikzpicture}[x=0.45cm,y=0.45cm]
\useasboundingbox (-12.5,-3.5) rectangle (12.5,6.9);
\fill[gray!22] (-1.5,-2.5) rectangle (1.5,-1.5)
  (-2.5,-1.5) rectangle (2.5,0.5)
  (-3.5,0.5) rectangle (3.5,3.5)
  (-2.5,3.5) rectangle (2.5,4.5);
\fill[orange!25] (-8.5,0.5) rectangle (-5.5,1.5) (5.5,0.5) rectangle (8.5,1.5)
  (-8.5,1.5) rectangle (-4.5,2.5) (4.5,1.5) rectangle (8.5,2.5)
  (-7.5,2.5) rectangle (-5.5,3.5) (5.5,2.5) rectangle (7.5,3.5)
  (-6.5,3.5) rectangle (-5.5,4.5) (5.5,3.5) rectangle (6.5,4.5)
  (-11.5,1.5) rectangle (-9.5,3.5) (9.5,1.5) rectangle (11.5,3.5)
  (-11.5,3.5) rectangle (-10.5,4.5) (10.5,3.5) rectangle (11.5,4.5);
\draw[gray!55,thin,step=1] (-12.5,-3.5) grid (12.5,6.5);
\draw[line width=1.3pt,line cap=round,line join=round]
  (-12.5,0.5)--(-8.5,0.5)--(-8.5,2.5)--(-7.5,2.5)--(-7.5,3.5)--(-6.5,3.5)
  --(-6.5,4.5)--(-5.5,4.5)--(-5.5,2.5)--(-4.5,2.5)--(-4.5,1.5)--(-5.5,1.5)
  --(-5.5,0.5)--(-3.5,0.5)
  (3.5,0.5)--(5.5,0.5)--(5.5,1.5)--(4.5,1.5)--(4.5,2.5)--(5.5,2.5)--(5.5,4.5)
  --(6.5,4.5)--(6.5,3.5)--(7.5,3.5)--(7.5,2.5)--(8.5,2.5)--(8.5,0.5)--(12.5,0.5)
  (-3.5,1.5)--(3.5,1.5)
  (-11.5,1.5)--(-9.5,1.5)--(-9.5,3.5)--(-10.5,3.5)--(-10.5,4.5)--(-11.5,4.5)--cycle
  (9.5,1.5)--(11.5,1.5)--(11.5,4.5)--(10.5,4.5)--(10.5,3.5)--(9.5,3.5)--cycle;
\draw[pmface,line width=2.2pt] (-3.5,0.5)--(-3.5,1.5) (3.5,0.5)--(3.5,1.5);
\foreach \x in {-12,...,12}{\foreach \y in {-3,...,6}{%
  \ifnum\y<1 \pmplus{\x}{\y}\else\pmminus{\x}{\y}\fi}}
\foreach \x/\y in {-11/2,-11/3,-11/4,-10/2,-10/3,-8/1,-8/2,-7/1,-7/2,-7/3,-6/1,-6/2,
  -6/3,-6/4,-5/2,-3/1,-2/1,-1/1,0/1,1/1,2/1,3/1,5/2,6/1,6/2,6/3,6/4,7/1,7/2,7/3,
  8/1,8/2,10/2,10/3,11/2,11/3,11/4}{\pmplus{\x}{\y}}
\pmring{-3}{1}\pmring{3}{1}
\node[font=\footnotesize] at (0,6.55) {$\Gamma$};
\node[font=\footnotesize] at (-8,6.55) {$\Lambda$};
\node[font=\footnotesize] at (8,6.55) {$\Lambda$};
\end{tikzpicture}
\caption{\footnotesize Perception along the first row of~$\Gamma$, which
switches to~$+1$. The sites of~$\Gamma$ lie in the grey region, and those of~$\Lambda$ in
the orange one. The thick line is the interface. At the two circled sites $\omega$ is boundary-perceptive, and the
blue faces lie on the boundary of~$\Gamma$.}
\label{fig:stickiness}
\end{figure}
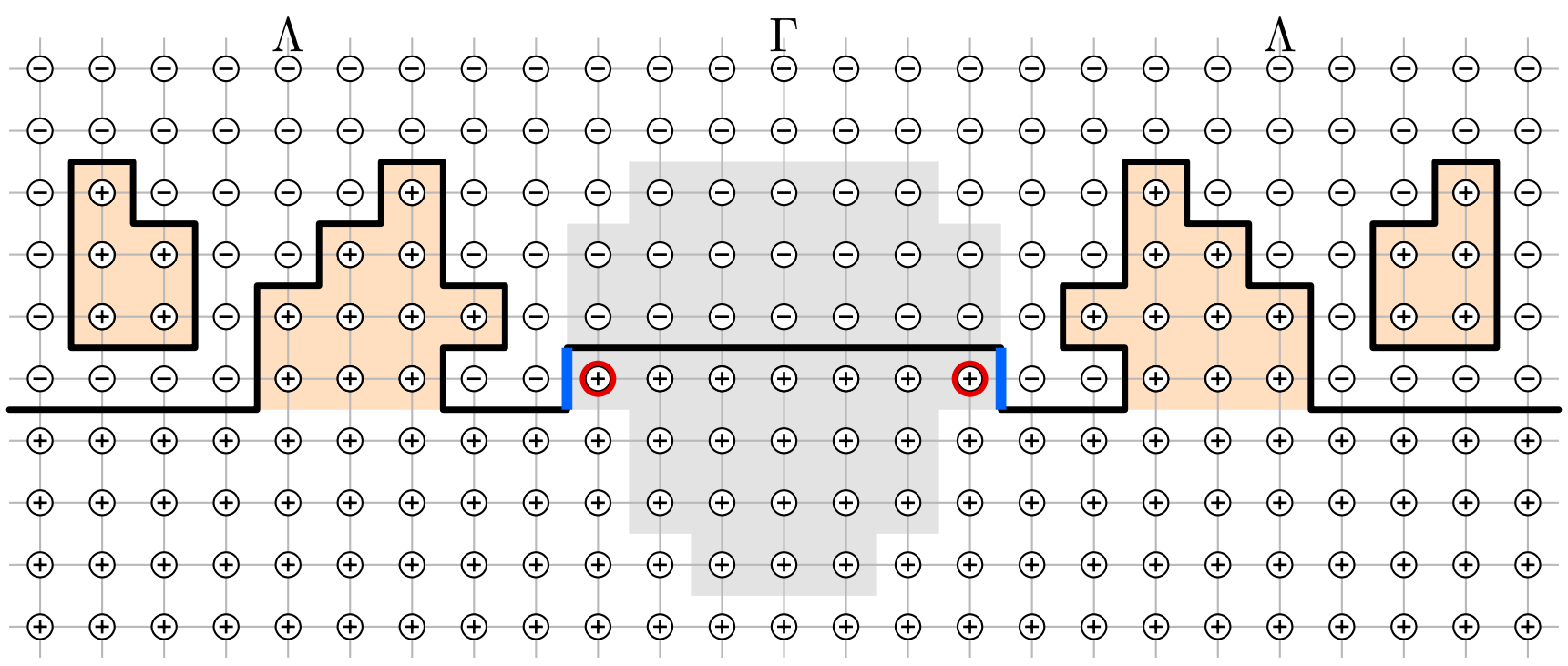

\Cref{fig:stickiness} shows a qualitative picture of~\Cref{th:stickiness} and~\Cref{cor:stickiness}. The whole first row of~$\Gamma$ switches to~$+1$, so that $\omega$
is perceptive at each of its sites. Only the two endpoints of the row belong to
$\partial_\omega\Gamma$, and $\omega$ is boundary-perceptive there. 

We observe that the assumptions of~\Cref{cor:stickiness} can be satisfied. For instance, for the
fractional kernel in the plane and $s$ small, a perturbation at distance $2$ from
$\Gamma$ makes $\omega$ boundary-perceptive at a site of~$\partial_\omega\Gamma$, see~\Cref{prop:fractional}.

We also remark that, in dimension~$d=1$, if~$1\in\Gamma$, the second condition
in~\eqref{eq:instabilitycondition} is always satisfied at the site~$i=1$, provided
that~$K(i,j)>0$ for every~$i\neq j$ and~$\Lambda\neq\varnothing$. Indeed, in this case,
by assumptions~\ref{it:transinv} and~\ref{it:even},
$$ -\sum_{j\neq1}K(1,j)\,\omega_j=\sum_{j\ge2}K(1,j)-\sum_{j\le0}K(1,j)=0.$$
Hence, in dimension~$1$, the half-space configuration is perceptive at the site~$1$
under every perturbation, no matter how small or how far, whenever~$1\in\Gamma$, as it happens in the continuous setting,
see~\cite[Theorem~1.4]{DIPIERRO20171791}.
This is in contrast with the case~$d\ge2$, in which small perturbations do not change the
minimizer, see~\Cref{th:rigidity}.

Our fourth result is inspired by~\cite{dipierro2020disconnectednessnonlocalminimalsurfaces},
where the connectedness of nonlocal minimal surfaces in a cylinder is investigated. We
address the analogous question for Ising-type energies.

Given \(L,h\in\mathbb N\), we consider the box
\begin{equation*}
\mathcal R:=\big\{i=(i',i_d)\in\mathbb Z^{d-1}\times\mathbb Z\ \text{such that}\
\|i'\|_\infty\le L\ \text{and}\ |i_d|\le h\big\},
\end{equation*}
whose base is a cube of side \(2L+1\) and whose height is \(2h+1\). More generally, for
$M\in\mathbb N\cup\{+\infty\}$ with $M\ge h$ we set
\begin{equation*}
\Gamma_M:=\big\{i=(i',i_d)\in\mathbb Z^{d-1}\times\mathbb Z\ \text{such that}\
\|i'\|_\infty\le L\ \text{and}\ |i_d|\le M\big\}.
\end{equation*}
Thus $\Gamma_h=\mathcal R$, and $\Gamma_{+\infty}=:\mathcal G$ is an infinite cylinder.

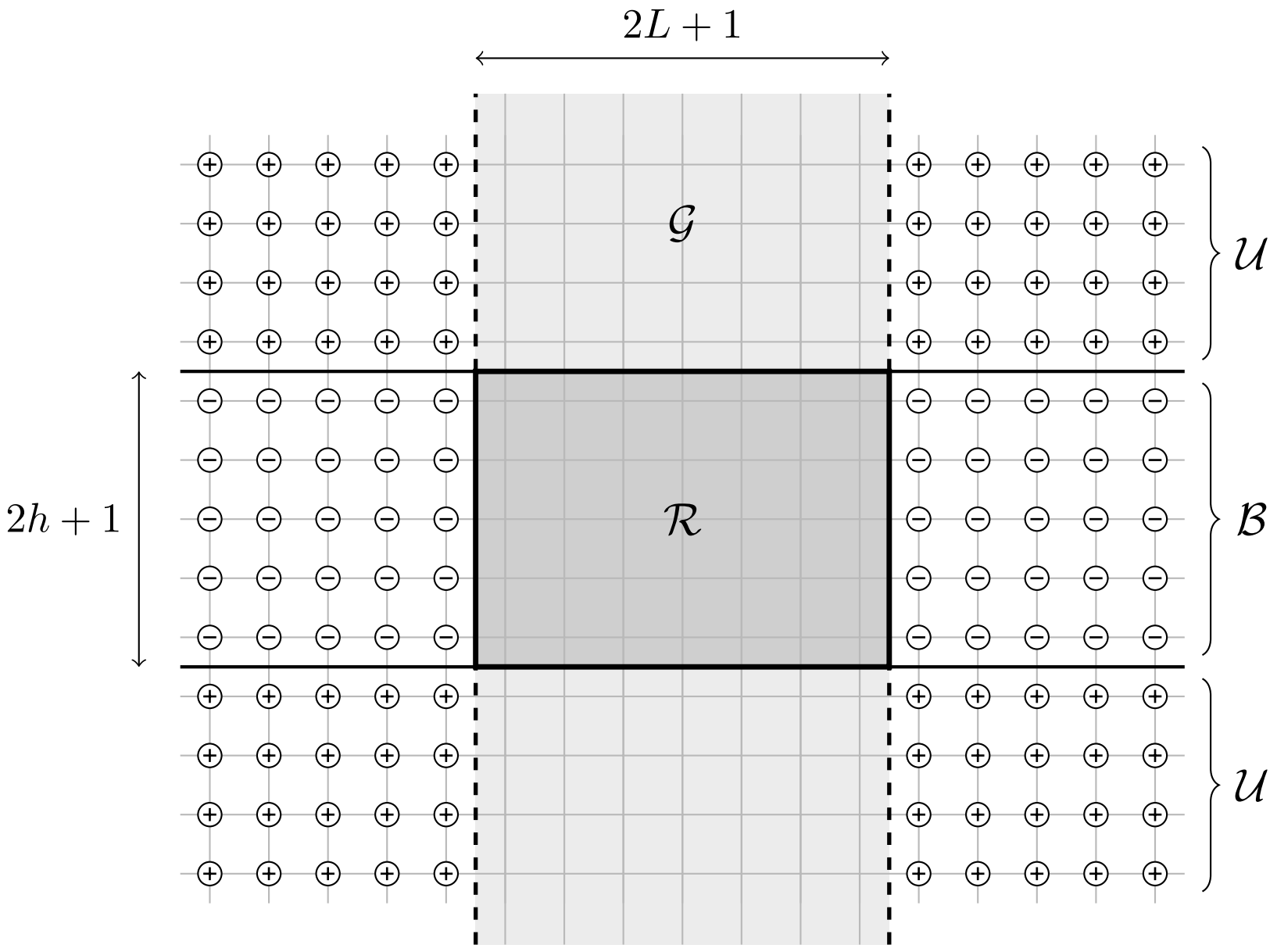
\begin{figure}[H]
\centering
\newcommand{\pmplus}[2]{\draw[fill=white,line width=0.4pt] (#1,#2) circle (0.2);
  \draw[line width=0.5pt] (#1-0.11,#2)--(#1+0.11,#2) (#1,#2-0.11)--(#1,#2+0.11);}
\newcommand{\pmminus}[2]{\draw[fill=white,line width=0.4pt] (#1,#2) circle (0.2);
  \draw[line width=0.5pt] (#1-0.11,#2)--(#1+0.11,#2);}
\begin{tikzpicture}[x=0.5cm,y=0.5cm]
\fill[gray!15] (-3.5,-7.2) rectangle (3.5,7.2);
\fill[gray!38] (-3.5,-2.5) rectangle (3.5,2.5);
\draw[gray!55,thin,step=1] (-8.5,-6.5) grid (8.5,6.5);
\foreach \x in {-3,...,3}{\draw[gray!55,thin] (\x,-7.2) -- (\x,7.2);}
\draw[line width=1.3pt] (-3.5,-2.5) rectangle (3.5,2.5);
\draw[dashed,line width=1pt] (-3.5,-7.2) -- (-3.5,7.2) (3.5,-7.2) -- (3.5,7.2);

\draw[line width=0.8pt] (-8.5,2.5) -- (-3.5,2.5) (-8.5,-2.5) -- (-3.5,-2.5)
  (3.5,2.5) -- (8.5,2.5) (3.5,-2.5) -- (8.5,-2.5);
\foreach \x in {4,...,8}{\foreach \y in {3,...,6}{\pmplus{\x}{\y}\pmplus{\x}{-\y}\pmplus{-\x}{\y}\pmplus{-\x}{-\y}}}
\foreach \x in {4,...,8}{\foreach \y in {-2,...,2}{\pmminus{\x}{\y}\pmminus{-\x}{\y}}}
\node[font=\small] at (0,0) {$\mathcal R$};
\node[font=\small] at (0,5) {$\mathcal G$};
\draw[decorate,decoration={brace,amplitude=4pt}] (8.8,6.3) -- (8.8,2.7)
  node[midway,right=4pt,font=\small] {$\mathcal U$};
\draw[decorate,decoration={brace,amplitude=4pt}] (8.8,2.3) -- (8.8,-2.3)
  node[midway,right=4pt,font=\small] {$\mathcal B$};
\draw[decorate,decoration={brace,amplitude=4pt}] (8.8,-2.7) -- (8.8,-6.3)
  node[midway,right=4pt,font=\small] {$\mathcal U$};
\draw[<->] (-3.5,7.8) -- (3.5,7.8) node[midway,above] {\footnotesize $2L+1$};
\draw[<->] (-9.2,-2.5) -- (-9.2,2.5) node[midway,left] {\footnotesize $2h+1$};
\end{tikzpicture}
\caption{\footnotesize The geometry of~\Cref{th:fumihiko} in the plane. The infinite
cylinder $\mathcal G$ is the grey region bounded by dashed lines. The box
$\mathcal R$ is the darker region bounded by a thick line. Outside $\mathcal G$, the
exterior datum $\eta$ is equal to~$+1$ on $\mathcal U$ and to~$-1$ on
$\mathcal B$. For $h\le M\le+\infty$, the domain $\Gamma_M$ is the part of~$\mathcal G$
with $|i_2|\le M$, and $\eta=+1$ on $\mathcal G\setminus\Gamma_M$.}
\label{fig:regions}
\end{figure}

We also set
\begin{equation}\label{eq:trisplit}\begin{split}
\mathcal U&:=\{j\in\mathbb Z^d\ \text{such that}\ |j_d|>h\}\\
{\mbox{and }}\quad
\mathcal B&:=\{j\in\mathbb Z^d\ \text{such that}\ \|j'\|_\infty>L\ \text{and}\ |j_d|\le h\}.\end{split}
\end{equation}
We call $\mathcal U$ the \emph{caps} and $\mathcal B$ the \emph{belt}, see~\Cref{fig:regions}. Hence
$\mathbb Z^d=\mathcal U\cup\mathcal R\cup\mathcal B$.

We prescribe the exterior datum~$\eta$ as
\begin{equation}\label{eq:slabconfig}
\eta_i:=
\begin{cases}
+1, & {\mbox{ if }}i \in \mathcal U,\\
-1, &{\mbox{ if }}  i \in \mathcal B.
\end{cases}
\end{equation}

The question that we address is whether the configuration
\begin{equation}\label{eq:ustar}
u^\star_i:=
\begin{cases}
+1, & {\mbox{ if }}i\in\mathcal G,\\
\eta_i, &{\mbox{ if }}  i\in\mathbb Z^d\setminus\mathcal G 
\end{cases}
\end{equation}
is a minimizer of~$H_{\Gamma_M}(\cdot\mid\eta)$. 
We recall that, if~$M=+\infty$,
the energy is understood through the differences of~\Cref{rmk:gammainfinite}.

Since $\eta=+1$ on $\mathcal G\setminus\Gamma_M$ for every~$M\ge h$,
the configuration~$u^\star$ is an admissible competitor. Hence the question is whether the minimizer\footnote{Otherwise, the negative exterior datum on $\mathcal B$ nucleates a
droplet of the negative phase inside $\Gamma_M$. In physics, \emph{nucleation} is the
formation of a droplet of one phase inside another, such as a liquid drop in a vapour. The
droplet forms when the energy it gains exceeds the cost of its interface.} is identically~$+1$ in~$\Gamma_M$.

We will see in~\Cref{th:fumihiko} that the answer to this minimization problem is governed by the amount of mass that
the kernel $K$ carries at the scale \(h\). If \(K\) carries enough mass at distances of order \(h\), for instance if it decays
slowly enough, then \(u^\star\) is the minimizer for every \(L\). For compactly supported kernels, instead, \(u^\star\) may fail to be minimal. To formalize
this concept we define the quantity
\begin{equation}\label{eq:PN}
f(l):=\sum_{j\in\mathbb Z^d\setminus\mathcal R}K(l,j)\,u^\star_j ,
\end{equation}
which is an external field
acting on a site \(l\in\mathcal R\), produced by~\(u^\star\) outside
\(\mathcal R\). The sign of \(f(l)\) indicates which phase this field favours at that site.

Besides $u^\star$, one further competitor will play a role, namely the \emph{slab
configuration}
\begin{equation}\label{eq:slabreference}
w_i:=
\begin{cases}
+1, &{\mbox{ if }} |i_d|>h,\\
-1, &{\mbox{ if }} |i_d|\le h,
\end{cases}
\qquad{\mbox{ for all }} i\in\mathbb Z^d .
\end{equation}
The configuration~$w$ is obtained from $u^\star$ by flipping every spin in~$\mathcal R$.

For the nearest-neighbour kernel the answer is complete, and it is decided by the
shape of~$\mathcal R$.

\begin{theorem}[The nearest-neighbour kernel]\label{prop:NN}
Let~$d\ge2$ and let~$K=K_{\rm nn}$ be the nearest-neighbour kernel
in~\eqref{eq:NNkernel}. Then the minimizers of~$H_{\mathcal R}(\cdot\mid\eta)$ are the
following.
\begin{enumerate}[label=\textup{(\roman*)}]
\item\label{it:NNi} The configuration~$u^\star$ alone, if $2L+1>(d-1)(2h+1)$.
\item\label{it:NNii} The configurations $u^\star$ and $w$, and no others, if $2L+1=(d-1)(2h+1)$.
\item\label{it:NNiii} The slab configuration~$w$ alone, if $2L+1<(d-1)(2h+1)$. In this case the
positive phase of the minimizer is $\{i\ \text{such that}\ |i_d|>h\}$, the union of two
disjoint half-spaces. It is therefore disconnected.
\end{enumerate}
\end{theorem}

In the plane, the three cases of~\Cref{prop:NN} are~$h<L$, $h=L$ and~$h>L$. The proof
of~\Cref{prop:NN} is given in~\Cref{sec:NN}. It relies on an explicit discrete vector field
which bounds the perimeter of the droplets from below, in the spirit of the
sub-calibrations of~\cite{DePhilippisPaolini2009}.

For general kernels, we prove the following result.

\begin{theorem}\label{th:fumihiko}
Let $d\ge2$. Then the following hold:
\begin{enumerate}[label=\textup{(\roman*)}]
\item\label{it:rigidity}
If the kernel $K$ satisfies
\begin{equation}\label{eq:mainthreshold}
f(l)>0\qquad\text{for every } l\in\mathcal R,
\end{equation}
then $u^\star$, defined in~\eqref{eq:ustar}, is a minimizer of~$H_{\Gamma_M}(\cdot\mid\eta)$.

Moreover, $u^\star$
is the unique minimizer if either~$M=h$, or~$M>h$ and $K(i,j)>0$ for every $i\ne j$.

For $M=+\infty$, uniqueness is intended among all configurations with exterior
datum~$\eta$ that minimize the energy with respect to finite modifications in~$\mathcal G$.

\item\label{it:fractional}
Let $K_s(i,j):=|i-j|^{-d-2s}$ with
$s\in(0,1)$. There exists $s_\star=s_\star(d,h)>0$ such that~\eqref{eq:mainthreshold}
holds for every $s\in(0,s_\star]$
and every $L\in\mathbb N$.
\end{enumerate}
\end{theorem}

We stress that
the positive phase of~$u^\star$ is
connected, since $u^\star\equiv+1$ on $\mathcal G$ joins the two components of~$\mathcal U$. Thus, 
\Cref{th:fumihiko}(i) says that the minimizer has connected
positive phase. A particular case in which this situation occurs is when the kernel is ``strongly'' nonlocal, as stated in \Cref{th:fumihiko}(ii).
If instead the kernel is ``short-range'', then the positive phase
of the minimizer may be disconnected, see~\Cref{prop:NN}\ref{it:NNiii}.

We remark that the parameter $M$ enters 
in Theorem~\ref{th:fumihiko}
only through the energy~$H_{\Gamma_M}(\cdot\mid\eta)$, and the threshold of
part~\ref{it:fractional} does not depend on it. Hence we prove~\Cref{th:fumihiko} for $M=h$, that is, in the box $\mathcal R$, and we postpone the general case to Appendix~\ref{sec:reduction}.

Part~\ref{it:rigidity} is only a sufficient criterion. Indeed, for the
nearest-neighbour kernel in the plane, \eqref{eq:mainthreshold} fails at the midpoint of
each lateral face of~$\mathcal R$, see~\Cref{rmk:shortrange}, while $u^\star$ is the
unique minimizer if $1\le h<L$, see~\Cref{ex:NNminimal}.

In part~\ref{it:fractional} one can take~$s_\star$ of order~$1/\log h$,
see~\Cref{lem:fracfield}.

Regarding~\Cref{th:fumihiko}(i),
we also stress that when~$M=+\infty$ and~$K$ is not strictly positive, uniqueness
may fail,
see Remark~\ref{CONTEXNU}.
The different features exhibited by the case~$M=+\infty$
are not merely cosmetic. The conceptual distinction
with respect to the case of~$M<+\infty$ is that
two minimal configurations can arise and one cannot always use
one configuration as a competitor for the other, since they may differ on infinitely many sites.
\medskip

In terms of~\Cref{def:sticky}, \Cref{th:fumihiko} shows that the configuration
identically $-1$ is perceptive at every site of~$\mathcal R$ when the spins of
$\mathcal U$ are reversed, see~\Cref{rmk:cylinderinstability}. If, instead, the spins on
$\{|j_d|=h+1\}$ are left equal to~$-1$, then, for fractional kernels with $s$ small, the
configuration identically $-1$ is boundary-perceptive at every site of~$\mathcal R$ having
a nearest neighbour outside $\mathcal R$,
as made precise by the following statement.

\begin{proposition}\label{prop:collar}
Let $d\ge2$, let $h \in\mathbb N$, and let $K=K_s$ be as in~\eqref{eq:frachyp}. Let
$\zeta\equiv-1$ and
\[
\Lambda:=\big\{j\in\mathbb Z^d\ \text{such that}\ |j_d|\ge h+2\big\}.
\]
There exists $s_1=s_1(d,h)>0$ such that, for every $s\in(0,s_1]$ and every
$L\in\mathbb N$, every minimizer of~$H_{\mathcal R}(\cdot\mid\zeta_\Lambda)$ is equal
to~$+1$ in $\mathcal R$. In particular, the configuration identically $-1$ is boundary-perceptive under
the perturbation $\Lambda$ at every site of~$\mathcal R$ having a nearest neighbour
outside $\mathcal R$.
\end{proposition}

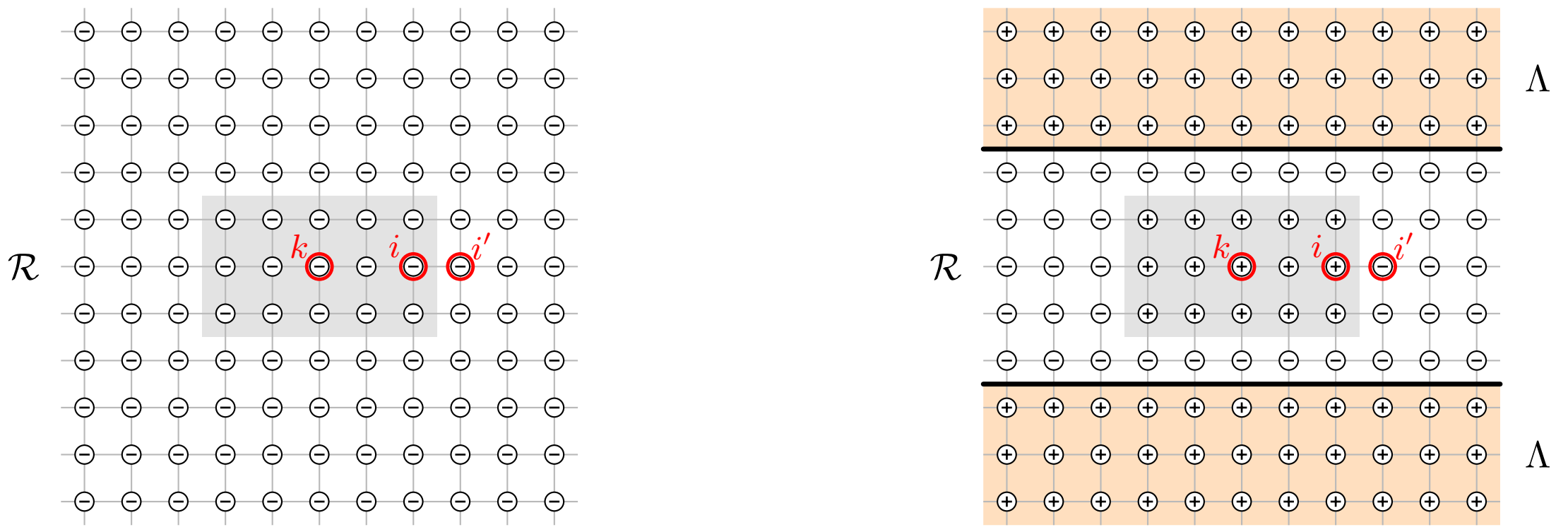
\begin{figure}[H]
\centering
\newcommand{\pmplus}[2]{\draw[fill=white,line width=0.4pt] (#1,#2) circle (0.2);
  \draw[line width=0.5pt] (#1-0.11,#2)--(#1+0.11,#2) (#1,#2-0.11)--(#1,#2+0.11);}
\newcommand{\pmminus}[2]{\draw[fill=white,line width=0.4pt] (#1,#2) circle (0.2);
  \draw[line width=0.5pt] (#1-0.11,#2)--(#1+0.11,#2);}
\newcommand{\pmcollar}[1]{%
\begin{tikzpicture}[x=0.44cm,y=0.44cm]
\useasboundingbox (-7,-5.5) rectangle (6.8,5.5);
\fill[gray!22] (-2.5,-1.5) rectangle (2.5,1.5);
\ifnum#1=1
\fill[orange!25] (-5.5,2.5) rectangle (5.5,5.5) (-5.5,-5.5) rectangle (5.5,-2.5);
\fi
\draw[gray!55,thin,step=1] (-5.5,-5.5) grid (5.5,5.5);
\ifnum#1=1
\draw[line width=1.3pt,line cap=round] (-5.5,2.5)--(5.5,2.5) (-5.5,-2.5)--(5.5,-2.5);
\foreach \x in {-5,...,5}{\foreach \y in {3,4,5}{\pmplus{\x}{\y}\pmplus{\x}{-\y}}}
\foreach \x in {-5,...,5}{\pmminus{\x}{2}\pmminus{\x}{-2}}
\foreach \x in {-5,-4,-3,3,4,5}{\foreach \y in {-1,0,1}{\pmminus{\x}{\y}}}
\foreach \x in {-2,...,2}{\foreach \y in {-1,0,1}{\pmplus{\x}{\y}}}
\node[font=\footnotesize] at (6.3,4) {$\Lambda$};
\node[font=\footnotesize] at (6.3,-4) {$\Lambda$};
\else
\foreach \x in {-5,...,5}{\foreach \y in {-5,...,5}{\pmminus{\x}{\y}}}
\fi
\foreach \x in {0,2,3}{\draw[red,line width=0.9pt] (\x,0) circle (0.27);}
\node[font=\scriptsize,text=red,inner sep=0pt] at (-0.43,0.44) {$k$};
\node[font=\scriptsize,text=red,inner sep=0pt] at (1.60,0.44) {$i$};
\node[font=\scriptsize,text=red,inner sep=0pt] at (3.46,0.44) {$i'$};
\node[font=\footnotesize] at (-6.3,0) {$\mathcal R$};
\end{tikzpicture}}
\begin{minipage}{0.49\textwidth}\centering
\pmcollar{0}
\par\smallskip{\footnotesize (a)}
\end{minipage}\hfill
\begin{minipage}{0.49\textwidth}\centering
\pmcollar{1}
\par\smallskip{\footnotesize (b)}
\end{minipage}
\caption{\footnotesize The setting of~\Cref{prop:collar} in the plane, with $L=2$ and
$h=1$. The sites of~$\mathcal R$ lie in the grey region, and those of~$\Lambda$ in the
orange one, bounded by a thick line. (a) The exterior datum $\zeta\equiv-1$ and the
corresponding minimizer, which is identically $-1$ in $\mathcal R$. (b) The exterior datum
$\zeta_\Lambda$, perturbed in $\Lambda$, and the corresponding minimizer, which is
identically~$+1$ in $\mathcal R$ for $s$ small.
The red circles mark $k=(0,0)$, $i=(2,0)$ and $i'=(3,0)$ in both panels.}
\label{fig:collar}
\end{figure}

See~\Cref{fig:collar} for a visualization of \Cref{prop:collar}.

The following result gives a necessary and sufficient condition for the minimality of
$u^\star$ in the box.

\begin{lemma}\label{prop:criterion}
For $E\subseteq\mathcal R$, let
\begin{equation}\label{eq:Phi}
\Phi(E):=\sum_{l\in E,\, j\in\mathcal R\setminus E}K(l,j)\;+\;\sum_{l\in E}f(l).
\end{equation}
Then $u^\star$ is a minimizer of~$H_{\mathcal R}(\cdot\mid\eta)$ if and only if
\begin{equation}\label{eq:sharp}
\Phi(E)\ge0\qquad\text{for every }E\subseteq\mathcal R .
\end{equation}
Moreover, $u^\star$ is the unique minimizer if and only if $\Phi(E)>0$ for every nonempty
$E\subseteq\mathcal R$.
\end{lemma}

Condition~\eqref{eq:sharp} is harder to check than~\eqref{eq:mainthreshold}, since it
involves every subset of~$\mathcal R$. In simple cases, however, it can be verified for
every $E$. This is the case of the nearest-neighbour kernel, see~\Cref{prop:NN}.

\subsection{Structure of the paper}

\Cref{sec:tools} collects the general tools for the ferromagnetic Ising model that we use
throughout, namely the Euler--Lagrange inequality, the maximum principle, and the notion
of calibration adapted to the lattice setting. In the same section, we prove~\Cref{prop:generic}, which shows that minimizers are unique for generic exterior data, as
in~\cite{DFV2026} for nonlocal minimal surfaces. \Cref{sec:rigidity} contains the proofs
of~\Cref{th:rigidity}, \Cref{th:stickiness} and~\Cref{cor:stickiness}. In~\Cref{sec:fumihiko} we prove~\Cref{prop:criterion}, and we
deduce from it~\Cref{th:fumihiko} in the box and~\Cref{prop:collar}. The proof
of~\Cref{prop:NN} is contained in~\Cref{sec:NN}.
Finally, in Appendix~\ref{sec:reduction} we prove~\Cref{th:fumihiko} in the general case
of the cylinder, by reducing it to the case of the box. In
Appendix~\ref{sec:NNapp} we construct the vector field used in the proof
of~\Cref{prop:NN}, and in
Appendix~\ref{sec:fracfield} we estimate the field of the fractional kernel, which
gives~\Cref{th:fumihiko}\ref{it:fractional} and~\Cref{prop:collar}.

\section{Some tools for the ferromagnetic Ising model}\label{sec:tools}

\subsection{Euler--Lagrange inequality for the Ising model}

We now derive the Euler--Lagrange inequality for the Ising model. Let
$\Gamma\subset\Z^d$ be finite, let $u$ be a minimizer of~$H_\Gamma$, let $l\in\Gamma$, and
let $u^{\{l\}}$ be obtained from $u$ by reversing the spin at $l$. Only the pairs
containing $l$ are affected, each with total weight one in $H_\Gamma$, and
$u^{\{l\}}_l u_j=-u_lu_j$ for $j\neq l$. Hence
\begin{equation}\label{eq:singleflip}
 H_\Gamma(u^{\{l\}})-H_\Gamma(u)
=\sum_{j\neq l}K(l,j)\big[(1+u_lu_j)-(1-u_lu_j)\big]
=2\,u_l\sum_{j\neq l}K(l,j)\,u_j .
\end{equation}
Since $u$ is a minimizer, the left-hand side of~\eqref{eq:singleflip} is nonnegative.
Hence
\begin{equation}\label{eq:EL}
u_l \sum_{j \neq l} K(l,j) u_j \ge 0.
\end{equation}

\begin{remark}
Inequality~\eqref{eq:EL} has a simple interpretation. The quantity
\begin{equation}\label{eq:localfield}
\phi_u(l):=\sum_{j\neq l}K(l,j)u_j
\end{equation}
is the local field\footnote{The terminology local field is standard in statistical mechanics. It denotes the effective field acting on a single spin due to its interaction with the surrounding spins. In our setting, the interactions are long-range, but the field still acts locally, namely at a single site.} acting on the spin at $l$. Therefore, the Euler--Lagrange inequality~\eqref{eq:EL} states that, at every point, the spin has the same sign as its local field (or the local field vanishes). In particular,
\[
u_l=\operatorname{sign}\phi_u(l)
\]
whenever the local field is nonzero.
\end{remark}

\subsection{Calibrations for the Ising model}

In this section, we introduce the notion of calibration for the Ising model, which can be viewed as the discrete counterpart of the calibrations for nonlocal minimal surfaces developed in \cite{Cabr__2020, Pagliari_2019}.

As in the continuum setting, calibrations provide a convenient tool for establishing the minimality properties of suitable configurations.

\begin{definition}\label{def:calib}
Given a configuration~$u$ and a map~$
g:\mathbb{Z}^d \times \mathbb{Z}^d \to \mathbb{R}$,
we say that $g$ is a calibration for $u$ if the following hold:
\begin{itemize}
\item[(i)] (\textit{Antisymmetry condition}) $g$ is antisymmetric, that is, $g(j,i)=-g(i,j)$, and $|g(i,j)| \le 1$ for every $i,j \in \mathbb{Z}^d$;
\item[(ii)] (\textit{Divergence-free condition}) for every $i \in \mathbb{Z}^d$,
\begin{equation}\label{eq:dfc}
\quad \sum_{j \in \mathbb{Z}^d} K(j-i)\,(g(i,j)-g(j,i)) = 0; 
\end{equation}
\item[(iii)] (\textit{Saturation condition}) for every $i \neq j$,
\begin{equation}\label{eq:sc}
    g(i,j)\,(u_j-u_i) = |u_j-u_i|.
\end{equation}
\end{itemize}
\end{definition}

\begin{remark}\label{rmk:anstysim}
The antisymmetry condition in~(i) is not restrictive. Indeed, let $g$ satisfy $|g|\le1$,
\eqref{eq:dfc} and~\eqref{eq:sc}, and set
\[
\tilde{g}(i,j) := \frac{g(i,j)-g(j,i)}{2}.
\]
Then $\tilde g$ is antisymmetric and $|\tilde g|\le1$. It satisfies~\eqref{eq:dfc},
since $\tilde g(i,j)-\tilde g(j,i)=g(i,j)-g(j,i)$. If $u_i\neq u_j$,
then~\eqref{eq:sc} applied to~$(i,j)$ and to~$(j,i)$ gives
$g(i,j)=\operatorname{sign}(u_j-u_i)=-g(j,i)$, hence $\tilde g(i,j)=g(i,j)$. If
$u_i=u_j$, both sides of~\eqref{eq:sc} vanish. Hence $\tilde g$ satisfies~\eqref{eq:sc},
and it is a calibration for $u$.
\end{remark}

Let $\Gamma\subset\mathbb Z^d$ be a finite set and let $\eta$ be a boundary datum.
We consider the class of admissible configurations
\[
\mathcal{F}_{\eta} := \{ v \in \C  \, | \; v = \eta \text{ in } \Z^d \setminus \Gamma \}
\]
and we study the minimization problem
\[
\inf_{v \in \mathcal{F}_{\eta}} H_{\Gamma}(v\mid\eta).
\]

Our calibration criterion is the following.

\begin{theorem}\label{th:calibrcrit}
Let $\Gamma\subset\mathbb Z^d$ be finite, let $\eta$ be a boundary datum, and let $u \in $ $\mathcal{F}_{\eta}$.
If there exists a calibration $g$ for $u$, then
\begin{equation}\label{vncmweiur34y7iWERTY5634}
H_{\Gamma}(u\mid\eta)\leq H_{\Gamma}(v\mid\eta)
\end{equation}
for every admissible competitor $v \in \mathcal{F}_{\eta}$.

Moreover, every minimizer $v\in\mathcal F_\eta$ satisfies the saturation
condition~\eqref{eq:sc} on every pair $\{i,j\}$ with $i\in\Gamma$ and $K(i,j)>0$. In
particular, if $K(i,j)>0$ for every $i\neq j$, then $g$ is a calibration for every
minimizer.
\end{theorem}

\begin{proof}
For $x,y\in\{-1,1\}$ we have $1-xy=|x-y|$. Hence, since $|g|\le1$,
\[
1-v_iv_j=|v_i-v_j|\ge |g(i,j)|\, |v_i-v_j| \ge g(i,j)(v_j-v_i) \]
and similarly
\[
1-v_i\eta_j \ge g(i,j)(\eta_j-v_i).
\]
Multiplying by $K(i,j)\ge0$ and summing, we obtain, for every admissible configuration~$v \in \F_\eta $,
\begin{equation}\label{ineq}
H_{\Gamma}(v\mid\eta) \ge a(v) + b_1(v) + b_0,
\end{equation}
where
\begin{align*}
&a(v):=\frac{1}{2}
\sum_{i, j\in\Gamma}
K(i,j)\,g(i,j)\,(v_j-v_i),\\
&b_1(v):=-\sum_{i\in\Gamma,\, j\in\Z^d\setminus\Gamma}
K(i,j)\,g(i,j)\,v_i,\\ {\mbox{and }}\quad
& b_0:=
\sum_{i\in\Gamma,\, j\in \Z^d\setminus\Gamma}
K(i,j)\,g(i,j)\,\eta_j.
\end{align*}

In order to establish~\eqref{vncmweiur34y7iWERTY5634},
it suffices to prove that
\begin{equation}\label{eq:ab1}
a(v)=-b_1(v)\qquad\text{for every } v\in\mathcal F_\eta .
\end{equation}
Indeed, by~\eqref{eq:ab1}, inequality~\eqref{ineq} gives $H_\Gamma(v\mid\eta)\ge b_0$ for
every $v\in\mathcal F_\eta$. Moreover, by the saturation condition~\eqref{eq:sc}, equality
holds in~\eqref{ineq} for $v=u$, and hence, again by~\eqref{eq:ab1},
\begin{equation}\label{eq:Hub0}
H_\Gamma(u\mid\eta)=b_0 .
\end{equation}

Hence, we now focus on the proof of~\eqref{eq:ab1}.
For this, we use the properties of~$K$ and the fact that~$g$~is antisymmetric
to see that
\begin{equation}\label{vncewuity3479t83fgiuhjhjd43thy}\begin{split}
a(v)
&= \frac{1}{2}
\sum_{i ,j \in \Gamma}
K(j-i)\,g(i,j)\,v_j
-
\frac{1}{2}
\sum_{i , j \in \Gamma}
K(j-i)\,g(i,j)\,v_i \\
&= -
\sum_{i,  j \in \Gamma}
K(j-i)\,g(i,j)\,v_i \\
&= -
\sum_{i \in \Gamma} v_i
\sum_{j \in \Gamma}
K(j-i)\,g(i,j).
\end{split}\end{equation}

By the divergence-free condition~\eqref{eq:dfc} and the antisymmetry of~$g$,
\[
\sum_{j \in \mathbb{Z}^d}
K(j-i)\,g(i,j)=0.
\]
Therefore
\[
\sum_{j \in \Gamma}
K(j-i)\,g(i,j)
=
-\sum_{j \in \Z^d\setminus\Gamma}K(j-i)\,g(i,j).
\]
Hence, plugging this information into~\eqref{vncewuity3479t83fgiuhjhjd43thy},
\begin{align*}
a(v)=
\sum_{i \in \Gamma} v_i
\sum_{j \in \Z^d\setminus\Gamma}
K(j-i)\,g(i,j) =
\sum_{i \in \Gamma,\, j \in \Z^d\setminus\Gamma}
K(j-i)\,g(i,j)\,v_i 
= -b_1(v).
\end{align*}
This proves~\eqref{eq:ab1} and hence the claim in~\eqref{vncmweiur34y7iWERTY5634}.

We now prove the second claim. For $v\in\mathcal F_\eta$ and $i,j\in\mathbb Z^d$, we set
\[
D_{ij}(v):=(1-v_iv_j)-g(i,j)\,(v_j-v_i).
\]
Since $1-v_iv_j=|v_j-v_i|$ and $|g|\le1$, we have 
\begin{equation}\label{DIJneg}
D_{ij}(v)\ge0.\end{equation} Moreover,
\begin{equation}\label{DIJneg2}
{\mbox{$D_{ij}(v)=0$ if and only if~$v$ satisfies the saturation condition~\eqref{eq:sc} on the
pair $\{i,j\}$.}}\end{equation}

Subtracting  $a(v)$, $b_1(v)$ and $b_0$ from $H_\Gamma(v\mid\eta)$ term by term, and recalling that $v=\eta$ outside~$\Gamma$, we obtain
\begin{equation*}
H_\Gamma(v\mid\eta)-a(v)-b_1(v)-b_0
=\frac12\sum_{i,j\in\Gamma}K(i,j)\,D_{ij}(v)
+\sum_{ i\in\Gamma,\, j\in\mathbb Z^d\setminus\Gamma}K(i,j)\,D_{ij}(v).
\end{equation*}
By~\eqref{eq:ab1}, we thus obtain that
\begin{equation}\label{eq:Dsum}
H_\Gamma(v\mid\eta)-b_0
=\frac12\sum_{i,j\in\Gamma}K(i,j)\,D_{ij}(v)
+\sum_{ i\in\Gamma,\, j\in\mathbb Z^d\setminus\Gamma}K(i,j)\,D_{ij}(v).
\end{equation}

If $v$ is a minimizer, then~$H_\Gamma(v\mid\eta)=H_\Gamma(u\mid\eta)=b_0$, in light of~\eqref{eq:Hub0}. Therefore~\eqref{eq:Dsum} boils down to
\begin{equation*}
0=\frac12\sum_{i,j\in\Gamma}K(i,j)\,D_{ij}(v)
+\sum_{ i\in\Gamma,\, j\in\mathbb Z^d\setminus\Gamma}K(i,j)\,D_{ij}(v).
\end{equation*}
Thus, we deduce from this and~\eqref{DIJneg}
that~$D_{ij}(v)=0$ whenever $i\in\Gamma$ and $K(i,j)>0$.
Accordingly, the desired claim follows from this fact and~\eqref{DIJneg2}.

Furthermore, if $K(i,j)>0$ for every~$i\neq j$, then a minimizer~$v$ satisfies the saturation condition on every pair $\{i,j\}$ with at
least one site in $\Gamma$. On the pairs $\{i,j\}$ with both sites outside
$\Gamma$, we have $v=\eta=u$, and the saturation condition holds because $g$ is a calibration for~$u$. Hence $g$ is a calibration for $v$ as well.
\end{proof}

\subsection{Maximum principle for the Ising model}

The main goal of this section is to prove a comparison principle for the Ising model.
The proof relies on the submodularity of the Hamiltonian, which is a well-known
property of ferromagnetic systems, and the comparison principle can be seen as a
zero-temperature version of the FKG and Holley inequalities,
see~\cite{FKG1971, Holley1974}.

Given two configurations \(u,w : \mathbb Z^d \to \{-1,1\}\), we define
\[
m:=u\wedge w:=\min\{u,w\}
\qquad{\mbox{and}}\qquad
M:=u\vee w:=\max\{u,w\}.
\]

\begin{lemma}[Submodularity]\label{lem:submodularity}
For every finite set \(\Gamma\subset \mathbb Z^d\),
\[
H_\Gamma(m)+H_\Gamma(M)
\le
H_\Gamma(u)+H_\Gamma(w).
\]
\end{lemma}

\begin{proof}
It suffices to verify the inequality for the contribution of each interaction pair $(i,j)$.

For every $i,j\in\mathbb Z^d$, one has
\begin{equation}\label{vnewty34ityou49876}
m_i m_j + M_i M_j
\ge
u_i u_j + w_i w_j.
\end{equation}
Indeed, equality holds unless the ordering between $u$ and $w$ is reversed from $i$ to~$j$, namely unless~$u_i=1$, $w_i=-1$,
$u_j=-1$, and~$w_j=1$,
or vice versa. In this case,
\[
u_i u_j+w_i w_j=-2,
\qquad{\mbox{and}}\qquad
m_i m_j+M_i M_j=2,
\] which entails~\eqref{vnewty34ityou49876}.

Therefore, by the ferromagnetic assumption~\ref{it:ferro},
\[
-K(i,j)\big(m_i m_j+M_i M_j\big)
\le
-K(i,j)\big(u_i u_j+w_i w_j\big).
\]
Summing over all interaction pairs yields the desired inequality.
\end{proof}

\begin{remark}[Necessity of the ferromagnetic assumption]\label{rmk:submodnecessary}
\Cref{lem:submodularity} fails without assumption~\ref{it:ferro}. For instance, let $K$
be the nearest-neighbour kernel with $K(i,j)=-1$ if $|i-j|=1$, and let $\Gamma=\{i_0,j_0\}$ be two
neighbouring sites. Let $u$ and $w$ coincide outside $\Gamma$, with~$(u_{i_0},u_{j_0})=(+1,-1)$
and $(w_{i_0},w_{j_0})=(-1,+1)$. Then
$H_\Gamma(m)+H_\Gamma(M)-H_\Gamma(u)-H_\Gamma(w)=4>0$.
\end{remark}

\begin{proposition}[No crossing]\label{prop:comparison}
Let $\Gamma\subset\mathbb Z^d$ be a finite set, and let~$u,v:\mathbb Z^d\to\{-1,1\}$ be minimizers in $\Gamma$. Assume that 
\begin{equation}\label{strictpos5843976}
{\mbox{$K(i,j)>0$ for every $i \neq j$}}\end{equation} and that
\[
u\le v
\qquad\text{in }\mathbb Z^d\setminus\Gamma.
\]
Let
\[
A:=\{i\in\mathbb Z^d:\ u_i>v_i\}
\qquad{\mbox{and}}\qquad
B:=\{i\in\mathbb Z^d:\ u_i<v_i\}.
\]
Then either~$A=\varnothing$ or~$B=\varnothing$.
\end{proposition}

\begin{proof}
Let $m:=u\wedge v$ and $M:=u\vee v$. Since $u\le v$ in $\mathbb Z^d\setminus\Gamma$, the configurations $m$ and $M$ are admissible competitors for~$u$ and~$v$, respectively.
By the submodularity property in~\Cref{lem:submodularity},
\[
H_{\Gamma}(u)+H_{\Gamma}(v)\ge H_{\Gamma}(m)+H_{\Gamma}(M).
\]
Thus, since \(u\) and \(v\) are minimizers, and \(m\) and \(M\) are admissible competitors,
we have
\[
H_{\Gamma}(u)\le H_{\Gamma}(m)
\qquad{\mbox{and}}\qquad
H_{\Gamma}(v)\le H_{\Gamma}(M).
\]
Therefore,
\begin{equation}\label{vnjeuity4i3ugjkvbesdj}
H_{\Gamma}(u)+H_{\Gamma}(v)=H_{\Gamma}(m)+H_{\Gamma}(M).
\end{equation}

We now compute the defect explicitly:
\[
\begin{aligned}
& H_{\Gamma}(u)+H_{\Gamma}(v)-H_{\Gamma}(m)-H_{\Gamma}(M) \\
&=
{\frac12}\sum_{i,j\in\mathbb Z^d}K(i,j)
\Big[
(1-u_iu_j)+(1-v_iv_j)
-(1-m_im_j)-(1-M_iM_j)
\Big] \\
&=
{\frac12}\sum_{i,j\in\mathbb Z^d}K(i,j)
\Big[
m_im_j+M_iM_j-u_iu_j-v_iv_j
\Big].
\end{aligned}
\]

If \(i\) and \(j\) are both either in \(A\) or in \(B\), then the contribution is zero,
since in \(A\) and \(B\) the pair
\((m,M)\) is just a rearrangement of the pair \((u,v)\).
Moreover, outside \(A\cup B\) one has \(u=v\), and then the contribution is zero.

The only nonzero case is when \(i\in A\) and \(j\in B\), or conversely.
If \(i\in A\), then
\[
u_i=1,\qquad v_i=-1,
\qquad
m_i=-1,\qquad M_i=1.
\]
If \(j\in B\), then
\[
u_j=-1,\qquad v_j=1,
\qquad
m_j=-1,\qquad M_j=1.
\]
Hence
\[
m_im_j+M_iM_j-u_iu_j-v_iv_j
=
(-1)(-1)+(1)(1)-(1)(-1)-(-1)(1)
=
4.
\]
Therefore, recalling also~\eqref{vnjeuity4i3ugjkvbesdj},
\begin{equation}\label{vnjewkuoty34793456789iofnjdvsbl}
0=H_{\Gamma}(u)+H_{\Gamma}(v)-H_{\Gamma}(m)-H_{\Gamma}(M)
=
4\sum_{i\in A,\, j\in B}K(i,j).
\end{equation}

Now, if both $A$ and $B$ were nonempty, there would exist
$i\in A$ and $j\in B$. Since $A\cap B=\varnothing$, we have
$i\neq j$, and therefore $K(i,j)>0$, contradicting the
identity in~\eqref{vnjewkuoty34793456789iofnjdvsbl}.
\end{proof}

\begin{remark}[Comparison with the strong comparison principle of Candel--de~la~Llave]
The comparison principle proved in~\Cref{prop:comparison} should be compared with Proposition~19 of~\cite{CandelDeLaLlave1998}, where configurations take values in a linear space and the interaction is ferromagnetic. Under the assumptions of~\cite{CandelDeLaLlave1998}, two ground states that agree at one point must coincide everywhere.

The difference between the two principles reflects a structural distinction between the underlying configuration spaces rather than a difference in the proof techniques. In the Aubry--Mather setting of \cite{CandelDeLaLlave1998}, configurations take values in a linear space, so sums, differences, and convex combinations of admissible configurations remain admissible. This additional algebraic structure is one of the key ingredients underlying the strong comparison principle.

By contrast, the Ising configuration space $\{-1,1\}^{\mathbb Z^d}$ is not linear: if $u$ and $v$ are configurations, neither $u-v$ nor $tu+(1-t)v$, with $t\in(0,1)$, is in general an admissible configuration. Consequently, the mechanism leading to a strong comparison principle is no longer available. What remains is the lattice structure of~$\{-1,1\}^{\mathbb Z^d}$ together with the submodularity of the Hamiltonian, which is sufficient to control the sign of~$u-v$ through the crossing sets, but provides no information once $u=v$ in part of the domain.

For this reason, our comparison principle cannot, in general, be strengthened beyond a no-crossing statement. Indeed, uniqueness of finite-volume minimizers with prescribed exterior datum is false for the Ising model, as shown in Example~\ref{ex:nonunique}, whereas it follows from the strong comparison principle in Corollary~4 of~\cite{CandelDeLaLlave1998}.
\end{remark}

As an immediate consequence of~\Cref{prop:comparison}, we obtain the following maximum principle.

\begin{corollary}[Maximum principle] \label{cor:max}
Let $\Gamma\subset\mathbb Z^d$ be a finite set, and let~$u,v:\mathbb Z^d\to\{-1,1\}$ be minimizers in $\Gamma$. Assume that~$
K(i,j)>0$ for every~$i\neq j$
and that
\[
u\le v
\qquad\text{in }\mathbb Z^d\setminus\Gamma.
\]
If 
\begin{equation}\label{assum8673po21fbjk}
{\mbox{there exists~$i\in\mathbb Z^d\setminus\Gamma$ such that~$u_i< v_i$,}}
\end{equation} then
\[
u\le v
\qquad\text{in }\mathbb Z^d.
\]
\end{corollary}

The strict inequality assumption in~\eqref{assum8673po21fbjk}
is essential to preserve the same ordering of~$u$ and~$v$ in~$\Gamma$. As a matter of fact, if~$u=v$ in $\Z^d\setminus\Gamma$,
then~\Cref{prop:comparison} would give that either~$u\le v$ or~$v\le u$ in~$\Gamma$, but~$u$ and~$v$ may differ, thus violating
not only the maximum principle in
Corollary~\ref{cor:max} but also uniqueness.
An example of this situation is shown here below.

\begin{example}[Non-uniqueness for a given exterior datum]\label{ex:nonunique}
{\rm{Let $\Gamma:=\{(0,0)\}\subset\mathbb Z^2$ and define the exterior
datum $\eta:\mathbb Z^2\setminus\Gamma\to\{-1,+1\}$ by
\[
\eta(i)=
\begin{cases}
+1 & \text{if } i_2<0,\\
-1 & \text{if } i_2>0,\\
+1 & \text{if } i_2=0,\ i_1<0,\\
-1 & \text{if } i_2=0,\ i_1>0,
\end{cases}
\]
see Figure~\ref{fig:nonunique}. Let $u^{\pm}$ be the two admissible
configurations given by~$u^{\pm}=\eta$ in~$\mathbb Z^2\setminus\Gamma$ and~$u^{\pm}(0)=\pm 1$. Then both~$u^{+}$ and~$u^{-}$ are minimizers of~$H_\Gamma(\,\cdot\mid\eta)$.

\begin{center}
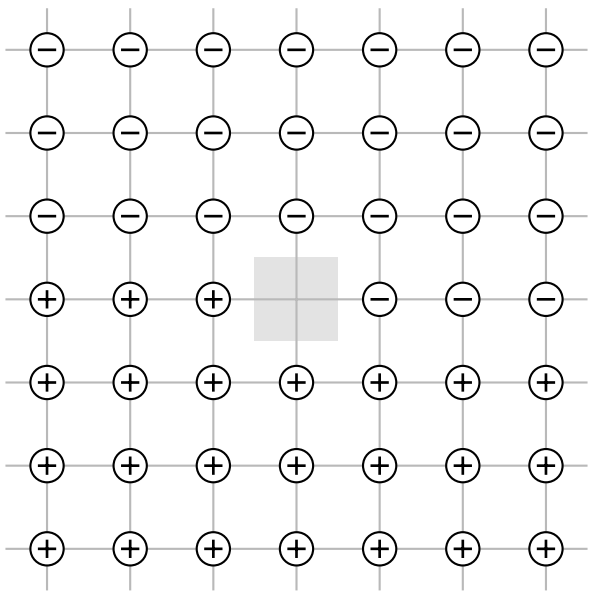

\begin{minipage}{0.78\textwidth}
\centering
\newcommand{\pmplus}[2]{\draw[fill=white,line width=0.4pt] (#1,#2) circle (0.2);
  \draw[line width=0.5pt] (#1-0.11,#2)--(#1+0.11,#2) (#1,#2-0.11)--(#1,#2+0.11);}
\newcommand{\pmminus}[2]{\draw[fill=white,line width=0.4pt] (#1,#2) circle (0.2);
  \draw[line width=0.5pt] (#1-0.11,#2)--(#1+0.11,#2);}
\begin{tikzpicture}[x=0.55cm,y=0.55cm]
\fill[gray!22] (-0.5,-0.5) rectangle (0.5,0.5);
\draw[gray!55,thin,step=1] (-3.5,-3.5) grid (3.5,3.5);
\foreach \x in {-3,...,3}{\foreach \y in {1,...,3}{\pmminus{\x}{\y}\pmplus{\x}{-\y}}}
\foreach \x in {1,...,3}{\pmminus{\x}{0}\pmplus{-\x}{0}}
\end{tikzpicture}
\captionof{figure}{\footnotesize The exterior datum of~\Cref{ex:nonunique}. The set
$\Gamma$ consists of the origin, drawn in grey, where the spin is not prescribed.}
\label{fig:nonunique}
\end{minipage}
\end{center}

To check this, we argue as follows.
The map $j\mapsto-j$ reverses the sign of
the exterior datum, that is, $\eta(-j)=-\eta(j)$ for every $j\neq0$. Moreover
$K(0,-j)=K(0,j)$, by the assumptions on the kernel. Hence
\[
\sum_{j\neq0}K(0,j)\,\eta(j)=0 ,
\]
and therefore
\[
H_\Gamma(u^{+}\mid\eta)-H_\Gamma(u^{-}\mid\eta)
=-2\sum_{j\neq0}K(0,j)\,\eta(j)=0 .
\]
Since $\Gamma$ consists of a single site, $u^{+}$ and $u^{-}$ are the only
admissible configurations, and both are therefore minimizers.
}}\end{example}

\Cref{ex:nonunique} shows that uniqueness may fail for a given exterior datum.
Following~\cite{DFV2026}, we now show that this happens only for exceptional data. More
precisely, along a monotone family of exterior data, the minimizer is unique except for
finitely many data.

\begin{proof}[Proof of~\Cref{prop:generic}]
We follow the argument of~\cite[Theorem~1.1]{DFV2026}. 
Let $\mathcal J$ be the set of~$t\in I$ for which $H_\Gamma(\cdot\mid\eta^{(t)})$ has at least two
distinct minimizers. For $t\in\mathcal J$, let~$u^{(t)}\neq v^{(t)}$ be two such
minimizers. Since they coincide outside $\Gamma$, by~\Cref{prop:comparison}, we
can suppose that $u^{(t)}\le v^{(t)}$.

We set
\begin{equation}\label{eq:Vt}
V^{(t)}:=\{i\in\mathbb Z^d\ \text{such that}\ u^{(t)}_i<v^{(t)}_i\}.
\end{equation}
Then $V^{(t)}$ is a nonempty subset of~$\Gamma$, and
\begin{equation}\label{eq:Vtsigns}
u^{(t)}_i=-1
\qquad\text{and}\qquad
v^{(t)}_i=+1
\qquad\text{ for all }i\in V^{(t)}.
\end{equation}

Let $t,\tau$ be in $\mathcal J$ with $t <\tau$. Outside~$\Gamma$ we have $v^{(t)}=\eta^{(t)}\le\eta^{(\tau)}=u^{(\tau)}$, with strict inequality
at some point by hypothesis. Hence~\Cref{cor:max}, applied to
$v^{(t)}$ and $u^{(\tau)}$, gives
\begin{equation}\label{eq:vtutau}
v^{(t)}\le u^{(\tau)}
\qquad\text{in }\mathbb Z^d.
\end{equation}

By~\eqref{eq:Vtsigns} and~\eqref{eq:vtutau}, $u^{(\tau)}_i=+1$ for all $i\in V^{(t)}$. On the other
hand, $u^{(\tau)}_i=-1$ for all~$i\in V^{(\tau)}$, by~\eqref{eq:Vtsigns}. Therefore
\begin{equation*}
V^{(t)}\cap V^{(\tau)}=\varnothing. 
\end{equation*}

Accordingly, the sets $V^{(t)}$, for $t\in\mathcal J$, are nonempty and pairwise
disjoint subsets of~$\Gamma$. Choosing a site in each of them, we obtain an injective map
from $\mathcal J$ to~$\Gamma$. Hence $\mathcal J$ has at most $\sharp\Gamma$ elements.
\end{proof}

\begin{example}[A continuous family of exterior data]\label{ex:irrational}
{\rm{Let~$d\ge2$, let~$\Gamma\subset\mathbb Z^d$ be finite, and let~$\theta\in\mathbb R^d$
be such that~$\theta_k/\theta_m$ is irrational for some~$k\neq m$ with~$\theta_m\neq0$.
For~$t\in\mathbb R$, we define the exterior
datum~$\eta^{(t)}:\mathbb Z^d\setminus\Gamma\to\{-1,+1\}$ by}
\[{
\eta^{(t)}(j):=
\begin{cases}
+1 & \text{if } \theta\cdot j\le t,\\
-1 & \text{if } \theta\cdot j> t.
\end{cases}}
\]
Then~$\eta^{(t_1)}\le\eta^{(t_2)}$ for every~$t_1<t_2$. Moreover, the
set~$\{\theta_k a+\theta_m b\ \text{such that}\ a,b\in\mathbb Z\}$ is dense in~$\mathbb R$, and so
is the set~$\{\theta\cdot j\ \text{such that}\ j\in\mathbb Z^d\}$. Hence there are
infinitely many sites~$j\in\mathbb Z^d$ such that~$t_1<\theta\cdot j\le t_2$. Since~$\Gamma$
is finite, at least one of them lies outside~$\Gamma$, and
therefore~$\eta^{(t_1)}\neq\eta^{(t_2)}$.

Accordingly, the family~$\{\eta^{(t)}\}_{t\in\mathbb R}$ satisfies the assumptions
of~\Cref{prop:generic} with~$I:=\mathbb R$, and the minimizer is unique for
every~$t\in\mathbb R$, except for at most~$\sharp\Gamma$ values of~$t$.}
\end{example}

We now discuss the strict positivity assumption on the kernel in formula~\eqref{strictpos5843976}. As a matter of fact,
without this assumption the claims in~\Cref{prop:comparison},
\Cref{cor:max} and~\Cref{prop:generic} may be false, as established by the following examples.

\begin{example}\label{lem:strictpos}{\rm{
Let $d=1$, let~$K=K_{\rm nn}$ be the nearest-neighbour kernel
in~\eqref{eq:NNkernel}, and let~$\Gamma:=\{0\}$.

\begin{center}
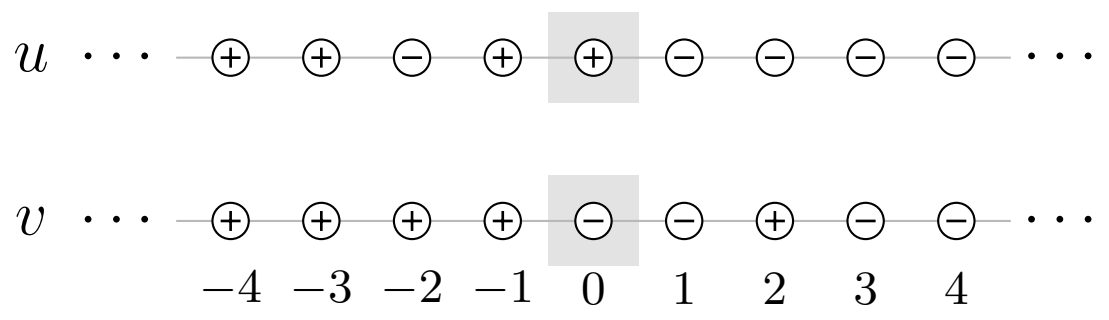

\begin{minipage}{0.85\textwidth}
\centering
\newcommand{\pmplus}[2]{\draw[fill=white,line width=0.4pt] (#1,#2) circle (0.2);
  \draw[line width=0.5pt] (#1-0.11,#2)--(#1+0.11,#2) (#1,#2-0.11)--(#1,#2+0.11);}
\newcommand{\pmminus}[2]{\draw[fill=white,line width=0.4pt] (#1,#2) circle (0.2);
  \draw[line width=0.5pt] (#1-0.11,#2)--(#1+0.11,#2);}
\newcommand{\pmrow}[4]{%
\fill[gray!22] (-0.5,#1-0.5) rectangle (0.5,#1+0.5);
\draw[gray!55,thin] (-4.6,#1) -- (4.6,#1);
\foreach \x in {#2}{\pmplus{\x}{#1}}
\foreach \x in {#3}{\pmminus{\x}{#1}}
\node at (-5.2,#1) {$\cdots$};\node at (5.2,#1) {$\cdots$};
\node at (-6.2,#1) {$#4$};}
\begin{tikzpicture}[x=0.6cm,y=0.6cm]
\pmrow{1.8}{-4,-3,-1,0}{-2,1,2,3,4}{u}
\pmrow{0}{-4,-3,-2,-1,2}{0,1,3,4}{v}
\foreach \x in {-4,...,4}{\node[font=\scriptsize] at (\x,-0.75) {$\x$};}
\end{tikzpicture}
\captionof{figure}{\footnotesize The configurations $u$ and $v$ in~\Cref{lem:strictpos}. The set $\Gamma=\{0\}$ is drawn in grey.}
\label{fig:counterexample}
\end{minipage}
\end{center}

Let $u$ and $v$ be the configurations defined by
\begin{eqnarray*}&&
u_{-2}=u_{2}=-1, \qquad v_{-2}=v_{2}=+1,\\
&& u_{-1}=v_{-1}=+1, \qquad u_{1}=v_{1}=-1,\\
&& u_j=v_j=\operatorname{sign}(-j)\quad  {\mbox{ for }} |j|\ge3,\\
&& u_0=+1\qquad{\mbox{and}}\qquad v_0=-1,
\end{eqnarray*}
see~\Cref{fig:counterexample}.

Since $\Gamma=\{0\}$ and $K$ has range one, the energy in $\Gamma$ depends only on the
interaction of the spin at $0$ with the spins at $\pm1$. These two spins have opposite
signs, and $K(0,-1)=K(0,1)$. Hence their contributions cancel, whatever the value of the
spin at $0$. Thus both values of the spin at $0$ give a minimizer, and $u$ and $v$ are
minimizers in $\Gamma$.

Moreover $u\le v$ in $\mathbb Z\setminus\Gamma$, with $u_{\pm2}<v_{\pm2}$, while~$u_0>v_0$. Hence the sets $\{u>v\}=\{0\}$ and~$\{u<v\}=\{\pm2\}$ are both nonempty, and therefore~\Cref{prop:comparison} fails. Since the inequality outside~$\Gamma$ is strict at~$\pm2$,
also~\Cref{cor:max} fails.}}
\end{example}

\begin{example}{\rm{
Let~$d=1$, let~$K=K_{\rm nn}$ be the nearest-neighbour kernel
in~\eqref{eq:NNkernel},
and let~$\Gamma=\{0\}$.

For $t\in\mathbb N$, let $\eta^{(t)}$
be equal to~$+1$ at the sites~$j\le-1$ and~$2\le j\le t+1$, and to~$-1$ elsewhere. This
family is strictly increasing in $t$. However, $\eta^{(t)}$ changes only at the sites $j\ge2$, and these sites do not interact
with the origin, since $K(0,j)=0$ for $|j|\ge2$. Moreover, the field at the origin is~$K(0,-1)-K(0,1)=0$ for every $t$, so
both values of the spin at the origin give a minimizer. Hence the minimizer is not unique
for every~$t\in\mathbb N$, thus providing a counterexample to~\Cref{prop:generic}.
}}
\end{example}

\section{Rigidity and perception for the half-space}\label{sec:rigidity}

In this section we prove~\Cref{th:rigidity}, \Cref{th:stickiness}, and~\Cref{cor:stickiness}.

\subsection{Minimality of the half-space configuration}
We first prove the minimality of the half-space configuration~$\omega$ (see~\eqref{eq:flatconfig}) by using calibrations.

\begin{lemma}\label{lem:calibH}
The map
\begin{equation}\label{defgcal065489}
g(i,j):=\operatorname{sign}(i_d-j_d)
\end{equation}
is a calibration for $\omega$.
\end{lemma}
\begin{proof}
We check the three conditions in~\Cref{def:calib}.
First, clearly
\[
|g(i,j)|\le 1
\qquad \text{for every } i,j\in\mathbb{Z}^d.
\]
Moreover, $g$ is antisymmetric, since
\[
g(i,j)=\operatorname{sign}(i_d-j_d)
=-\operatorname{sign}(j_d-i_d)
=-g(j,i).
\]

We now check the divergence-free condition~\eqref{eq:dfc}. 
Since $g$ is antisymmetric, for every $i\in\mathbb{Z}^d$,
\begin{eqnarray*} &&
\sum_{j\in\mathbb{Z}^d}
K(j-i)\big(g(j,i)-g(i,j)\big)=2\sum_{j\in\mathbb{Z}^d} K(j-i)g(j,i)\\
&&\qquad=2\sum_{h\in\mathbb{Z}^d} K(h)\operatorname{sign}(h_d)\\
&&\qquad =2\sum_{{h\in\mathbb{Z}^d}\atop{h_d>0}} K(h)\operatorname{sign}(h_d)
+2\sum_{{h\in\mathbb{Z}^d}\atop{h_d<0}} K(h)\operatorname{sign}(h_d)\\
&&\qquad =2\sum_{{h\in\mathbb{Z}^d}\atop{h_d>0}} K(h)\big(\operatorname{sign}(h_d)
+\operatorname{sign}(-h_d)\big)\\
&&\qquad =2\sum_{{h\in\mathbb{Z}^d}\atop{h_d>0}} K(h)\big(\operatorname{sign}(h_d)-\operatorname{sign}(h_d)\big)
=0,
\end{eqnarray*}
as desired.

Finally, we verify the saturation condition~\eqref{eq:sc}. If $\omega_i=\omega_j$, then
\[
g(i,j)\big(\omega_j-\omega_i\big)=0
=
|\omega_j-\omega_i|.
\]
If instead $\omega_i\neq \omega_j$, then $i_d-j_d$ has the same sign as
$\omega_j-\omega_i$. Therefore
\begin{eqnarray*}
g(i,j)\big(\omega_j-\omega_i\big)
= \operatorname{sign}(i_d-j_d)\big(\omega_j-\omega_i\big)
= \operatorname{sign}(\omega_j-\omega_i)\big(\omega_j-\omega_i\big)
=
|\omega_j-\omega_i|.
\end{eqnarray*}
Thus $g$ is a calibration for $\omega$.
\end{proof}

\begin{corollary}\label{cor:flatmin}
The half-space configuration~$\omega$ is a ground state for $H$.   
\end{corollary}
\begin{proof}
The claim follows immediately from~\Cref{th:calibrcrit} and~\Cref{lem:calibH}.
\end{proof}

\subsection{Uniqueness and rigidity of the half-space configuration}

In general, a minimizer with prescribed exterior datum need not be unique, as shown in~\Cref{ex:nonunique}. However, the calibration of~\Cref{lem:calibH} also implies that, in dimension $d\ge2$, the half-space configuration is the unique minimizer.

\begin{lemma}[Uniqueness of the half-space minimizer]\label{lem:uniqueness}
Let $d\ge2$ and let $\Gamma\subset\mathbb Z^d$
be finite. Assume that $K(i,j)>0$ for every $i\neq j$.
Then $\omega$ is the unique minimizer of~$H_\Gamma(\cdot\mid\omega)$.
\end{lemma}

\begin{proof}
Let $u$ be a minimizer of~$H_\Gamma(\cdot\mid\omega)$, and let $g$ be the calibration
of~\Cref{lem:calibH}. By~\Cref{th:calibrcrit}, we have that~$g$ is a calibration also for~$u$.

Let~$i\in\Gamma$. Since $\Gamma$ is finite and $d\ge2$, there exists~$j\in\mathbb Z^d\setminus\Gamma$ with $j_d=i_d$. Then,
recalling the definition of~$g$ in~\eqref{defgcal065489}, we
get that~$g(i,j)=0$.
We note that~$i\neq j$, therefore
the saturation
condition~\eqref{eq:sc} gives that~$u_i=u_j$. 

Hence $u$ is constant on every horizontal
hyperplane. Since $u$ coincides with
$\omega$ outside $\Gamma$, we conclude that~$u=\omega$ in $\mathbb Z^d$.
\end{proof}

\begin{remark}\label{rmk:below}
We point out that perturbations of the exterior datum~$\omega$ supported in the upper half-space~$\Pi^+$ cannot propagate below the flat interface.

Indeed, let $\Lambda\subseteq\Pi^+\setminus\Gamma$ be
nonempty, and let~$\sigma$ be a minimizer of~$H_\Gamma(\cdot\mid\omega_\Lambda)$ (recall the notation in~\eqref{ndjkwhriuy32234567SDFGHVdye4ut}). In~$\Z^d\setminus\Gamma$ we have $\omega\le\omega_\Lambda=\sigma$, with strict inequality in $\Lambda$.
Hence, if~$K(i,j)>0$ for every~$i\neq j$, \Cref{cor:max}
gives~$
\omega\le\sigma$ in~$\mathbb Z^d$.
Since $\omega=+1$ in $\Pi^-$, it follows that $\sigma=+1$ in $\Gamma\cap\Pi^-$.
\end{remark}

Now let $\Lambda\subseteq\Pi^+\setminus\Gamma$. 
For $i\in\Gamma$, we define
\begin{equation}\label{eq:alpha1}
\alpha_i(\Lambda)
:=
\sum_{j\in\Lambda}
K(i,j)
\big[(\omega_\Lambda)_j-\omega_j\big]
=
2\sum_{j\in\Lambda}K(i,j),
\end{equation}
where the last equality of \eqref{eq:alpha1} follows since $\omega_j=-1$ and $(\omega_\Lambda)_j=1$ for every
$j\in\Lambda$. In particular, $\alpha_i(\Lambda)\ge0$ and, recalling the notation
in~\eqref{eq:alpha},
\begin{equation}\label{eq:alpha1BIS}\sum_{i\in\Gamma} \alpha_i(\Lambda) =\alpha_\Gamma(\Lambda).\end{equation}
Also, since $\omega_\Lambda$ and $\omega$ differ only in
$\Lambda\subseteq\mathbb Z^d\setminus\Gamma$, for every configuration~$w$ we have
\begin{equation}\label{eq:datumchange}
H_\Gamma(w\mid\omega_\Lambda)-H_\Gamma(w\mid\omega)
=\sum_{i\in\Gamma,\, j\in\Lambda}K(i,j)\,w_i\big[\omega_j-(\omega_\Lambda)_j\big]
=-\sum_{i\in\Gamma}w_i\,\alpha_i(\Lambda).
\end{equation}

\begin{theorem}[Stability of the flat problem]\label{th:stability}
Let $d\ge2$
and let $\Gamma\subset\mathbb Z^d$
be finite.
Assume that~$K(i,j)>0$ for every~$i\neq j$.
Let~$\delta_\Gamma$ be as in~\Cref{def:gap}. Then~$\delta_\Gamma>0$. Moreover, if
\[
2\,\alpha_\Gamma(\Lambda)<\delta_\Gamma,
\]
then every minimizer $\sigma$ of the perturbed problem coincides with $\omega$
in~$\Gamma$.
\end{theorem}

\begin{proof}
We first show that $\delta_\Gamma>0$. By~\Cref{lem:uniqueness}, $\omega$ is the unique
minimizer of~$H_\Gamma(\cdot\mid\omega)$. Hence every configuration~$v\neq\omega$ has
strictly larger energy. Since the minimum in~\eqref{eq:gap} is taken over finitely many
such configurations, it is positive, as desired.

Now let $v$ be a configuration with $v\neq\omega$ in $\Gamma$. Using~\eqref{eq:datumchange} with~$w:=v$ and with~$w:=\omega$, we obtain, respectively,
\begin{equation*}
H_\Gamma(v\mid\omega_\Lambda)-H_\Gamma(\omega\mid\omega_\Lambda)
=H_\Gamma(v\mid\omega)-H_\Gamma(\omega\mid\omega_\Lambda)
-\sum_{i\in\Gamma}v_i\,\alpha_i(\Lambda).
\end{equation*}
and
\begin{equation*}
H_\Gamma(v\mid\omega)-H_\Gamma(\omega\mid\omega_\Lambda)
=H_\Gamma(v\mid\omega)-H_\Gamma(\omega\mid\omega)
+\sum_{i\in\Gamma}\omega_i\,\alpha_i(\Lambda).
\end{equation*}
Combining these two identities, we find that
\begin{equation}\label{eq:stab3}
H_\Gamma(v\mid\omega_\Lambda)-H_\Gamma(\omega\mid\omega_\Lambda)
=H_\Gamma(v\mid\omega)-H_\Gamma(\omega\mid\omega)
-\sum_{i\in\Gamma}(v_i-\omega_i)\,\alpha_i(\Lambda).
\end{equation}

Since $v\neq\omega$ in $\Gamma$, by~\eqref{eq:gap} we have
\begin{equation}\label{eq:stab4}
H_\Gamma(v\mid\omega)-H_\Gamma(\omega\mid\omega)\ \ge\ \delta_\Gamma .
\end{equation}
Moreover $v_i-\omega_i\le2$ and $\alpha_i(\Lambda)\ge0$, so that, by~\eqref{eq:alpha1BIS},
\begin{equation*}
\sum_{i\in\Gamma}(v_i-\omega_i)\,\alpha_i(\Lambda)\ \le\ 2\sum_{i\in\Gamma}\alpha_i(\Lambda)
=2\,\alpha_\Gamma(\Lambda).
\end{equation*}
Plugging this and~\eqref{eq:stab4} into~\eqref{eq:stab3}, we conclude that
\[
H_\Gamma(v\mid\omega_\Lambda)-H_\Gamma(\omega\mid\omega_\Lambda)
\ \ge\ \delta_\Gamma-2\,\alpha_\Gamma(\Lambda)\ >\ 0 .
\]
Hence $\omega$ is the unique minimizer of the perturbed problem, that is, $\sigma=\omega$
in $\Gamma$ for every minimizer $\sigma$.
\end{proof}

The proof of~\Cref{th:rigidity} now follows from~\Cref{cor:flatmin}, \Cref{lem:uniqueness}, and~\Cref{th:stability}.

\begin{remark}\label{rmk:sharpsite}
For a single site, the thresholds of~\Cref{th:rigidity} and~\Cref{th:stickiness}
are sharp and identify the precise quantity that makes the system
perceptive. Indeed, let~$\Gamma:=\{i\}$ with~$i\in\Pi^+$. 
In this setting, recalling the notation in~\eqref{eq:localfield},
we have that
$$ \alpha_\Gamma(\Lambda)=\alpha_i(\Lambda) \qquad
{\mbox{and}}\qquad \phi_\omega(i)=\sum_{j\neq i}K(i,j)\omega_j.
$$
Moreover, by~\eqref{eq:singleflip},
$$ \delta_\Gamma=H_\Gamma(+1\mid\omega)-H_\Gamma(-1\mid\omega)
=2\omega_i \phi_\omega(i) =-2\phi_\omega(i) .
$$
Accordingly, the thresholds of~\Cref{th:rigidity} and~\Cref{th:stickiness} boil down, respectively, to
\begin{eqnarray*}
2\alpha_i(\Lambda)<\delta_\Gamma \qquad{\mbox{and }}\qquad 2
\alpha_i(\Lambda)> \delta_\Gamma. \end{eqnarray*}

Also, by~\eqref{eq:datumchange},
\begin{eqnarray*}
H_\Gamma(+1\mid\omega_\Lambda)-H_\Gamma(-1\mid\omega_\Lambda)&=&
H_\Gamma(+1\mid\omega)- H_\Gamma(-1\mid\omega)
-2\,\alpha_i(\Lambda)
\\&=&
\delta_\Gamma-2\alpha_i(\Lambda).
\end{eqnarray*}
Hence the minimizer of the perturbed problem is $\omega$ if
$\alpha_i(\Lambda)<\delta_\Gamma/2$, it is equal to~$+1$ at $i$ if~$\alpha_i(\Lambda)>\delta_\Gamma/2$, and both values of the spin at $i$ give a minimizer if
$\alpha_i(\Lambda)=\delta_\Gamma/2$.
\end{remark}

\subsection{Perception for the half-space configuration}
By~\Cref{th:rigidity}, we know that
the half-space configuration is the unique minimizer
with respect to its own exterior datum and,
if~$2\alpha_\Gamma(\Lambda)<\delta_\Gamma$,
it is still a minimizer with respect to a small perturbation of its exterior datum. 
We now consider perturbations above this threshold,
and we show that they may reverse the spin of every minimizer at a given site.

Recall that, for any configuration~$u$, the local field at $i\in\mathbb Z^d$ is defined
in~\eqref{eq:localfield} as
\[
\phi_u(i)=\sum_{j\neq i}K(i,j)u_j.
\]
By~\eqref{eq:localfield} and~\eqref{eq:alpha1}, the second condition
in~\eqref{eq:instabilitycondition} reads
\begin{equation}\label{eq:ipotesi}
\alpha_i(\Lambda)>-\phi_\omega(i).
\end{equation}

\begin{proof}[Proof of~\Cref{th:stickiness}]
Set $M:=\sigma\vee\omega$ and $m:=\sigma\wedge\omega$. In~$\Z^d\setminus\Gamma$ we have
$\sigma=\omega_\Lambda\ge\omega$, hence
\begin{equation}\label{eq:Mmout}
M=\omega_\Lambda
\qquad\text{and}\qquad
m=\omega
\qquad\text{in }\mathbb Z^d\setminus\Gamma.
\end{equation}
Hence, $m$ is a
competitor for~\eqref{eq:flatproblem} and~$M$ is a competitor for~\eqref{eq:perturbedproblem}. 

Also, by~\Cref{lem:submodularity},
\begin{equation}\label{eq:vee1}
H_\Gamma(M\mid\omega_\Lambda)+H_\Gamma(m\mid\omega)
\le H_\Gamma(\sigma\mid\omega_\Lambda)+H_\Gamma(\omega\mid\omega).
\end{equation}
Moreover, since, by~\Cref{cor:flatmin}, we know that~$\omega$ solves~\eqref{eq:flatproblem},
\begin{equation*}
H_\Gamma(m\mid\omega)\ge H_\Gamma(\omega\mid\omega).
\end{equation*}
From this and~\eqref{eq:vee1}, we obtain that
\[
H_\Gamma(M\mid\omega_\Lambda)\le
H_\Gamma(\sigma\mid\omega_\Lambda),
\]
thus, since $\sigma$ is a minimizer
of~\eqref{eq:perturbedproblem}, so is $M$.

By definition, $M\ge\omega$ in $\Gamma$. Moreover $\omega_\Lambda=\omega$ in $\Gamma$, since
$\Lambda\cap\Gamma=\varnothing$. This and~\eqref{eq:Mmout} give
that~$
M\ge\omega_\Lambda$ in~$\mathbb Z^d$.
As a consequence, 
since $K\ge0$, recalling the definition in~\eqref{eq:localfield},
$$ \phi_M(i)=\sum_{j\neq i}K(i,j)M_j \ge
\sum_{j\neq i}K(i,j)(\omega_\Lambda)_j
= \phi_{\omega_\Lambda}(i).
$$

Furthermore, recalling~\eqref{eq:alpha1},
\begin{eqnarray*}&&
\phi_{\omega_\Lambda}(i)=\sum_{j\neq i}K(i,j)(\omega_\Lambda)_j
=\sum_{j\neq i}K(i,j)\omega_j+ \sum_{j\neq i}K(i,j)\big[(\omega_\Lambda)_j-\omega_j\big]\\&&\qquad
=\phi_\omega(i)+ \sum_{j\in\Lambda}K(i,j)\big[(\omega_\Lambda)_j-\omega_j\big]=\phi_\omega(i)+\alpha_i(\Lambda).
\end{eqnarray*}
Combining these pieces of information we conclude that
\begin{equation}\label{eq:fieldsigma}
\phi_M(i) \ge\phi_\omega(i)+\alpha_i(\Lambda).
\end{equation}

We now claim that
\begin{equation}\label{cnjewr4y3t794uigfuSDFGH}
M_i=+1.\end{equation}
Suppose, by contradiction, that $M_i=-1$. Then, we can employ~\eqref{eq:EL} with~$u:=M$ and~$l:=i$
and obtain that
$$ 0\le M_i \sum_{j \neq i} K(i,j) M_j=-\sum_{j \neq i} K(i,j) M_j,
$$ namely~$\phi_M(i)\le0$.

On the other hand, by~\eqref{eq:ipotesi} and~\eqref{eq:fieldsigma},
$$\phi_M(i) \ge\phi_\omega(i)+\alpha_i(\Lambda)>0
$$
and the desired contradiction is established. This proves~\eqref{cnjewr4y3t794uigfuSDFGH}.

Now we recall that, by assumption, $\omega_i=-1$, hence~$M_i=\max\{\sigma_i,-1\}=\sigma_i$.
Therefore, in light of~\eqref{cnjewr4y3t794uigfuSDFGH}, we obtain that~$\sigma_i=+1\neq\omega_i$, as desired.
\end{proof}

\begin{proof}[Proof of~\Cref{cor:stickiness}]
By~\Cref{th:stickiness}, we know that~$\omega$ is perceptive at~$i$ under the perturbation~$\Lambda$.
Since~$\operatorname{dist}(\Lambda,\Gamma)\ge2$ and~$i\in\partial_\omega\Gamma$, the claim
follows immediately from~\Cref{def:pointsticky}.
\end{proof}

We now apply~\Cref{cor:stickiness} to the fractional kernel in the plane.

\begin{proposition}\label{prop:fractional}
Let $d=2$ and $K=K_s$ be as in~\eqref{eq:frachyp}. Let $L\in\mathbb Z$, let $\Gamma$ be a
finite subset of~$\{j\in\Z^2\ \text{such that}\ j_1\le L\}$ containing $i:=(L,1)$, and,
for~$R\in\mathbb N$, let
\[
\Lambda_R:=\big\{j\in\Z^2\ \text{such that}\ L+2\le j_1\le L+R+1\ \text{and}\
1\le j_2\le R\big\}.
\]
Then, if $R$ is large enough, there exists $s^*=s^*(R)>0$ such that, for every
$s\in(0,s^*]$, $\omega$ is boundary-perceptive at $i$ under the perturbation $\Lambda_R$.
\end{proposition}

\begin{proof}
The claim will follow once we show that
the assumptions of~\Cref{cor:stickiness} are satisfied. 

Since $\Gamma\subseteq\{j_1\le L\}$
and $\Lambda_R\subseteq\{j_1\ge L+2\}$, we have $\operatorname{dist}(\Lambda_R,\Gamma)\ge2$.
Moreover, the site~$i':=(L+1,1)$ lies outside $\Gamma$, it is adjacent to~$i$, and~$\omega_{i'}=-1=\omega_i$, hence~$i\in\partial_\omega\Gamma$
(recall the definition in~\eqref{3456cmeuiwdef59ui54bordo}). 

It remains to
check~\eqref{eq:instabilitycondition}. For $s\in(0,1)$ we set
\begin{equation*}
G_R(s):=2\sum_{j\in\Lambda_R}\frac1{|i-j|^{2+2s}}
\qquad{\mbox{and}}\qquad
F(s):=-\sum_{j\neq i}\frac{\omega_j}{|i-j|^{2+2s}} 
\end{equation*} and we observe that,
with $K=K_s$,
condition~\eqref{eq:instabilitycondition} reads
\begin{equation}\label{eq:fracsticky}
G_R(s)>F(s).
\end{equation}
Accordingly, we now focus on proving~\eqref{eq:fracsticky}.

For this, for any~$k\in\mathbb Z$, let
$P_k:=\{j\in\mathbb Z^2\ \text{such that}\ j_2=1+k\}$. By~\eqref{eq:flatconfig} we have
\begin{equation}\label{eq:omegarows}
\omega=-1\ \text{on } P_k \text{ for } k\ge0
\qquad{\mbox{and}} \qquad
\omega=+1\ \text{on } P_k \text{ for } k\le-1 .
\end{equation}

Moreover, for $k\ge1$ we consider
the map $j\mapsto(j_1,2-j_2)$, which is a bijection from $P_k$ onto~$P_{-k}$ and
satisfies
$$ |i-j|^2= (L-j_1)^2+(1-j_2)^2=(L-j_1)^2+(1-(2-j_2))^2=|i-(j_1,2-j_2)|^2.
$$
Therefore, exploiting~\eqref{eq:omegarows}, for every~$k\ge1$,
\begin{equation}\label{eq:rowcancel}\begin{split}&
\sum_{j\in P_k}\frac{\omega_j}{|i-j|^{2+2s}}+\sum_{j\in P_{-k}}\frac{\omega_j}{|i-j|^{2+2s}}
=-\sum_{j\in P_k}\frac{1}{|i-j|^{2+2s}}+\sum_{j\in P_{-k}}\frac{1}{|i-j|^{2+2s}}
\\&\qquad
=-\sum_{j\in P_k}\frac{1}{|i-j|^{2+2s}}+\sum_{j\in P_{k}}\frac{1}{|i-(j_1,2-j_2)|^{2+2s}}
=0
 .\end{split}
\end{equation}

We also note that the series defining $F(s)$ converges absolutely, since~$2+2s>2$. 
Accordingly, we can write
\begin{eqnarray*}
F(s)&=&-\sum_{k\in\Z}\sum_{j\in P_k,\,j\neq i}\frac{\omega_j}{|i-j|^{2+2s}} 
\\
&=&-\sum_{j\in P_0,\,j\neq i}\frac{\omega_j}{|i-j|^{2+2s}} 
-\sum_{k\in\Z\setminus\{0\}} \sum_{j\in P_k }\frac{\omega_j}{|i-j|^{2+2s}}\\
&=&-\sum_{j\in P_0,\,j\neq i}\frac{\omega_j}{|i-j|^{2+2s}} ,
\end{eqnarray*} where~\eqref{eq:rowcancel} has been used in the last line.

Recalling~\eqref{eq:omegarows}, we thus find that
$$ F(s)=\sum_{j\in P_0,\,j\neq i}\frac{1}{|i-j|^{2+2s}}.$$
Notice also that, if~$j\in P_0$, then~$|i-j|=|j_1-L|$, hence
\begin{equation}\label{eq:fracfield}
F(s)=\sum_{k\in\mathbb Z\setminus\{0\}}\frac1{|k|^{2+2s}}=
2\sum_{k=1}^{+\infty}\frac1{k^{2+2s}}.
\end{equation}

The right-hand side of~\eqref{eq:fracfield} is continuous in~$[0,1)$, since~$k^{-2-2s}\le k^{-2}$ and the series converges uniformly, hence $F$ can be extended to a continuous function in~$[0,1)$, by setting
\begin{equation}\label{eq:fracfieldBIS}
F(0):=2\sum_{k=1}^{+\infty}\frac1{k^{2}}<+\infty. \end{equation}

Now we observe that the function $G_R$
is defined by a finite sum, and therefore it is continuous in~$[0,1)$.
We also write
$$\Lambda_R=\{i+(a,b)\ \text{such that}\ 2\le a\le R+1
\ \text{and}\ 0\le b\le R-1\}.$$
Consequently, defining, for every~$n\in\{2,\dots,R\}$, the sets
\[
A_n:=\big\{i+(n,b)\ \text{such that}\ 0\le b\le n-1\big\},
\]
we see that~$A_n$ are pairwise disjoint subsets of~$\Lambda_R$. 

Moreover $\sharp A_n=n$ and, for every $j\in A_n$,
$$|i-j|^2= n^2+b^2\le n^2+(n-1)^2\le2n^2 .$$
Gathering these pieces of information, we conclude that
\begin{equation*}
G_R(0)\ \ge\ 2\sum_{n=2}^{R}\sum_{j\in A_n}\frac1{|i-j|^2}
\ \ge\ 2\sum_{n=2}^{R}\frac{n}{2n^2}=\sum_{n=2}^{R}\frac1n .
\end{equation*}

As a result of this and~\eqref{eq:fracfieldBIS}, 
there exists~$R_0\ge 2$
such that~$G_R(0) > F(0)$ for every~$R\ge R_0$. 
By continuity, for every~$R\ge R_0$
there exists~$s^*=s^*(R)>0$ such that~$G_R(s)>F(s)$ for every~$s\in(0,s^*]$.
This establishes~\eqref{eq:fracsticky} and completes the proof.
\end{proof}

\section{Nucleation of the negative phase in a box}\label{sec:fumihiko}

In this section we prove~\Cref{th:fumihiko} in the box, \Cref{prop:collar}, \Cref{prop:criterion}
and~\Cref{prop:NN}. For this, we will establish some preliminary facts.
 
\subsection{A criterion for minimality}\label{sec:criterionbox}

\begin{lemma} \label{lem:variationformula}
Let $\Gamma\subset\mathbb Z^d$ be finite, let $u:\mathbb Z^d\to\{-1,+1\}$
be a configuration, and let $E\subseteq \Gamma$.
Let~$v$ be the configuration given by 
\[
v_l:=
\begin{cases}
-u_l & \text{if } l\in E,\\
u_l & \text{if } l\notin E.
\end{cases}
\]
Then,
\begin{equation}\label{eq:Glocalfield}
H_{\Gamma}(v)-H_{\Gamma}(u)
=2\sum_{ l\in E,\, j\in \mathbb Z^d \setminus E}K(l,j)\,u_l u_j .
\end{equation}
\end{lemma}
 
\begin{proof}
Since $v$ differs from $u$ only in $E$, we have $v_lv_j=u_lu_j$ whenever
$l$ and $j$ both belong to~$E$ or both belong to~$E^c$. Hence only the
pairs with exactly one endpoint in $E$ contribute to the energy
difference. Every such pair $\{l,j\}$ with $l\in E$ carries total weight
one in $H_\Gamma$: if $j\in\Gamma$, it appears twice in the first double
sum, each time with the factor $\tfrac12$, while if
$j\in\mathbb Z^d\setminus\Gamma$ it appears exactly once in the second
sum. Therefore
\[
H_{\Gamma}(v)-H_{\Gamma}(u)
=\sum_{ l\in E,\, j\in \Z^d\setminus E }K(l,j)\big(u_lu_j-v_lv_j\big)
=2\sum_{ l\in E,\, j\in \Z^d\setminus E}K(l,j)\,u_lu_j ,
\]
where we used that $v_l=-u_l$ for $l\in E$ and $v_j=u_j$ for $j\in \Z^d\setminus E$.
\end{proof}

The quantity appearing in the right-hand side of~\eqref{eq:Glocalfield} is the total field exerted on $E$ by the configuration outside $E$. We formalize this notion in the following definition.

\begin{definition}[Generalized local field]\label{def:field}
For a configuration~$u:\mathbb Z^d\to\{-1,+1\}$ and a finite set
$E\subset\mathbb Z^d$ we set
\[
\Phi_u(E):=\sum_{l\in E,\, j\in \mathbb Z^d \setminus E}K(l,j)\,u_l u_j.
\]
\end{definition}

We stress that if~$E=\{l\}$ one finds that
$$\Phi_u(\{l\})=
\sum_{j\neq l }K(l,j)\,u_l u_j=u_l\,\phi_u(l),$$
where the notation in~\eqref{eq:localfield} has been used.

More generally, for every finite~$E\subset\mathbb Z^d$,
\begin{equation}\label{eq:Phinonadd}
\Phi_u(E)=\sum_{l\in E}u_l\,\phi_u(l)-2\sum_{\{l,j\}\subseteq E}K(l,j)\,u_lu_j ,
\end{equation}
since each pair~$\{l,j\}\subseteq E$ is counted twice in~$\sum_{l\in E}u_l\,\phi_u(l)$
and does not contribute to~$\Phi_u(E)$.

Now we introduce the following notation.
For $X,Y\subseteq\mathbb Z^d$ we write
\begin{equation}\label{eq:IXY}
I(X,Y):=\sum_{l\in X,\,j\in Y}K(l,j)\ \in[0,+\infty].
\end{equation}

The next result shows that, for a defect $E\subseteq\mathcal R$, the field consists of the external field~$f$, introduced in~\eqref{eq:PN},
and of the internal interactions $I$.

\begin{lemma}\label{lem:fieldsplit}
For every $E\subseteq\mathcal R$,
\begin{equation}\label{eq:defectsplit}
\Phi_{u^\star}(E)
=\sum_{l\in E}f(l)\;+\;I\big(E,\mathcal R\setminus E\big).
\end{equation}
In particular, $\Phi_{u^\star}(E)=\Phi(E)$, with $\Phi$ as in~\eqref{eq:Phi}.
\end{lemma}

\begin{proof}
Since $E\subseteq\mathcal R$ and $u^\star=+1$ on $\mathcal R$, we deduce from~\Cref{def:field} that
\begin{equation}\label{cnjsekutiy4t4y35y4y45y}
\Phi_{u^\star}(E)=\sum_{l\in E,\, j\in\mathbb Z^d\setminus E}K(l,j)\,u^\star_j .
\end{equation}
Recalling~\eqref{eq:trisplit}, we see that the complement of~$E$ splits into the three disjoint pieces
$\mathbb Z^d\setminus E=\mathcal U\cup(\mathcal R\setminus E)\cup\mathcal B$, in which
$u^\star$ equals $+1$, $+1$ and $-1$ respectively. Recalling the definition of~$f$ in~\eqref{eq:PN}, we obtain the desired result from~\eqref{cnjsekutiy4t4y35y4y45y}.\end{proof}

We now start proving our main statements. For this,
we observe that, since $\mathcal R$ is finite, the minimum of
$H_{\mathcal R}(\cdot\mid\eta)$ is attained. 

\begin{proof}[Proof of~\Cref{prop:criterion}]
For $E\subseteq\mathcal R$ let $(u^\star)^E$ be the configuration obtained from $u^\star$
by flipping the spins in $E$, namely
\[
((u^\star)^E)_i:=
\begin{cases}
-u^\star_i & \text{if } i\in E,\\
u^\star_i & \text{if } i\notin E.
\end{cases}
\]
We observe that the map
$E\mapsto(u^\star)^E$ is a bijection between the subsets of~$\mathcal R$ and the
competitors, and $(u^\star)^E=u^\star$ in $\mathcal R$ if and only if~$E=\varnothing$.

By~\Cref{lem:variationformula} and~\Cref{lem:fieldsplit},
for every~$E\subseteq\mathcal R$,
\begin{eqnarray*}&&
H_{\mathcal R}\big((u^\star)^E\mid\eta\big)-H_{\mathcal R}(u^\star\mid\eta)
= 2 \sum_{ l\in E,\, j\in \mathbb Z^d \setminus E}K(l,j)\,u^\star_l u^\star_j 
=2\Phi_{u^\star}(E)=2\Phi(E).
\end{eqnarray*}
Hence $u^\star$ is a minimizer if and only if $\Phi(E)\ge0$ for every $E$, and it is the
unique minimizer if and only if $\Phi(E)>0$ for every nonempty $E$, as desired.
\end{proof}

\begin{remark}\label{rmk:mincut}
By~\eqref{eq:Phi} and~\eqref{eq:IXY}, we have that
$$\Phi(E)=I(E,\mathcal R\setminus E)+\sum_{l\in E}f(l).$$
The first term is a nonlocal perimeter of~$E$
relative to~$\mathcal R$ (see in particular~\cite{caffarelli2010nonlocal} for the case~$K=K_s$). The second term is a forcing
term. Hence, by~\Cref{prop:criterion}, $u^\star$ is a minimizer if and only if the
perimeter cost $I(E,\mathcal R\setminus E)$ of every droplet $E\subseteq\mathcal R$ of
the negative phase is at least its gain $-\sum_{l\in E}f(l)$. In combinatorial
optimisation, this is a ``minimum cut problem'', see~\cite{PicardRatliff1975}.
\end{remark}

\begin{theorem}\label{th:connected}
Assume~\eqref{eq:mainthreshold}. Then
$u^\star$ is the unique minimizer of~$H_{\mathcal R}(\,\cdot\mid\eta)$, and its positive
phase is connected.
\end{theorem}

\begin{proof}
By~\Cref{prop:criterion}, it suffices to show that $\Phi_{u^\star}(E)>0$ for every
nonempty $E\subseteq\mathcal R$. Since~$K\ge0$ by assumption, the term
$I(E,\mathcal R\setminus E)$ in~\eqref{eq:defectsplit} is nonnegative, and therefore, thanks to~\eqref{eq:mainthreshold},
\[
\Phi_{u^\star}(E)\ \ge\ \sum_{l\in E}f(l) > 0 ,
\] as desired.

Moreover, recalling~\eqref{eq:trisplit}, we have that~$\{u^\star=+1\}=\mathcal U\cup\mathcal R$,
which entails that the positive
phase of~$u^\star$ is connected.
\end{proof}

\begin{remark}\label{rmk:singleflip}
Stability under single flips does not imply minimality. Indeed, by~\eqref{eq:Phinonadd}
with~$u=u^\star$, for every~$E\subseteq\mathcal R$,
\[
\Phi(E)=\sum_{l\in E}\Phi(\{l\})-2\sum_{\{l,j\}\subseteq E}K(l,j),
\]
so that~$\Phi(E)$ may be negative even if~$\Phi(\{l\})>0$ for every~$l\in E$,
as it happens for the nearest-neighbour kernel in the plane when~$h>L\ge1$, see~\Cref{prop:NN}, which is proved in~\Cref{sec:NN} below.

Moreover, taking~$E=\mathcal R$ in~\eqref{eq:Phi}, we see that the
condition~$\sum_{l\in\mathcal R}f(l)\ge0$ is necessary for the minimality of~$u^\star$.
\end{remark}

\begin{remark}\label{rmk:shortrange}
If $K$ has range at most $h$, assumption~\eqref{eq:mainthreshold} fails.
Indeed, assume that $K(i,j)=0$ whenever $|i-j|>T$, for some $T\le h$. Let $x=(x',0)$,
with $\|x'\|_\infty=L$, be the midpoint of a lateral face of~$\mathcal R$, see~\Cref{fig:wall}. Then
\[
\operatorname{dist}(x,\mathcal B)=1
\qquad\text{and}\qquad
\operatorname{dist}(x,\mathcal U)=h+1>T.
\]

\begin{figure}[H]
\centering
\definecolor{pmmark}{RGB}{230,0,0}
\newcommand{\pmplus}[2]{\draw[fill=white,line width=0.4pt] (#1,#2) circle (0.2);
  \draw[line width=0.5pt] (#1-0.11,#2)--(#1+0.11,#2) (#1,#2-0.11)--(#1,#2+0.11);}
\newcommand{\pmminus}[2]{\draw[fill=white,line width=0.4pt] (#1,#2) circle (0.2);
  \draw[line width=0.5pt] (#1-0.11,#2)--(#1+0.11,#2);}
\begin{tikzpicture}[x=0.6cm,y=0.6cm]
\fill[gray!22] (-4.5,-1.5) rectangle (0.5,1.5);
\draw[gray!55,thin,step=1] (-4.5,-3.5) grid (3.5,3.5);
\draw[line width=1.3pt] (-4.5,1.5) -- (0.5,1.5) -- (0.5,-1.5) -- (-4.5,-1.5);
\draw[line width=1.3pt] (-1,0) -- (1,0) (0,-1) -- (0,1);
\foreach \x in {-4,...,3}{\foreach \y in {2,3}{\pmplus{\x}{\y}\pmplus{\x}{-\y}}}
\foreach \x in {-4,...,0}{\foreach \y in {-1,0,1}{\pmplus{\x}{\y}}}
\foreach \x in {1,2,3}{\foreach \y in {-1,0,1}{\pmminus{\x}{\y}}}
\draw[pmmark,line width=1.1pt] (0,0) circle (0.25);
\node[font=\scriptsize,text=pmmark,inner sep=0.6pt] at (-0.5,-0.5) {$x$};
\draw[decorate,decoration={brace,amplitude=4pt}] (-4.8,-1.3) -- (-4.8,1.3)
  node[midway,left=4pt,font=\small] {$\mathcal R$};
\draw[decorate,decoration={brace,amplitude=4pt}] (3.8,3.3) -- (3.8,1.7)
  node[midway,right=4pt,font=\small] {$\mathcal U$};
\draw[decorate,decoration={brace,amplitude=4pt}] (3.8,1.3) -- (3.8,-1.3)
  node[midway,right=4pt,font=\small] {$\mathcal B$};
\draw[decorate,decoration={brace,amplitude=4pt}] (3.8,-1.7) -- (3.8,-3.3)
  node[midway,right=4pt,font=\small] {$\mathcal U$};
\end{tikzpicture}
\caption{\footnotesize The configuration~$u^\star$ near a site $x$ at the midpoint of a
lateral face of~$\mathcal R$, for $h=1$. Only the right part of~$\mathcal R$ is shown.
The sites of~$\mathcal R$ lie in the grey region, and the thick segments join $x$ to its
nearest neighbours. The set $\mathcal B$ is at distance $1$ from $x$, and $\mathcal U$ is at distance
$h+1$.}
\label{fig:wall}
\end{figure}
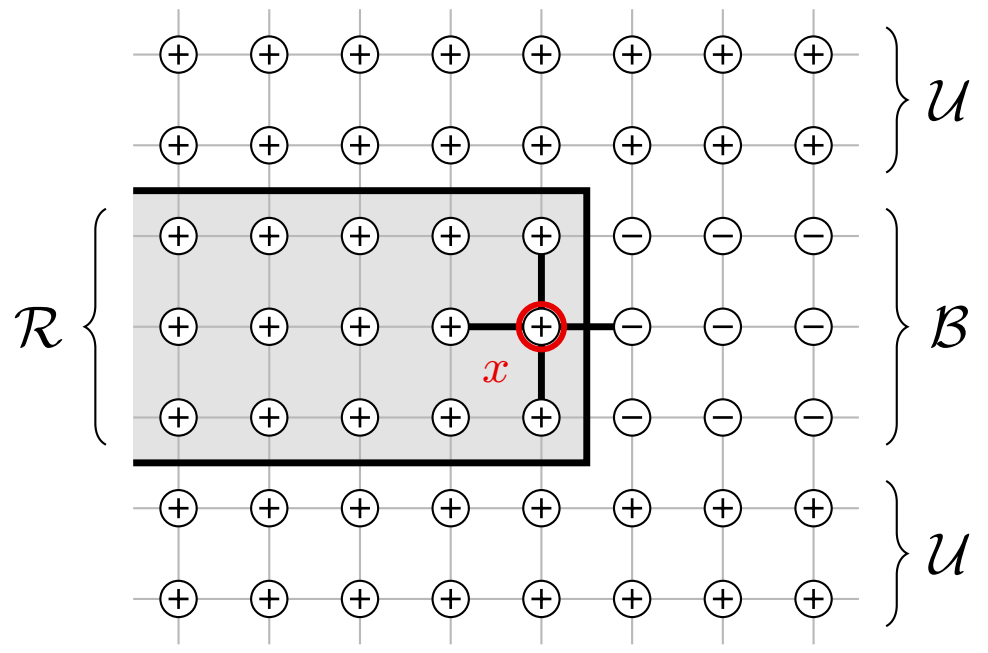

Hence~$\mathcal U$ does not interact with~$x$, and $f(x)=-\sum_{j\in\mathcal B}K(x,j)$. Thus~$f(x)<0$ as soon as~$K$ does not vanish on nearest-neighbour bonds, and
part~\ref{it:rigidity} of~\Cref{th:fumihiko} does not apply.

Nevertheless, $u^\star$ may still be the unique minimizer, see~\Cref{ex:NNminimal}.
\end{remark}

\begin{example}\label{ex:NNminimal}
Let~$d=2$, let~$K=K_{\rm nn}$ be the nearest-neighbour kernel
in~\eqref{eq:NNkernel},
and let~$1\le h<L$. 

Then~$u^\star$ is the unique
minimizer of~$H_{\mathcal R}(\cdot\mid\eta)$, by~\Cref{prop:NN}\ref{it:NNi}, which is
proved in~\Cref{sec:NN} below.
\end{example}

\subsection{The nearest-neighbour kernel}\label{sec:NN}

In this subsection we prove~\Cref{prop:NN}. Hence,~$K=K_{\rm nn}$ is the
nearest-neighbour kernel in~\eqref{eq:NNkernel}, and we denote by~$e_1,\dots,e_d$ the
canonical basis of~$\Z^d$.

We recall that, by~\eqref{eq:defectsplit},
\[
\Phi(E)=I(E,\mathcal R\setminus E) +\sum_{l\in E}f(l).
\]
We bound the perimeter $I(E, \mathcal R \setminus E)$ from below by means of a discrete vector
field, in analogy with the definition of perimeter by duality, see e.g.~\cite{Giusti1984}.

A discrete vector field in~$\mathcal R$ is a function~$g:\mathcal R\times\mathcal R\to\mathbb R$
such that~$g(l,j)=-g(j,l)$ for every~$l$, $j\in\mathcal R$, and~$g(l,j)=0$
whenever~$|l-j|\ne1$. Its divergence\footnote{Compare with the divergence-free condition
in~\Cref{def:calib}.} is
\begin{equation}\label{eq:Dg}
D_g(l):=\sum_{j\in\mathcal R}g(l,j)\qquad{\mbox{for every }}l\in\mathcal R.
\end{equation}

The basic tool of this subsection is the following discrete Gauss--Green formula. For
every discrete vector field~$g$ in~$\mathcal R$ and every~$E\subseteq\mathcal R$,
\begin{equation}\label{eq:GaussGreen}
\sum_{l\in E}D_g(l)=\sum_{{l\in E,\,j\in\mathcal R\setminus E}\atop{|l-j|=1}}g(l,j).
\end{equation}
Indeed, by~\eqref{eq:Dg},
\[
\sum_{l\in E}D_g(l)=\sum_{l,j\in E}g(l,j)+\sum_{l\in E,\,j\in\mathcal R\setminus E}g(l,j).
\]
Now, by the antisymmetry of~$g$,
\[
\sum_{l,j\in E}g(l,j)=\frac12\sum_{l,j\in E}\big(g(l,j)+g(j,l)\big)=0.
\]
Moreover, since~$g(l,j)=0$ if~$|l-j|\neq1$,
\[
\sum_{l\in E,\,j\in\mathcal R\setminus E}g(l,j)
=\sum_{{l\in E,\,j\in\mathcal R\setminus E}\atop{|l-j|=1}}g(l,j),
\]
which completes the proof of~\eqref{eq:GaussGreen}.

The following duality formula\footnote{This is the discrete counterpart of the
notion of perimeter introduced by De Giorgi in~\cite{DeGiorgi1954}, according to which
the perimeter of a measurable set~$E\subseteq\mathbb R^n$ in an open
set~$\Omega\subseteq\mathbb R^n$ is
$$ P(E, \Omega):=\sup\biggl\{\int_E\operatorname{div}\xi\,dx\ :\ \xi\in C^1_c(\Omega;\mathbb R^n)
{\mbox{ and }}|\xi|\le1\biggr\},$$
see~\cite[Definition~3.35]{AmbrosioFuscoPallara2000}. In~\eqref{eq:DeGiorgi}, the role of the compact support
in~$\Omega$ is played by the fact that a discrete vector field in~$\mathcal R$ is defined
only on the pairs of sites of~$\mathcal R$.} is known in the general setting of weighted
graphs. We include its proof, since it is elementary in the case at hand.

\begin{lemma}[see~{\cite[Proposition~3.1]{MazonSoleraToledo2020}}]\label{lem:DeGiorgi}
For every~$E\subseteq\mathcal R$,
\begin{equation}\label{eq:DeGiorgi}
I(E,\mathcal R\setminus E)=\max\Bigl\{\,\sum_{l\in E}D_g(l)\ :\ g{\mbox{ is a discrete vector field in~$\mathcal R$ and }}|g|\le1\Bigr\}.
\end{equation}
\end{lemma}

\begin{proof} We provide a direct and self-contained proof for the convenience of the reader.
Let~$g$ be a discrete vector field in~$\mathcal R$ with~$|g|\le1$. Then, by Gauss-Green formula~\eqref{eq:GaussGreen}, 
\[
\sum_{l\in E}D_g(l) = \sum_{{l\in E,\,j\in\mathcal R\setminus E}\atop{|l-j|=1}}g(l,j)\le I(E,\mathcal R\setminus E).
\]
The equality in \eqref{eq:DeGiorgi} holds if we take
\begin{equation*}
g(l,j):=
\begin{cases}
1, & {\mbox{ if }}l\in E,\ j\in\mathcal R\setminus E{\mbox{ and }}|l-j|=1,\\
-1, & {\mbox{ if }}l\in\mathcal R\setminus E,\ j\in E{\mbox{ and }}|l-j|=1,\\
0, & {\mbox{ otherwise.}}
\end{cases}
\end{equation*}This completes the proof.
\end{proof}

Now we observe that, by~\eqref{eq:defectsplit} and the Gauss--Green formula~\eqref{eq:GaussGreen}, for every discrete vector
field~$g$ in~$\mathcal R$ and every~$E\subseteq\mathcal R$,
\begin{align}
\Phi(E)&=\sum_{l\in E}f(l)+I(E,\mathcal R\setminus E)\notag\\
&=\sum_{l\in E}\big(D_g(l)+f(l)\big)+\Bigl[\,I(E,\mathcal R\setminus E)-\sum_{l\in E}D_g(l)\Bigr]\notag\\
&=\sum_{l\in E}\big(D_g(l)+f(l)\big)+\sum_{{l\in E,\,j\in\mathcal R\setminus E}\atop{|l-j|=1}}\big(1-g(l,j)\big).\label{eq:Phicalib}
\end{align}
In particular, if~$|g|\le1$, then
$$ \Phi(E)\ge\sum_{l\in E}\big(D_g(l)+f(l)\big).$$
Thus, by~\Cref{prop:criterion}, $u^\star$ is a minimizer as soon as there exists a discrete
vector field\footnote{A field with these properties plays the role of a sub-calibration, in
the sense of~\cite[Definition~1.4]{DePhilippisPaolini2009}, see
also~\cite[Theorem~2.6]{Cabr__2020} for nonlocal perimeters.}~$g$ in~$\mathcal R$ such
that~$|g|\le1$ and~$D_g+f\ge0$ in~$\mathcal R$.

We now reduce the inequality~$D_g+f\ge0$ to an equation for $g$ that is simpler to satisfy. For this, we observe that
the average of~$D_g+f$ over~$\mathcal R$ does not depend on~$g$. Indeed, let~$g$ be a
discrete vector field in~$\mathcal R$. By the Gauss--Green
formula~\eqref{eq:GaussGreen} with~$E:=\mathcal R$,
\[
\sum_{l\in\mathcal R}D_g(l)
=\sum_{{l\in\mathcal R,\,j\in \emptyset}\atop{|l-j|=1}}g(l,j)=0,
\]
since the empty set has no near neighbours.
As a consequence,
\[
\frac{1}{\sharp\mathcal R}\sum_{l\in\mathcal R}\big(D_g(l)+f(l)\big)
=\frac{1}{\sharp\mathcal R}\sum_{l\in\mathcal R}f(l).
\]
Since the minimum of a function does not exceed its average, we deduce that
\[
\min_{l\in\mathcal R}\big(D_g(l)+f(l)\big)\le\frac{1}{\sharp\mathcal R}\sum_{l\in\mathcal R}f(l),
\]
with equality if and only if~$D_g+f$ is constant. Hence, the best lower bound
for~$D_g+f$ that a discrete vector field can provide is the average of~$f$, namely
\begin{equation}\label{RALGLE}
\langle f\rangle:=\frac{1}{\sharp\mathcal R}\sum_{l\in\mathcal R}f(l),
\end{equation}
and this bound is attained if and only if~$g$ solves the equation
\begin{equation}\label{eq:NNdiv}
D_g(l)+f(l)=\langle f\rangle\qquad{\mbox{for every }}l\in\mathcal R.
\end{equation}
Accordingly, we look for a discrete vector field~$g$ in~$\mathcal R$ which solves this
equation and satisfies~$|g|\le1$.

We now compute~$\langle f\rangle$. First,
\[
\sharp\mathcal R=(2L+1)^{d-1}(2h+1).
\]
Moreover, by~\eqref{eq:PN}, for every~$l\in\mathcal R$,
\[
f(l)=\chi_{\{l_d=h\}}(l)+\chi_{\{l_d=-h\}}(l)
-\sum_{k=1}^{d-1}\Bigl(\chi_{\{l_k=L\}}(l)+\chi_{\{l_k=-L\}}(l)\Bigr),
\]
and thus
\[
\sum_{l\in\mathcal R}f(l)=2(2L+1)^{d-1}-2(d-1)(2L+1)^{d-2}(2h+1).
\]
Therefore, by the definition of~$\langle f\rangle$,
\begin{equation}\label{eq:F0value}
\langle f\rangle=\frac{2}{2h+1}-\frac{2(d-1)}{2L+1}.
\end{equation}

The next result states that~\eqref{eq:NNdiv} has a solution~$g$ with~$|g|\le1$. Vector
fields with constant divergence, prescribed normal trace and bounded components define
the so-called calibrable sets in the theory of the total variation flow and of the
crystalline mean curvature flow, see~\cite{BellettiniNovagaPaolini2001,
BellettiniCasellesNovaga2002, CasellesChambolleMollNovaga2008}, and
also~\cite{LasicaMollMucha2017} for rectangles with boundary data of both signs
and~\cite[Theorem~5.8]{MazonSoleraToledo2020} for weighted graphs, when the exterior
datum has only one sign. In Appendix~\ref{sec:NNapp} we construct~$g$ explicitly,
see~\eqref{eq:gNN}, from the solution of a Poisson equation with Neumann boundary
conditions.

\begin{lemma}\label{lem:NNcalib}
There exists a discrete vector field~$g$ in~$\mathcal R$ which
solves~\eqref{eq:NNdiv} and satisfies
$$ |g(l,j)|<1\qquad{\mbox{for every }}l,\,j\in\mathcal R.$$
\end{lemma}

\begin{proof}[Proof of~\Cref{prop:NN}]
Let~$g$ be as in~\Cref{lem:NNcalib}. By~\eqref{eq:Phicalib} and~\eqref{eq:NNdiv},
for every~$E\subseteq\mathcal R$,
\begin{equation}\label{eq:PhiNN}
\Phi(E)=\langle f\rangle\,\sharp E+\sum_{{l\in E,\,j\in\mathcal R\setminus E}\atop{|l-j|=1}}\big(1-g(l,j)\big).
\end{equation}
Since~$|g(l,j)|<1$ whenever~$|l-j|=1$, every term of the sum in~\eqref{eq:PhiNN} is
positive. Moreover, since~$\mathcal R$ is connected by nearest-neighbour paths, this sum
is empty if and only if either~$E=\varnothing$ or~$E=\mathcal R$.

We recall, by~\Cref{prop:criterion},  that every competitor is of the form~$(u^\star)^E$ with~$E\subseteq\mathcal R$,
that~$(u^\star)^{\varnothing}=u^\star$ and~$(u^\star)^{\mathcal R}=w$, and that
\begin{equation}\label{eq:competitors}
H_{\mathcal R}\big((u^\star)^E\mid\eta\big)-H_{\mathcal R}(u^\star\mid\eta)=2\Phi(E).
\end{equation}

If~$2L+1>(d-1)(2h+1)$, then~$\langle f\rangle>0$, and so~$\Phi(E)>0$ for every nonempty
set~$E\subseteq\mathcal R$, by~\eqref{eq:PhiNN}. This and~\eqref{eq:competitors}
establish~\ref{it:NNi}.

If~$2L+1=(d-1)(2h+1)$, then~$\langle f\rangle=0$. Thus, by~\eqref{eq:PhiNN}, we have that~$\Phi(E)>0$
for every~$\varnothing\ne E\subsetneq\mathcal R$,
while~$\Phi(\varnothing)=0=\Phi(\mathcal R)$. As a consequence, $u^\star$ and~$w$ are the
only minimizers, which proves~\ref{it:NNii}.

Finally, if~$2L+1<(d-1)(2h+1)$, then~$\langle f\rangle<0$. Hence, by~\eqref{eq:PhiNN},
$\Phi(\mathcal R)=\langle f\rangle\,\sharp\mathcal R$ and
$$ \Phi(E)-\Phi(\mathcal R)\ge-\langle f\rangle\,\sharp(\mathcal R\setminus E)>0
\qquad{\mbox{for every }}E\subsetneq\mathcal R.$$
Therefore, $w$ is the only minimizer, by~\eqref{eq:competitors}. Its positive phase
is~$\{i\ \text{such that}\ |i_d|>h\}$, by~\eqref{eq:slabreference}, and the proof
of~\ref{it:NNiii} is complete.
\end{proof}

We point out that the construction of this subsection relies on the locality of the
nearest-neighbour kernel. Indeed, the locality allows one to
split~\eqref{eq:NNdiv} into one-dimensional problems, one for each direction, see
Appendix~\ref{sec:NNapp}. For a kernel with long-range interactions, every site
of~$\mathcal R$ interacts with the whole exterior datum and with all the other sites
of~$\mathcal R$. For this reason, in the long-range
case we do not construct a vector field, and we rely instead on the direct estimates
of~$f$ in~\Cref{th:connected} and~\Cref{lem:fractional}.

\subsection{The fractional kernel}\label{sec:fracbox}

We now deal with the case of the fractional kernel:

\begin{lemma} \label{lem:fractional}
Let~$K=K_s$ be as in~\eqref{eq:frachyp}.
Then, there exists~$s_1\in(0,1)$, depending only on~$d$ and~$h$, such that, if~$s\in(0, s_1]$, then~$f(l)>0$ for every~$l\in\mathcal R$.
\end{lemma}

\begin{proof}
Let 
\begin{equation}\label{dfeseeselamvda876543}\begin{split}
\Sigma&:=\{j\in\mathbb Z^d\ \text{such that}\ |j_d|=h+1\},\\
\Lambda&:=\{j\in\mathbb Z^d\ \text{such that}\ |j_d|\ge h+2\},\\
{\mbox{and }}\qquad
S&:=\{j\in\mathbb Z^d\ \text{such that}\ |j_d|\le h+1\}\setminus\mathcal R.
\end{split}\end{equation}
 For~$l\in\mathcal R$, we set
\begin{equation}\label{dfeseeselamvda8765432}
{\underline f}(l):=I(\{l\},\Lambda)-I(\{l\},S).
\end{equation}
Then
$\mathcal U=\Lambda\cup\Sigma$ and $S=\mathcal B\cup\Sigma$, 
and so~$\Z^d\setminus\mathcal R=\Lambda\cup\Sigma\cup  \mathcal B$.
Since~$u^\star=+1$ in~$\mathcal U$ and~$u^\star=-1$ in~$\mathcal B$,
\begin{eqnarray*}
&&f(l)=\sum_{j\in\mathbb Z^d\setminus\mathcal R}K(l,j)\,u^\star_j =\sum_{j\in \Lambda}K(l,j) 
+ \sum_{j\in\Sigma}K(l,j) 
- \sum_{j\in \mathcal B}K(l,j) 
\\&&\qquad \quad=
I(\{l\},\Lambda)+I(\{l\},\Sigma)-I(\{l\},\mathcal B)\\
&&\qquad \quad={\underline f}(l)+I(\{l\},S)+I(\{l\},\Sigma)-I(\{l\},\mathcal B)
={\underline f}(l)+2\,I(\{l\},\Sigma)\ \ge\ {\underline f}(l).
\end{eqnarray*}
Accordingly, the claim follows from~\Cref{lem:fracfield}.
\end{proof}

\begin{proof}[Proof of~\Cref{th:fumihiko}]
Part~\ref{it:rigidity} for~$M=h$, including the uniqueness of~$u^\star$, is~\Cref{th:connected}. The case~$M>h$ is proved in Appendix~\ref{sec:reduction}.

Part~\ref{it:fractional} is a consequence of~\Cref{lem:fractional}
and~\Cref{lem:fracfield}, with~$s_\star:=s_1$ as in~\eqref{eq:s1def}.
\end{proof}

\begin{remark}\label{rmk:cylinderinstability}
We point out that~\Cref{th:fumihiko}\ref{it:rigidity} provides an example of perception,
in the sense of~\Cref{def:sticky}, but not of perception at the boundary, in the sense
of~\Cref{def:pointsticky}.

Indeed, let~$\zeta\equiv-1$. Since~$K\ge0$, we have that~$H_{\mathcal R}(v\mid\zeta)\ge0$
for every configuration~$v$, and~$H_{\mathcal R}(\zeta\mid\zeta)=0$. Hence, the
configuration identically~$-1$ is a minimizer of~$H_{\mathcal R}(\cdot\mid\zeta)$.

Moreover, we observe that the exterior datum~$\eta$ defined in~\eqref{eq:slabconfig} satisfies 
\[
\eta = \zeta_{\mathcal U}= \zeta (1-2 \chi_\mathcal U) \quad \text{in $\Z^d \setminus \mathcal R$},
\]
(where we recall that $\zeta_{\mathcal U}$ is defined in \eqref{ndjkwhriuy32234567SDFGHVdye4ut})
 Therefore, if we assume that $f(l)>0$ for each~$l \in R$, we obtain, by applying~\Cref{th:connected}, that the unique minimizer
of~$H_{\mathcal R}(\cdot\mid\zeta_{\mathcal U})$ is~$u^\star$. Thus, since~$u^\star$ is equal to~$+1$
in~$\mathcal R$, $\zeta$ is perceptive at
every site of~$\mathcal R$ under the perturbation~$\mathcal U$.

However, for every~$i'\in\Z^{d-1}$ with~$\|i'\|_\infty\le L$, the
site~$(i',h+1)\in\mathcal U$ is adjacent to the site~$(i',h)\in\mathcal R$, and
so~$\operatorname{dist}(\mathcal U,\mathcal R)=1$. Thus, the perturbation~$\mathcal U$ does
not satisfy the requirement~$\operatorname{dist}(\Lambda,\Gamma)\ge2$
in~\Cref{def:pointsticky}.
\end{remark}

\begin{remark}\label{CONTEXNU}
We point out that when~$M=+\infty$ and~$K$ is not strictly positive, uniqueness may fail,
and therefore~\Cref{th:fumihiko}(i) is sharp. As an example
in which another minimizer (with disconnected positive phase) exists, let~$d\ge2$, $L,h\ge1$, and~$e_d=(0,\ldots,0,1)$. We define~$q:=(2h+2)e_d$ and
$$K(z):=
\begin{cases}
1 & \text{if }z=q\text{ or }z=-q,\\
0 & \text{otherwise}.
\end{cases}$$
This kernel is even and nonnegative, $K(0)=0$, and 
$$\sum_{z\in\Z^d}K(z)=2.$$
Besides, for every $l\in\mathcal R$, we have that~$-h\le l_d\le h$, therefore
$$(l+q)_d\ge -h+(2h+2)=h+2\qquad{\mbox{and}}\qquad
(l-q)_d\le h-(2h+2)=-h-2.$$
Consequently, $l+q$ and $l-q$ both belong to $\mathcal U$, where $u^\star=+1$.

Moreover, $l+q$ and $l-q$ are the only two sites that interact with~$l$. By~\eqref{eq:PN},
$$f(l)=K(l,l+q)u^\star_{l+q}+K(l,l-q)u^\star_{l-q}=1+1=2.$$
Hence, condition~\eqref{eq:mainthreshold} holds in~$\mathcal R$.

Now we set $M:=+\infty$ and define 
$$v_i:=
\begin{cases}
-1 & \text{if }i\in\mathcal G,\\
\eta_i & \text{if }i\notin\mathcal G.
\end{cases}$$
Both $u^\star$ and $v$ have the same exterior datum $\eta$. To check the minimality of both~$u^\star$ and~$v$, we take
$u\in\{u^\star,v\}$ and reverse its spins on any finite set~$E\subseteq\mathcal G$.
We denote the resulting configuration by $u^E$. 

If $i\in E$, both $i+q$ and $i-q$
have the same first $d-1$ coordinates as $i$, so they also belong to $\mathcal G$.
Since $u$ is constant in $\mathcal G$, we have $u_i u_{i+q}=u_i u_{i-q}=1$.

If $i,j\in E$, then $u_i^E u_j^E=(-u_i)(-u_j)=u_i u_j$.
If $i,j\notin E$, then $u_i^E=u_i$ and $u_j^E=u_j$.
In both cases the contribution of $\{i,j\}$ to the energy difference is zero.

If instead~$i\in E$ and~$j\notin E$, then $u_i^E u_j^E=(-u_i)u_j=-u_i u_j$. As a consequence,
\begin{align*}
H_{\mathcal G}(u^E)-H_{\mathcal G}(u)
&=\sum_{i\in E}\sum_{j\notin E}K(j-i)
  \big[(1-u_i^E u_j^E)-(1-u_i u_j)\big]\\
&=2\sum_{i\in E}\sum_{j\notin E}K(j-i)\,u_i u_j .
\end{align*}

Additionally,
given~$i\in E$, the only possible nonzero terms in the inner sum over~$j\notin E$ here above correspond
to~$j=i+q$ and~$j=i-q$. In either case $K(j-i)=1$ and $u_i u_j=1$.
The term with $j=i+q$ therefore contributes $2$ to the energy difference
exactly when $i+q\notin E$. The same calculation applies to $j=i-q$. Hence,
$$H_{\mathcal G}(u^E)-H_{\mathcal G}(u)
=2\,\sharp\{i\in E:i+q\notin E\}
 +2\,\sharp\{i\in E:i-q\notin E\}\ge0.$$
This proves minimality for both configurations~$u^\star$ and~$v$ in the sense of~\Cref{rmk:gammainfinite}.

Finally, by inspection, one observes that the
positive phase of $v$ is disconnected.
\end{remark}

We now prove~\Cref{prop:collar}. This is the case of a perturbation at distance~$2$ from~$\mathcal R$, see~\Cref{rmk:cylinderinstability}.

\begin{proof}[Proof of~\Cref{prop:collar}]
Let~$\Lambda$, $S$ and~${\underline f}$ be as in~\eqref{dfeseeselamvda876543}
and~\eqref{dfeseeselamvda8765432}, let~$s_1$ be as in~\Cref{lem:fracfield}, and
let~$s\in(0,s_1]$. We observe that
\begin{equation}\label{eq:collarsplit}
\Z^d\setminus\mathcal R=\Lambda\cup S,
\qquad{\mbox{with}}\qquad
\zeta_\Lambda=+1 {\mbox{ in }}\Lambda
\qquad{\mbox{and}}\qquad
\zeta_\Lambda=-1 {\mbox{ in }}S,
\end{equation}
and the sets~$\Lambda$ and~$S$ are disjoint.

Let~$u^+$ be
\begin{equation}
u_i^+=
\begin{cases}
+1 & {\mbox{ if }  }i\in \mathcal R,\\
(\zeta_{\Lambda})_i, &{\mbox{ if }} i \in \Z^d \setminus \mathcal R.
\end{cases}
\end{equation}
Given a nonempty set~$E\subseteq\mathcal R$,
let~$v$ be defined as
\begin{equation}\label{defcompet}
v_i=
\begin{cases}
-u_i, & {\mbox{ if }  }i\in E,\\
u_i, &{\mbox{ if }} i \in \Z^d \setminus E.
\end{cases}
\end{equation}

By~\Cref{lem:variationformula} and~\eqref{eq:collarsplit} we have that
\begin{eqnarray*}
&&H_{\mathcal R}(v\mid\zeta_\Lambda)-H_{\mathcal R}(u^+\mid\zeta_\Lambda)
=2\sum_{l\in E,\, j\in\Z^d\setminus E}K(l,j)\,u^+_l u^+_j\\
&&\qquad=2\sum_{l\in E}\Bigg(\sum_{j\in\mathcal R\setminus E}K(l,j)
+\sum_{j\in\Lambda}K(l,j)-\sum_{j\in S}K(l,j)\Bigg)\\
&&\qquad=2\,I(E,\mathcal R\setminus E)+2\sum_{l\in E}{\underline f}(l).
\end{eqnarray*}
Since~$K\ge0$, we have that~$I(E,\mathcal R\setminus E)\ge0$. Moreover, by applying~\Cref{lem:fracfield}, we obtain that ${\underline f}(l)>0$ for
every~$l\in\mathcal R$ and every~$L\in\mathbb N$.
Thus,
\begin{equation}\label{diffcompet}
  H_{\mathcal R}(v\mid\zeta_\Lambda)-H_{\mathcal R}(u^+\mid\zeta_\Lambda) >0
\end{equation}
Since the set $E \subset \mathcal R$ taken in the definition of~$v$ in~\eqref{defcompet} is arbitrary, we can observe that~\eqref{diffcompet} holds true for every competitor $v$ of $u^+$ in $\mathcal R$. In particular, we conclude that~$u^+$ is the
unique minimizer of~$H_{\mathcal R}(\cdot\mid\zeta_\Lambda)$. This establishes the first
claim in~\Cref{prop:collar}.

We now prove the second claim. We showed  in~\Cref{rmk:cylinderinstability} that $\zeta$ is a minimizer of~$H_{\mathcal R}(\cdot\mid\zeta)$. Moreover, we observe that~$\operatorname{dist}(\Lambda,\mathcal R)\ge2$. Indeed, if~$j\in\Lambda$ and~$i\in\mathcal R$, then~$|j_d|\ge h+2$ and~$|i_d|\le h$.

Let now~$i\in\mathcal R$ be a site having a nearest neighbour~$i'\in\Z^d\setminus\mathcal R$.
Since $\zeta_i = \zeta_{i'}$, we have that~$i\in\partial_{\zeta}\mathcal R$ (see~\eqref{3456cmeuiwdef59ui54bordo}).
Since every minimizer of~$H_{\mathcal R}(\cdot\mid\zeta_\Lambda)$ is equal
to~$+1$ at~$i$, we have that~$\zeta$ is boundary-perceptive at~$i$ under the
perturbation~$\Lambda$.
\end{proof}

\appendix
\section{Reduction to the box}\label{sec:reduction}

In this appendix we complete the proof of~\Cref{th:fumihiko}, by reducing the
case~$M>h$ to the case~$M=h$, which has been addressed in~\Cref{sec:fumihiko}. The
case~$M=+\infty$ is motivated
by~\cite{dipierro2020disconnectednessnonlocalminimalsurfaces}, where the ambient domain
is an infinite cylinder.

Throughout this appendix, we pick~$M\in\mathbb N\cup\{+\infty\}$ with~$M\ge h$, and we
denote by~$\mathcal F$ the class of the configurations which differ from~$u^\star$,
defined in~\eqref{eq:ustar}, on a finite subset of~$\Gamma_M$. If~$M$ is finite, these
are all the configurations which coincide with~$\eta$ outside~$\Gamma_M$.
Instead, when~$M=+\infty$,
\Cref{rmk:gammainfinite} tests each configuration against its own finite modifications
in~$\mathcal G$. 

In both cases, we only need to compare configurations which differ at finitely many
sites. Given two such configurations~$u$ and~$v$, we denote by~$H(u)-H(v)$ their energy
difference, as in~\Cref{sec:assumptions}, namely
\begin{equation}\label{eq:Hdiff}
H(u)-H(v):=\sum_{\{l,j\}}K(l,j)\,\big(v_lv_j-u_lu_j\big),
\end{equation}
where the sum is over the unordered pairs~$\{l,j\}$ with~$l\neq j$. Only the pairs with
at least one site in the finite set~$\{u\neq v\}$ contribute, and so the series converges
absolutely, by assumption~\ref{it:summable}. Hence, a configuration~$\sigma\in\mathcal F$ is a minimizer
of~$H_{\Gamma_M}(\cdot\mid\eta)$ if and only if~$H(\sigma)-H(v)\le0$ for
every~$v\in\mathcal F$. Moreover, the configurations which coincide with~$\eta$
outside~$\mathcal R$ belong to~$\mathcal F$, and for them~$H(u)-H(v)=H_{\mathcal
R}(u\mid\eta)-H_{\mathcal R}(v\mid\eta)$.
\medskip

We stress that the submodularity of the Hamiltonian in~\Cref{lem:submodularity}, namely
\[
H_\Gamma(v \wedge \psi)+H_\Gamma(v \vee \psi)
\le
H_\Gamma(v)+H_\Gamma(\psi),
\]
is stated for a finite domain~$\Gamma$, since otherwise the energies are not finite.
However, the domain~$\Gamma_M$ is infinite when~$M=+\infty$. To this end, we
observe that, if the energies were finite, moving~$H_\Gamma(v)$ to the left hand side
and~$H_\Gamma(v\wedge\psi)$ to the right hand side would give
\begin{equation}\label{eq:Gsubmod}
H(v \vee \psi)- H(v)
\le
H(\psi) -H(v \wedge \psi).
\end{equation}
This inequality only involves energy differences, and therefore it makes sense in every
domain, in the sense of~\eqref{eq:Hdiff}, as soon as the set~$\{v<\psi\}$ is finite.
It is the main tool in the proof of the generalised comparison principle
in~\Cref{lem:crossing}, which we will use to reduce the cylinder to the box.

\medskip

We prove a quantitative version of~\eqref{eq:Gsubmod}, in which the defect is the
interaction between the two sets where the configurations cross. 
For this, we recall
the notation in~\eqref{eq:IXY}, namely, for~$X$, $Y\subseteq\Z^d$,
\[
I(X,Y)=\sum_{l\in X,\,j\in Y}K(l,j).
\]

\begin{lemma}\label{lem:Gsubmod}
Let~$v$ and~$\psi$ be configurations, let
\[
A:=\{i\in\Z^d\ \text{such that}\ v_i<\psi_i\}
\qquad{\mbox{and}}\qquad
B:=\{i\in\Z^d\ \text{such that}\ v_i>\psi_i\},
\]
and assume that~$A$ is finite. Then,
\begin{equation}\label{eq:subdiffAB}
H(v\vee\psi)-H(v)=H(\psi)-H(v\wedge\psi)-4\,I(A,B).
\end{equation}
In particular, \eqref{eq:Gsubmod} holds true.
\end{lemma}

\begin{proof}
The configuration~$v\vee\psi$ differs from~$v$ only in~$A$, and~$v\wedge\psi$ differs
from~$\psi$ only in~$A$. Since~$A$ is finite, by~\eqref{eq:Hdiff} we have that
\begin{equation}\label{eq:defectsum}
\big[H(v\vee\psi)-H(v)\big]-\big[H(\psi)-H(v\wedge\psi)\big]=-\sum_{\{l,j\}}K(l,j)\,D_{lj},
\end{equation}
where
$$ D_{lj}:=(v\wedge\psi)_l(v\wedge\psi)_j+(v\vee\psi)_l(v\vee\psi)_j-v_lv_j-\psi_l\psi_j.$$
If~$v\le\psi$ at both~$l$ and~$j$, then~$v\wedge\psi=v$ and~$v\vee\psi=\psi$ at both
sites, and so~$D_{lj}=0$. The same holds if~$v\ge\psi$ at both sites. Otherwise, one of
the two sites belongs to~$A$ and the other one to~$B$. In this case, $v\wedge\psi=-1$
and~$v\vee\psi=+1$ at both sites, while~$v_lv_j=\psi_l\psi_j=-1$, and so~$D_{lj}=4$.
Since~$A$ and~$B$ are disjoint, we see that
the right hand side of~\eqref{eq:defectsum} is equal
to~$-4\,I(A,B)$, which gives~\eqref{eq:subdiffAB}. Since~$I(A,B)\ge0$, this
implies~\eqref{eq:Gsubmod}.
\end{proof}

We can now state and prove the comparison principle.

\begin{lemma}\label{lem:crossing}
Let~$\Gamma\subseteq\Z^d$ and let~$\psi$ be a minimizer in~$\Gamma$, in the sense
of~\Cref{rmk:gammainfinite}. Let~$v$ be a configuration such that the
set~$A:=\{i\in\Z^d\ \text{such that}\ v_i<\psi_i\}$ is a finite subset of~$\Gamma$. Then,
\begin{equation}\label{eq:crossing}
H(v\vee\psi)-H(v)\le0.
\end{equation}
\end{lemma}

\begin{proof}
The configuration~$v\wedge\psi$ differs from~$\psi$ only in the finite
set~$A\subseteq\Gamma$, and therefore it is an admissible competitor for~$\psi$.
Since~$\psi$ is a minimizer in~$\Gamma$, this gives that~$H(\psi)-H(v\wedge\psi)\le0$,
which, together with~\eqref{eq:Gsubmod}, yields~\eqref{eq:crossing}.
\end{proof}

The reduction to the box
consists of two steps. First, we show that every minimizer in the box is a minimizer
in~$\Gamma_M$ (see~\Cref{lem:reduction}). This does not exclude that there are other
minimizers in~$\Gamma_M$, which take the value~$-1$ somewhere outside the box,
since~\eqref{eq:crossing} is not a strict inequality. Then, we show that,
if~$K>0$, there are no such minimizers, namely every minimizer in~$\Gamma_M$ is equal
to~$+1$ in~$\Gamma_M\setminus\mathcal R$ (see~\Cref{cor:cylinder}). Hence, the
minimizers in~$\Gamma_M$ coincide with those in the box, and uniqueness in the box
implies uniqueness in~$\Gamma_M$.

\begin{proposition}\label{lem:reduction}
Let~$v\in\mathcal F$ and let~$w$ be as in~\eqref{eq:slabreference}. Then, the
configuration~$v\vee w$ belongs to~$\mathcal F$, it coincides with~$v$ in~$\mathcal R$
and with~$+1$ in~$\Gamma_M\setminus\mathcal R$, and
\begin{equation}\label{eq:reduction}
H(v\vee w)-H(v)\le0.
\end{equation}
In particular, every minimizer of~$H_{\mathcal R}(\cdot\mid\eta)$ is a minimizer
of~$H_{\Gamma_M}(\cdot\mid\eta)$.
\end{proposition}

\begin{proof}
Let~$\omega$ be as in~\eqref{eq:flatconfig} and, for every~$i=(i',i_d)\in\Z^d$, let
$$ \psi_2(i):=\omega_{(i',\,i_d+h+1)}\qquad{\mbox{and}}\qquad\psi_1(i):=\psi_2(-i),$$
so that~$\psi_1=+1$ if and only if~$i_d>h$, and~$\psi_2=+1$ if and only if~$i_d<-h$.
In particular, $\psi_1\vee\psi_2=w$. Moreover, $\omega$ is a ground state,
by~\Cref{cor:flatmin}, and the energy differences are invariant under translations and
under the map~$i\mapsto-i$, by assumptions~\ref{it:transinv} and~\ref{it:even}.
Hence, $\psi_1$ and~$\psi_2$ are ground states, and thus minimizers in~$\Gamma_M$.

Let now~$A_1:=\{v<\psi_1\}$. If~$i\in A_1$, then~$v_i=-1$ and~$i_d>h$, and
therefore~$u^\star_i=+1\neq v_i$, by~\eqref{eq:slabconfig} and~\eqref{eq:ustar}.

Thus, $A_1$ is contained in the set where~$v$ differs from~$u^\star$, which is a finite
subset of~$\Gamma_M$. By~\Cref{lem:crossing}, the configuration~$v_1:=v\vee\psi_1$
satisfies~$H(v_1)-H(v)\le0$. Also, $v_1$ coincides with~$v$ outside~$A_1$ and
with~$u^\star$ in~$A_1$, and so~$v_1\in\mathcal F$.

The same argument, applied to~$v_1$ and~$\psi_2$, gives
that~$v_1\vee\psi_2=v\vee w$ belongs to~$\mathcal F$ and
that~$H(v\vee w)-H(v_1)\le0$. Since energy differences are additive, we
obtain~\eqref{eq:reduction}. The claims on the values of~$v\vee w$ follow from the fact
that~$w=-1$ in~$\mathcal R$ and~$w=+1$ in~$\Gamma_M\setminus\mathcal R$.

Now, let~$\sigma$ be a minimizer of~$H_{\mathcal R}(\cdot\mid\eta)$ and
let~$v\in\mathcal F$. The configuration~$v\vee w$ coincides with~$\eta$
outside~$\mathcal R$, and so it is a competitor for~$\sigma$ in the box, which gives
that~$H(\sigma)-H(v\vee w)\le0$. From this and~\eqref{eq:reduction} we deduce that
$$ H(\sigma)-H(v)=\big[H(\sigma)-H(v\vee w)\big]+\big[H(v\vee w)-H(v)\big]\le0,$$
as desired.
\end{proof}

\begin{proposition}\label{cor:cylinder}
Assume that~$K(i,j)>0$ for every~$i\neq j$. Then, every minimizer
of~$H_{\Gamma_M}(\cdot\mid\eta)$ is equal to~$+1$ in~$\Gamma_M\setminus\mathcal R$, and it
is a minimizer of~$H_{\mathcal R}(\cdot\mid\eta)$.
\end{proposition}

\begin{proof}
Let $\sigma$ be any minimizer in
$\Gamma_M$ with exterior datum $\eta$, in the sense of~\Cref{rmk:gammainfinite}.
We claim that
\begin{equation}\label{CLS:aim}
\sigma\in\mathcal F.\end{equation}
For finite $M$, this follows from the fact that~$\{\sigma\ne u^\star\}\subseteq\Gamma_M$ and the finiteness of~$\Gamma_M$.

Hence, to prove~\eqref{CLS:aim}, we can suppose that~$M=+\infty$. For each integer~$N\ge h$, we set
\begin{equation}\label{eq2c.yliKA92s}
E_N:=\{i\in\mathcal G:\ |i_d|\le N,\ \sigma_i=-1\}
\end{equation}
and observe that
\begin{equation}\label{eq2c.yliKA92s:02}
E_N\subseteq E_{N+1}.\end{equation}

Let now~$v$ be equal to~$+1$ in $E_N$ and~$\sigma$ elsewhere. Since $E_N$ is finite,
$v$ is an admissible competitor. For $i\in E_N$ and $j\notin E_N$ we have
$\sigma_i\sigma_j-v_iv_j=-2\sigma_j$, while all other spin products are unchanged.

Thus, by minimality and~\eqref{eq:Hdiff},
\begin{equation}\label{eq:cylinder-flip-energy}
0\le\frac{H(v)-H(\sigma)}{2}=\sum_{i\in E_N}S_N(i),
\end{equation}
where
$$
S_N(i):=\sum_{\substack{j\notin E_N\\\sigma_j=-1}}K(i,j)
        -\sum_{\substack{j\notin E_N\\\sigma_j=+1}}K(i,j).$$

Let $e_1$ denote the first coordinate vector. We define~$q:=(2L+1)e_1$. We choose~$r\in\N\setminus\{0\}$ such that
$$\sum_{|z|>r}K(z)\le\frac{ K(q)}2.$$ Let also~$N>h+2r$. Then, if $i\in E_N$ satisfies~$h+r<|i_d|\le N-r$, it follows that
$$(i+q)_1=i_1+2L+1\ge L+1\qquad{\mbox{and}}\qquad |(i+q)_d|=|i_d|>h.$$
Hence $i+q$ lies outside $\mathcal G$ and $\sigma_{i+q}=\eta_{i+q}=+1$.

Additionally, if~$j\notin E_N$ and~$\sigma_j=-1$, 
then either~$j\in\mathcal B$, in which case~$|j_d-i_d|\ge |i_d|-h>r$, or~$j\in\mathcal G$ with~$|j_d|>N$, in which case
$|j_d-i_d|\ge |j_d|-|i_d|>N-(N-r)=r$. 

As a result,
\begin{equation}\label{QDSOP:01}
S_N(i)\le \sum_{|j-i|>r}K(i,j)-K(i,i+q)
         =\sum_{|z|>r}K(z)-K(q)\le-\frac{K(q)}{2}.\end{equation}

We now define
$$F_N:=\{i\in E_N:\ |i_d|\le h+r\ \text{or}\ N-r<|i_d|\le N\}$$
and we observe that
\begin{equation}\label{QDSOP:02}S_N(i)\le\sum_{z\in\Z^d}K(z).\end{equation} Each horizontal section of $\mathcal G$
contains $(2L+1)^{d-1}$ sites. The interval $|i_d|\le h+r$ contains
$2(h+r)+1$ integer heights, and $N-r<|i_d|\le N$ contains $2r$. Therefore,
\begin{equation}\label{SKD:SHp-2er4tvewo}\sharp F_N\le (2L+1)^{d-1}\big(2(h+r)+1+2r\big).\end{equation}
By~\eqref{eq:cylinder-flip-energy}, the sum of $S_N(i)$ over $E_N$ is nonnegative. Hence,
we use~\eqref{QDSOP:01} on~$E_N\setminus F_N$ and~\eqref{QDSOP:02} on~$F_N$ and we obtain that
\begin{align*}
0\le\sum_{i\in E_N}S_N(i)
&=\sum_{i\in E_N\setminus F_N}S_N(i)+\sum_{i\in F_N}S_N(i)\\
&\le -\frac{K(q)}{2}\big(\sharp E_N-\sharp F_N\big)
     +\sharp F_N\,\sum_{z\in\Z^d}K(z).
\end{align*}
Hence, solving for $\sharp E_N$ and using~\eqref{SKD:SHp-2er4tvewo},
we arrive at
\begin{equation}\label{eq2c.yliKA92s:o0oejfmew}
\sharp E_N\le \left(1+\frac{2}{K(q)}\sum_{z\in\Z^d}K(z)\right)(2L+1)^{d-1}(2h+4r+1).
\end{equation}
The right-hand side of~\eqref{eq2c.yliKA92s:o0oejfmew} is independent of $N$,
so, in light of~\eqref{eq2c.yliKA92s:02}, there exists~$N_0>h+2r$ such that~$\sharp E_N=\sharp E_{N_0}$ for all~$N\ge N_0$,
and therefore~$E_N=E_{N_0}$ for every $N\ge N_0$.

Now take any $i\in\mathcal G$ with $\sigma_i=-1$. Choosing $N\ge\max\{N_0,|i_d|\}$
in~\eqref{eq2c.yliKA92s}, we have that~$i\in E_N=E_{N_0}$.
Consequently, $E_{N_0}$ contains every negative site of $\sigma$ in $\mathcal G$. Moreover,
by~\eqref{eq:ustar}, we know that~$u^\star=+1$ in~$\mathcal G$.
Also, outside $\mathcal G$, both $\sigma$ and $u^\star$ equal $\eta$.
It thereby follows that~$\{\sigma\ne u^\star\}=E_{N_0}$.
This is a finite subset of~$\mathcal G$, hence~$\sigma\in\mathcal F$, as claimed in~\eqref{CLS:aim}.

Let now~$\psi_1$ be as in the proof
of~\Cref{lem:reduction}. Let also~$A:=\{\sigma<\psi_1\}$ and~$B:=\{\sigma>\psi_1\}$.
We have that~$A$ is a finite subset of~$\Gamma_M$
and~$\sigma\vee\psi_1\in\mathcal F$.

By~\eqref{eq:subdiffAB}, with~$v:=\sigma$ and~$\psi:=\psi_1$, and the minimality
of~$\psi_1$,
$$ H(\sigma\vee\psi_1)-H(\sigma)= H(\psi_1)-H(\sigma\wedge\psi_1)-4\,I(A,B)
\le-4\,I(A,B).$$
Since~$\sigma$ is a minimizer and~$\sigma\vee\psi_1\in\mathcal F$, the left hand side
is nonnegative, and thus~$I(A,B)=0$. Moreover, if~$j\in\Z^d\setminus\Gamma_M$ and~$j_d<-h$,
then~$\sigma_j=\eta_j=+1$ and~$\psi_1(j)=-1$, and so~$B\neq\varnothing$.
Since~$K>0$, this gives that~$A=\varnothing$, namely~$\sigma=+1$ in the sites
of~$\Gamma_M$ with~$i_d>h$. The same argument with~$\psi_2$ shows that~$\sigma=+1$ in the
sites of~$\Gamma_M$ with~$i_d<-h$.

As a result, $\sigma$ coincides with~$\eta$ outside~$\mathcal R$. Since the competitors
in the box belong to~$\mathcal F$, we conclude that~$\sigma$ is a minimizer
of~$H_{\mathcal R}(\cdot\mid\eta)$.
\end{proof}

\begin{proof}[Proof of~\Cref{th:fumihiko} when~$M>h$]
Part~\ref{it:fractional} does not involve~$M$, and so we focus on
part~\ref{it:rigidity}. By~\Cref{th:connected}, $u^\star$ is the unique minimizer
of~$H_{\mathcal R}(\cdot\mid\eta)$. Hence, it is a minimizer
of~$H_{\Gamma_M}(\cdot\mid\eta)$, by~\Cref{lem:reduction}. If in addition~$K(i,j)>0$
for every~$i\neq j$, every minimizer of~$H_{\Gamma_M}(\cdot\mid\eta)$ is a minimizer
of~$H_{\mathcal R}(\cdot\mid\eta)$, by~\Cref{cor:cylinder}, and so it coincides
with~$u^\star$.
\end{proof}

\section{Proof of~\Cref{lem:NNcalib}}\label{sec:NNapp}

In this appendix we prove~\Cref{lem:NNcalib}. In this setting, $K=K_{\rm nn}$ is the
nearest-neighbour kernel in~\eqref{eq:NNkernel}.
We first show how a vector field in the continuum induces a discrete vector field
in~$\mathcal R$. Then, we solve a Poisson equation with Neumann boundary conditions,
and we take the gradient of its solution.

For every~$l\in\Z^d$, we consider the closed unit cube centred at~$l$, namely
$$ C_l:=l+\Bigl[-\frac12,\frac12\Bigr]^d,$$
and we set
\begin{equation}\label{eq:Qdef}
Q:=\bigcup_{l\in\mathcal R}C_l=\Bigl[-L-\frac12,\,L+\frac12\Bigr]^{d-1}\times\Bigl[-h-\frac12,\,h+\frac12\Bigr].
\end{equation}
For every~$l\in\Z^d$,
$$ \partial C_l=\bigcup_{{j\in\Z^d}\atop{|j-l|=1}}\big(C_l\cap C_j\big),$$
where each face~$C_l\cap C_j$ is a~$(d-1)$-dimensional unit cube, on which the exterior
normal of~$C_l$ is
$$ \nu_{l,j}:=j-l.$$
Similarly,
$$ \partial Q=\bigcup_{{l\in\mathcal R,\,j\in\Z^d\setminus\mathcal R}\atop{|l-j|=1}}\big(C_l\cap C_j\big),$$
and the exterior normal of~$Q$ on~$C_l\cap C_j$ is~$\nu_{l,j}$.

\begin{proposition}\label{prop:cells}
Let~$\xi=(\xi_1,\dots,\xi_d)\in C^1(Q;\mathbb R^d)$. For every~$l$, $j\in\mathcal R$, let
\begin{equation}\label{eq:gflux}
g_\xi(l,j):=\begin{cases}
\displaystyle\int_{C_l\cap C_j}\xi\cdot\nu_{l,j}\,d\mathcal H^{d-1}, &{\mbox{ if }}|l-j|=1,\\
0, &{\mbox{ otherwise.}}
\end{cases}
\end{equation}
Then, the following hold:
\begin{enumerate}[label=\textup{(\roman*)}]
\item\label{it:cellsone} The function~$g_\xi$ is a discrete vector field in~$\mathcal R$.
\item\label{it:cellstwo} If~$|\xi_k|\le1$ in~$Q$ for every~$k\in\{1,\dots,d\}$, then~$|g_\xi|\le1$.
\item\label{it:cellsthree} If~$\xi\cdot\nu_{l,j}=u^\star_j$ on~$C_l\cap C_j$ for every~$l\in\mathcal R$
and~$j\in\Z^d\setminus\mathcal R$ with~$|l-j|=1$, then,
for every~$l\in\mathcal R$,
\begin{equation}\label{eq:celldiv}
D_{g_\xi}(l)+f(l)=\int_{C_l}\operatorname{div}\xi\,dx.
\end{equation}
\end{enumerate}
\end{proposition}

\begin{proof}
Since~$\nu_{j,l}=-\nu_{l,j}$, we deduce from~\eqref{eq:gflux} that~$g_\xi(j,l)=-g_\xi(l,j)$, which
gives~\ref{it:cellsone}. Moreover, if~$j=l\pm e_k$, then~$\nu_{l,j}=\pm e_k$, and therefore
$$ |g_\xi(l,j)|\le\int_{C_l\cap C_j}|\xi_k|\,d\mathcal H^{d-1}\le\mathcal H^{d-1}\big(C_l\cap C_j\big)=1,$$
which proves~\ref{it:cellstwo}.

Now we prove~\ref{it:cellsthree}. Let~$l\in\mathcal R$, by the divergence theorem applied in~$C_l$,
by~\eqref{eq:Dg} and by~\eqref{eq:PN}, we have that
\begin{align*}
\int_{C_l}\operatorname{div}\xi\,dx&=\sum_{{j\in\Z^d}\atop{|j-l|=1}}\int_{C_l\cap C_j}\xi\cdot\nu_{l,j}\,d\mathcal H^{d-1}
=\sum_{{j\in\mathcal R}\atop{|j-l|=1}}g_\xi(l,j)+\sum_{{j\in\Z^d\setminus\mathcal R}\atop{|j-l|=1}}u^\star_j\\
&=D_{g_\xi}(l)+f(l),
\end{align*}
which is~\eqref{eq:celldiv}.
\end{proof}

Let now~$l\in\mathcal R$ and~$j\in\Z^d\setminus\mathcal R$ with~$|l-j|=1$. If~$j=l\pm e_d$,
then~$j\in\mathcal U$, and so~$u^\star_j=+1$. If instead~$j=l\pm e_k$ with~$k\le d-1$,
then~$j\in\mathcal B$, and so~$u^\star_j=-1$. Therefore, in view
of~\Cref{prop:cells}, we look for~$\xi$ of the form~$\xi=-\nabla\Psi$, where~$\Psi$
solves, in the notation of~\eqref{RALGLE},
\begin{equation}\label{eq:Neumann}
\begin{cases}
-\Delta\Psi=\langle f\rangle & {\mbox{ in }}Q,\\
\partial_\nu\Psi=-1 & {\mbox{ on }}\partial Q\cap\big\{|x_d|=h+\frac12\big\},\\
\partial_\nu\Psi=1 & {\mbox{ on }}\partial Q\cap\big\{|x_k|=L+\frac12\big\}{\mbox{ for every }}k\le d-1,\\
|\partial_k\Psi|\le1 & {\mbox{ in }}Q{\mbox{ for every }}k\in\{1,\dots,d\},
\end{cases}
\end{equation}
and~$\nu$ is the exterior normal of~$Q$. We observe that, if~$\Psi$ satisfies the first three conditions
in~\eqref{eq:Neumann}, then, by the divergence theorem in~$Q$,
\begin{align*}
\langle f\rangle\,\mathcal H^d(Q)&=-\int_Q\Delta\Psi\,dx=-\int_{\partial Q}\partial_\nu\Psi\,d\mathcal H^{d-1}\\
&=\mathcal H^{d-1}\big(\partial Q\cap\{|x_d|=h+\tfrac12\}\big)
-\sum_{k=1}^{d-1}\mathcal H^{d-1}\big(\partial Q\cap\{|x_k|=L+\tfrac12\}\big)\\
&=2(2L+1)^{d-1}-2(d-1)(2L+1)^{d-2}(2h+1).
\end{align*}
Since~$\mathcal H^d(Q)=(2L+1)^{d-1}(2h+1)$, this is~\eqref{eq:F0value}. Hence, the
first three conditions in~\eqref{eq:Neumann} are compatible.

We look for a solution of~\eqref{eq:Neumann} of the form
$$\Psi(x)=\psi_1(x_1)+\dots+\psi_d(x_d).$$
In this case, the first equation in~\eqref{eq:Neumann} reads
$$ -\sum_{k=1}^d\psi_k''(x_k)=\langle f\rangle .$$
Since the~$k$-th term depends only on~$x_k$, each~$\psi_k''$ is constant. Moreover, the
boundary conditions in~\eqref{eq:Neumann} read
$$ \psi_k'\Bigl(\pm\Bigl(L+\frac12\Bigr)\Bigr)=\pm1\quad{\mbox{for every }}k\le d-1
\qquad{\mbox{and}}\qquad
\psi_d'\Bigl(\pm\Bigl(h+\frac12\Bigr)\Bigr)=\mp1.$$
Since each~$\psi_k'$ is affine, we obtain that
$$ \psi_k'(t)=\frac{2t}{2L+1}\quad{\mbox{for every }}k\le d-1
\qquad{\mbox{and}}\qquad
\psi_d'(t)=-\frac{2t}{2h+1}.$$
As a result, up to an additive constant,
\begin{equation}\label{eq:PsiNN}
\Psi(x)=\frac{x_1^2+\dots+x_{d-1}^2}{2L+1}-\frac{x_d^2}{2h+1}.
\end{equation}
The function in~\eqref{eq:PsiNN} satisfies all the conditions in~\eqref{eq:Neumann}.
Indeed, by~\eqref{eq:F0value},
$$ -\Delta\Psi=-\frac{2(d-1)}{2L+1}+\frac{2}{2h+1}=\langle f\rangle .$$
Moreover, for every~$x\in Q$,
$$ |\partial_k\Psi(x)|=\frac{2|x_k|}{2L+1}\le1\quad{\mbox{for every }}k\le d-1
\qquad{\mbox{and}}\qquad
|\partial_d\Psi(x)|=\frac{2|x_d|}{2h+1}\le1.$$

We now go back to the lattice. The discrete vector field induced by~$-\nabla\Psi$ has
an explicit form, which is the content of the next result.

\begin{proposition}\label{prop:gNN}
Let~$\Psi$ be as in~\eqref{eq:PsiNN} and, for every~$l$, $j\in\mathcal R$, let
\begin{equation}\label{eq:gNN}
g(l,j):=\begin{cases}
\Psi(l)-\Psi(j), &{\mbox{ if }}|l-j|=1,\\
0, &{\mbox{ otherwise.}}
\end{cases}
\end{equation}
Then, $g$ is a discrete vector field in~$\mathcal R$ which solves~\eqref{eq:NNdiv}, and
$$ |g(l,j)|<1\qquad{\mbox{for every }}l,\,j\in\mathcal R.$$
\end{proposition}

\begin{proof}
Let~$\xi:=-\nabla\Psi$, and let~$g_\xi$ be as in~\eqref{eq:gflux}. Since~$\Psi$
satisfies~\eqref{eq:Neumann}, by~\Cref{prop:cells} the function~$g_\xi$ is a discrete
vector field in~$\mathcal R$ and, for every~$l\in\mathcal R$,
$$ D_{g_\xi}(l)+f(l)=\int_{C_l}\operatorname{div}\xi\,dx=-\int_{C_l}\Delta\Psi\,dx=\langle f\rangle.$$
Hence, it is enough to show that
\begin{equation}\label{eq:gtilde}
g_\xi=g
\end{equation}
and that~$|g|<1$.

We first prove~\eqref{eq:gtilde}. Since~$g_\xi$ and~$g$ are antisymmetric and vanish on
the pairs which are not nearest neighbours, it is enough to consider the
pairs~$(l,l+e_k)$ with~$l$, $l+e_k\in\mathcal R$ and~$k\in\{1,\dots,d\}$.
By~\eqref{eq:PsiNN},
$$ \partial_k\Psi(x)=\frac{2x_k}{2L+1}\quad{\mbox{for every }}k\le d-1
\qquad{\mbox{and}}\qquad
\partial_d\Psi(x)=-\frac{2x_d}{2h+1}.$$
In particular, $\partial_k\Psi$ is constant on~$C_l\cap C_{l+e_k}\subseteq\{x_k=l_k+\frac12\}$,
and it is affine on the segment from~$l$ to~$l+e_k$. Therefore,
by~\eqref{eq:gflux} and~\eqref{eq:gNN},
\begin{equation}\label{eq:gNNexpl}\begin{split}&
g_\xi(l,l+e_k)=-\partial_k\Psi\Bigl(l+\frac{e_k}2\Bigr)
=-\int_0^1\partial_k\Psi(l+te_k)\,dt\\&\quad=\Psi(l)-\Psi(l+e_k)=g(l,l+e_k),\end{split}
\end{equation}
which proves~\eqref{eq:gtilde}.

We now prove that~$|g|<1$. Let~$l$, $l+e_k\in\mathcal R$. The face~$C_l\cap C_{l+e_k}$
lies in the interior of~$Q$, namely~$|l_k+\frac12|\le L-\frac12$ if~$k\le d-1$
and~$|l_d+\frac12|\le h-\frac12$. Hence, by~\eqref{eq:PsiNN} and~\eqref{eq:gNNexpl},
\begin{equation*} |g(l,l+e_k)|=\Bigl|\partial_k\Psi\Bigl(l+\frac{e_k}2\Bigr)\Bigr|
\le\max\Bigl\{\frac{2L-1}{2L+1},\frac{2h-1}{2h+1}\Bigr\}<1.\qedhere\end{equation*}
\end{proof}

\Cref{lem:NNcalib} now follows from~\Cref{prop:gNN}.

\section{The case of the fractional kernel}\label{sec:fracfield}

We recall that~$\Lambda$ and $S$ are defined in~\eqref{dfeseeselamvda876543}
and that~$\underline f$ is defined in~\eqref{dfeseeselamvda8765432}.

We also observe that, for every~$m$, $k\in\mathbb N$
with~$m$, $k\ge1$,
\begin{equation}\label{eq:spherecount}
\begin{split}
\sharp\big\{z\in\Z^m\ \text{such that}\ \|z\|_\infty=k\big\}
&=(2k+1)^m-(2k-1)^m\\
&\le2m\,(2k+1)^{m-1}\\
&\le2m\,3^{m-1}k^{m-1}.
\end{split}
\end{equation}
Moreover, if $m \in \{d-1,d\}$, we have that $2m\,3^{m-1}k^{m-1} \le C_d k^{m-1}$, where $C_d:=2d\,3^{d}$. In particular,
\begin{equation}\label{eq:Cd}
\sharp\big\{z\in\Z^m\ \text{such that}\ \|z\|_\infty=k\big\} \le C_d k^{m-1}.
\end{equation}

\begin{lemma}\label{lem:fracfield}
Let~$d\ge2$, let~$h,L\in\mathbb N$, and let~$K=K_s$ be as in~\eqref{eq:frachyp}. Let
\begin{equation}\label{eq:s1def}
c_d:=\min\left\{\frac12,\ \frac{e^{-1}d^{-(d+2)/2}}{16\,C_d}\right\}
\qquad{\mbox{and}}\qquad
s_1:=\frac{c_d}{1+\log(2h+2)}.
\end{equation}
Then, for every~$s\in(0,s_1]$ and every~$l\in\mathcal R$,
\begin{equation}\label{eq:gpositive}
\underline f(l)\ge4C_d\big(1+\log(2h+2)\big)>0.
\end{equation}
\end{lemma}

\begin{proof}
Let~$l\in\mathcal R$ and define, for every~$k\ge1$,
\[
S_k(l):=\big\{j\in\Z^d\ \text{such that}\ \|j-l\|_\infty=k\big\}.
\]
\noindent\emph{Step 1: Upper bound on~$I(\{l\},S)$.} Since~$\|z\|_\infty\le|z|$ for
every~$z\in\Z^d$, we have that
\begin{equation}\label{eq:Kupper}
K_s(l,j)\le k^{-d-2s}\qquad{\mbox{for every }}j\in S_k(l).
\end{equation}
Now, let $j \in S$, then ~$|j_d|\le h+1$.  Since $l\in \mathcal R$, we also have that~$|l_d|\le h$. Therefore,
\[
|j_d-l_d| \le 2h+1.
\]
As a consequence, if~$k\ge2h+2$ and~$j\in S\cap S_k(l)$,
then
$$ \|j'-l'\|_\infty=k\qquad{\mbox{and}}\qquad j_d\in\{-h-1,\dots,h+1\}.$$
Thus, by applying~\eqref{eq:Cd}, we have that
\begin{equation}\label{eq:Scount}
\sharp(S\cap S_k(l))\le
\begin{cases}
C_d\,k^{d-1}, & {\mbox{ if }} 1\le k\le2h+1,\\
(2h+3)\,C_d\,k^{d-2}, & {\mbox{ if }} k\ge2h+2.
\end{cases}
\end{equation}
Moreover, since~$l\notin S$, every~$j\in S$ belongs to~$S_k(l)$ with~$k=\|j-l\|_\infty\ge1$.
Hence, the sets~$S\cap S_k(l)$
form a partition of~$S$, and
\[
I(\{l\},S)=\sum_{j\in S}K_s(l,j)=\sum_{k=1}^{+\infty}\,\sum_{j\in S\cap S_k(l)}K_s(l,j).
\]

Now, using~\eqref{eq:Kupper} and~\eqref{eq:Scount}, we have
\begin{eqnarray*}
I(\{l\},S)&\le&\ C_d\sum_{k=1}^{2h+1}\frac1{k^{1+2s}}+(2h+3)\,C_d\sum_{k=2h+2}^{+\infty}\frac1{k^{2+2s}} \\
&\le& C_d\sum_{k=1}^{2h+1}\frac1{k}+(2h+3)\,C_d\sum_{k=2h+2}^{+\infty}\frac1{k^{2}}.
\end{eqnarray*}

We also observe that
$$ \sum_{k=1}^{2h+1}\frac1k\le1+\log(2h+1)
\qquad{\mbox{and}}\qquad
\sum_{k=2h+2}^{+\infty}\frac1{k^2}\le\frac1{2h+1}.$$
Finally,
\begin{equation}\label{eq:fracS}
I(\{l\},S) \le C_d\big(1+\log(2h+1)\big)+\frac{2h+3}{2h+1}\,C_d  \le 4C_d\big(1+\log(2h+2)\big),
\end{equation}
where we used the fact that~$\frac{2h+3}{2h+1}\le3$. 

\noindent\emph{Step 2: Lower bound on~$I(\{l\},\Lambda)$.} Since~$|z|\le\sqrt d\,\|z\|_\infty$ for
every~$z\in\Z^d$, and~$s<1$, we have that
\begin{equation}\label{eq:Klower}
K_s(l,j)\ge d^{-(d+2)/2}\,k^{-d-2s}\qquad{\mbox{for every }}j\in S_k(l).
\end{equation}
For every~$k\ge2h+2$, we set
$$ T_k:=\big\{j\in\Z^d\ \text{such that}\ j_d=l_d+k\ \text{and}\
\|j'-l'\|_\infty\le k\big\}.$$
We observe that~$T_k\subseteq S_k(l)$ and~$\sharp T_k=(2k+1)^{d-1}\ge k^{d-1}$. Moreover,
if~$j\in T_k$, then~$j_d\ge k-h\ge h+2$, and so~$T_k\subseteq\Lambda$. Since the
sets~$T_k$ are pairwise disjoint, we deduce from~\eqref{eq:Klower} that
\begin{equation}\label{eq:fracLambda}\begin{split}
I(\{l\},\Lambda)&\ge\sum_{k=2h+2}^{+\infty}\,\sum_{j\in T_k}K_s(l,j)
\ge d^{-(d+2)/2}\sum_{k=2h+2}^{+\infty}\frac{1}{k^{1+2s}}\\
&\ge d^{-(d+2)/2}\int_{2h+2}^{+\infty}\frac{dr}{r^{1+2s}}
=\frac{d^{-(d+2)/2}}{2s}\,(2h+2)^{-2s}.
\end{split}\end{equation}

\noindent\emph{Step 3: Conclusion.} Let~$s\in(0,s_1]$ and let~$c_d$ be the quantity defined in~\eqref{eq:s1def}. We have
that
$$ 2s\log(2h+2)\le2s_1\big(1+\log(2h+2)\big)=2c_d\le1,$$
which gives that
\begin{equation}\label{eq:expbound}
(2h+2)^{-2s}\ge e^{-1}.
\end{equation}
Furthermore, by~\eqref{eq:s1def},
\begin{equation}\label{eq:oneovers}
\frac1s\ge\frac{1+\log(2h+2)}{c_d}\ge16\,e\,d^{(d+2)/2}\,C_d\big(1+\log(2h+2)\big).
\end{equation}
Plugging~\eqref{eq:expbound} and~\eqref{eq:oneovers} into~\eqref{eq:fracLambda}, we obtain that
\begin{equation}\label{intestimate}
  I(\{l\},\Lambda)\ge\frac{e^{-1}d^{-(d+2)/2}}{2s}\ge8C_d\big(1+\log(2h+2)\big).
\end{equation}

Thus, by~\eqref{intestimate} and~\eqref{eq:fracS}, we conclude that
$$ \underline f(l)\ge4C_d\big(1+\log(2h+2)\big),$$
which is~\eqref{eq:gpositive}.
\end{proof}

\section*{Acknowledgments}
SD and EV are members of the Australian Mathematical Society (AustMS).

This work has been supported by 
the Australian Future Fellowship FT230100333
and the Australian Laureate Fellowship FL190100081.

\vfill
\end{document}